\documentclass[a4paper]{amsart}
\usepackage{color}
\usepackage{amssymb,amsmath,amsthm,amstext,amsfonts,amscd}
\usepackage{psfrag}
\usepackage{url}
\usepackage{amsfonts}
\usepackage{mathrsfs}
\usepackage{graphicx}
\usepackage{enumerate}
\usepackage{fancyhdr}
\usepackage{enumitem}
\usepackage[colorlinks,linkcolor=black,anchorcolor=black,citecolor=black,hyperindex=true,CJKbookmarks=true]{hyperref}

\usepackage{epsfig}
\usepackage[T1]{fontenc}

\makeatletter
\@addtoreset{equation}{section}
\@addtoreset{figure}{section}
\@addtoreset{table}{section}

\renewcommand\thetable{\thesection.\@arabic\c@table}
\renewcommand\thefigure{\thesection.\@arabic\c@figure}
\makeatother

\usepackage{tikz}
\usetikzlibrary{shapes, arrows, positioning, shadows}

\usepackage{booktabs}

\makeatletter
\@namedef{subjclassname@2020}{2020 Mathematics Subject Classification}
\makeatother

\newtheorem{lemma}{Lemma}[section]
\newtheorem{theorem}[lemma]{Theorem}
\newtheorem{proposition}[lemma]{Proposition}
\newtheorem{corollary}[lemma]{Corollary}
\newtheorem{question}{Question}

\newtheorem{maintheorem}{Theorem}

\newtheorem{maincorollary}{Corollary}

\theoremstyle{remark}
\newtheorem{definition}[lemma]{Definition}
\newtheorem{remark}[lemma]{Remark}

\makeatletter
\def\lst@lettertrue{\let\lst@ifletter\iffalse}
\makeatother

\title{Variations of topological theory and ergodic theory via gap function in non-uniform specification}

\author{Wanshan Lin, Xueting Tian and Chenwei Yu}

\address{Wanshan Lin, School of Mathematical Sciences,  Fudan University\\Shanghai 200433, People's Republic of China}
\email{wanshanlin@fudan.edu.cn}

\address{Xueting Tian, School of Mathematical Sciences,  Fudan University\\Shanghai 200433, People's Republic of China}
\email{xuetingtian@fudan.edu.cn}

\address{Chenwei Yu, School of Mathematical Sciences,  Fudan University\\Shanghai 200433, People's Republic of China}
\email{23110180052@m.fudan.edu.cn}
\begin{document}

\renewcommand{\thepage}{\arabic{page}}
\setcounter{page}{1}

\footnotetext{All authors are co-first authors of the article.}

\begin{abstract}
	In our previous work \cite{Lin-Tian-Yu-2024}, we studied qualitative differences between specification and non-uniform specification.
	In this paper, we investigate quantitative variations depending on the gap function.
	We first obtain lower bounds for the Bowen topological entropy of irregular sets in terms of the lower linear growth of gap function.
	In contrast, we prove that every non-empty over-saturated set has full packing topological entropy under non-uniform specification.
	We also establish a quantitative lower bound for the Bowen topological entropy of transitive points under non-uniform specification, and a lower bound for the exponential growth of periodic orbits under its periodic version.
	Besides, we construct symbolic systems with non-uniform specification with a given gap growth rate which contain an arbitrary subshift.
	These systems show that certain bounds concerning irregular sets and periodic orbits are optimal.
	They also show that positive linear gap growth may destroy full Bowen entropy of transitive points, the conditional variational principle, the intermediate entropy and pressure properties, and the genericity of continuous functions whose unique maximizing measure has zero entropy.
\end{abstract}

\maketitle

\thispagestyle{empty}

\tableofcontents

\section{Introduction}
Throughout this paper, $(X,f)$ denotes a topological dynamical system, where $(X,d)$ is a compact metric space and $f:X\rightarrow X$ is continuous.
We assume that $X$ has at least two points.
Let $C(X)$ be the Banach space of real-valued continuous functions on $X$ with the supremum norm $\|\varphi\|=\sup_{x\in X}|\varphi(x)|$.
We write $\mathbb{Z}$, $\mathbb{N}$, and $\mathbb{N}_0$ for the integers, the positive integers, and the non-negative integers, respectively.

The specification-like properties, which were first introduced by Bowen \cite{Bowen-1971} play important roles in the study of the complexity of dynamical systems, including multifractal analysis, classification of recurrent points, growth of periodic points, intermediate entropy and pressure, ergodic optimization and so on. In this paper, we mainly consider the non-uniform specification, a variation of the Bowen's specification property.

\begin{definition}
	We say that a dynamical system $(X,f)$ satisfies the \emph{non-uniform specification property} with \emph{gap function} $M(n,\varepsilon):\mathbb{N}\times(0,+\infty)\rightarrow\mathbb{N}$ if the following conditions hold:
	\begin{itemize}
		\item $M(n,\varepsilon)$ is non-decreasing with $n$ and non-increasing with $\varepsilon$;
		\item for any positive integer $k\geq2$, any $x_1,\cdots,x_k\in X$ and any non-negative integers $a_1,b_1,\cdots,a_k,b_k$ with
		\begin{equation*}
			a_1\leq b_1<\cdots<a_k\leq b_k
		\end{equation*}
		and
		\begin{equation*}
			a_{i+1}-b_i\geq M(b_i-a_i+1,\varepsilon)\quad\text{for}\quad 1\leq i\leq k-1,
		\end{equation*}
		there exists a point $z\in X$ such that
		\begin{equation*}
			d(f^{n-a_i}x_i,f^nz)\leq\varepsilon,\quad\text{for}\quad a_i\leq n\leq b_i,1\leq i\leq k.
		\end{equation*}
	\end{itemize}
	If further, for any positive integer $p$ with $p\geq b_k-a_1+M(b_k-a_k+1,\varepsilon)$, the tracing point $z$ above can be chosen as a periodic point with $f^pz=z$, then $(X,f)$ is said to satisfy the \emph{non-uniform periodic specification property}.
\end{definition}

The definition of the non-uniform specification property is modified from the definitions of the weak specification property \cite[Definition 14]{Kwietniak-Lacka-Oprocha-2016} and non-uniform specification for subshifts \cite[Definition 2.14]{Pavlov-2016}.
In \cite{Kwietniak-Lacka-Oprocha-2016}, the weak specification property is also called the forward weak specification property, while its dual version is called the backward weak specification property.
In fact, the backward weak specification property coincides with the non-uniform specification property under the condition that
\begin{equation}\label{eq-intro-nonuniform-spec}
	\sup_{\varepsilon>0}\lim_{n\rightarrow\infty}\frac{M(n,\varepsilon)}{n}=0.
\end{equation}

The weak specification property was first utilized by Marcus in \cite{Marcus-1980} (though unnamed) to investigate ergodic toral diffeomorphisms.
Subsequently, Dateyama called it the almost weak specification property in \cite{Dateyama-1983}.
Quas and Soo \cite{Quas-Soo-2016} used this property together with asymptotic entropy expansiveness and the small boundary property to study universality.
Later, Burguet \cite{Burguet-2020} demonstrated that these additional conditions can be removed.

Let $\mathcal{M}(X)$, $\mathcal{M}_f(X)$ and $\mathcal{M}_f^{erg}(X)$ denote the space of Borel probability measures, $f$-invariant Borel probability measures and $f$-ergodic Borel probability measures, respectively. Let $h_{top}(f)$ denote the topological entropy of $(X,f)$ and $h_{\mu}(f)$ denote the measure-theoretic entropy of $\mu\in\mathcal{M}_f(X)$.
The variational principle gives
\begin{equation*}
	h_{top}(f)=\sup\{h_{\mu}(f):\mu\in\mathcal{M}_f(X)\},
\end{equation*}
and hence $h_{\mu}(f)\leq h_{top}(f)$ for all $\mu\in\mathcal{M}_f(X)$.
An ergodic measure $\mu\in\mathcal{M}_f^{erg}(X)$ is called a \emph{measure of maximal entropy} if $h_{\mu}(f)=h_{top}(f)$.
A dynamical system is called \emph{intrinsically ergodic} if there is a unique measure of maximal entropy.
In \cite{Pavlov-2016}, Pavlov investigated the variations of intrinsic ergodicity for subshifts under two weak versions of specification, where the backward weak specification property for subshifts is named as non-uniform specification.
It was shown in \cite{Pavlov-2016} by giving an example that we cannot guarantee intrinsic ergodicity when $\inf_{\varepsilon>0}\liminf_{n\rightarrow\infty}(M(n,\varepsilon)/\log n)>0$.
In \cite{Pavlov-2019}, Pavlov further introduced the notion of controlled specification with gap function $f(n)$.
It was proved in \cite{Pavlov-2019} that $\liminf_{n\rightarrow\infty}(f(n)/\log n)=0$ is the critical condition to guarantee the intrinsic ergodicity for subshifts with controlled specification.

For each $x\in X$, let $\delta_x$ denote the Dirac measure at $x$.
For $x\in X$ and $n\in\mathbb{N}$, the \emph{empirical measure} associated with the orbit segment $\{x,fx,\cdots,f^{n-1}x\}$ is defined as
\begin{equation*}
	\delta_x^n=\frac 1n\sum_{k=0}^{n-1}\delta_{f^kx}.
\end{equation*}
Let $V_f(x)$ be the set of accumulation points of $\{\delta_x^n:n\geq 1\}$.
For any subset $K\subset\mathcal{M}_f(X)$, define the \emph{saturated set} $G_K$ and \emph{over-saturated set} $G^K$ by
\begin{equation*}
	G_K=\{x\in X:V_f(x)=K\}\quad\text{and}\quad G^K=\{x\in X:V_f(x)\supset K\},
\end{equation*}
respectively.
By \cite[Proposition 3.8]{DGS1976}, every $V_f(x)$ is a non-empty compact connected subset of $\mathcal{M}_f(X)$.
Consequently, $G_K\neq\varnothing$ requires $K\subset\mathcal{M}_f(X)$ to be a non-empty compact connected set.
For simplicity, let $G_{\mu}=G_{\{\mu\}}$ and $G^{\mu}=G^{\{\mu\}}$.
Each point in $G_{\mu}$ is called a \emph{generic point} of $\mu$.

For any $n\in\mathbb{N}$, let $\mathrm{Per}_n(f):=\{x\in X:f^nx=x\}$ denote the set of periodic points of period $n$.
For any $x\in\mathrm{Per}_n(f)$, the atomic measure supported on $\mathrm{Orb}(x,f)$ is $\delta_x^n$, called a \emph{CO-measure} or \emph{periodic measure}.
The set of all periodic measures is denoted by $\mathcal{M}_f^{co}(X)$.
It is clear that
\begin{equation*}
	\mathcal{M}_f^{co}(X)\subset\mathcal{M}_f^{erg}(X)\subset\mathcal{M}_f(X).
\end{equation*}

In \cite{Lin-Tian-Yu-2024}, Lin, Tian and Yu investigated how the asymptotic behavior of gap functions in systems satisfying the non-uniform specification property, and its periodic version, influences certain ergodic properties, including:
\begin{itemize}
	\item [(R1)] the density of $\mathcal{M}_f^{co}(X)$ in $\mathcal{M}_f(X)$;
	\item [(R2)] the density of $G_K$ in $X$ for every non-empty compact connected set $K\subset\mathcal{M}_f(X)$;
	\item [(R3)] the residuality of $\mathcal{M}_f^{erg}(X)$ in $\mathcal{M}_f(X)$;
	\item [(R4)] the connectedness of $\mathcal{M}_f^{erg}(X)$;
	\item [(R5)] the entropy-dense property of $(X,f)$.
\end{itemize}
Consider a dynamical system $(X,f)$ with non-uniform specification (resp. non-uniform periodic specification).
Let $\mathscr{G}(X,f)$ be the collection of all gap functions for which $(X,f)$ satisfies the non-uniform specification property (resp. the non-uniform periodic specification property).
Define
\begin{equation}\label{1.2}
	\tau_*(X,f)=\inf_{M\in\mathscr{G}(X,f)}\sup_{\varepsilon>0}\liminf_{n\rightarrow\infty}\frac{M(n,\varepsilon)}{n}.
\end{equation}
We say that $(X,f)$ satisfies the \emph{non-uniform $\tau$-specification property} (resp. the \emph{non-uniform periodic $\tau$-specification property}) if $\tau_*(X,f)=\tau$ and this infimum is attained by an admissible gap function.
Equivalently, there is an $M\in\mathscr{G}(X,f)$ with $\tau(M)=\tau$, and every $M'\in\mathscr{G}(X,f)$ satisfies $\tau(M')\geq\tau$, where
\begin{equation}
	\tau(M):=\sup_{\varepsilon>0}\liminf_{n\rightarrow\infty}\frac{M(n,\varepsilon)}{n}.
\end{equation}
Once the gap function $M$ is fixed, we write $\tau$ for $\tau(M)$ below.
The following proposition can be easily checked by the definition.

\begin{proposition}\label{intro-prop:conjugacy-invariant}
	Suppose that $(X,f)$ satisfies the non-uniform specification property and $(X',f')$ is topologically conjugate to $(X,f)$.
	Then the following hold:
	\begin{enumerate}
		\item $(X',f')$ satisfies the non-uniform specification property;
		\item for any admissible gap function $M\in\mathscr{G}(X,f)$, then there exists an admissible gap function $M'\in\mathscr{G}(X',f')$ with $\tau(M')=\tau(M)$;
		\item $\tau_*(X,f)=\tau_*(X',f')$.
	\end{enumerate}
\end{proposition}

It was proved in \cite{Lin-Tian-Yu-2024} that (R3)--(R5) hold for dynamical systems with non-uniform $0$-specification, while (R1)--(R2) hold when the non-uniform periodic $0$-specification property is assumed. The same paper gives explicit examples demonstrating the sharpness of the zero-growth condition: there exist systems satisfying the non-uniform periodic specification property with a gap function $M$ such that $\tau(M)>0$, for which all properties (R1)--(R5) fail.
These results demonstrate that the asymptotic behavior of gap function can change the measure and orbit structure of the system under the non-uniform specification property.
The purpose of the present paper is to determine quantitative effects of the gap growth.
A natural question arises.

\begin{question}\label{Que1}
	How does the asymptotic growth of the gap function quantitatively affect entropy and ergodic properties under non-uniform specification?
\end{question}

%Our goal of this paper is to deepen the understanding of this question by deriving a series of results related to it, especially in quantitative views.
Table \ref{tab: variations} and \ref{tab: estimates} illustrate some variations in topological theory and ergodic theory under non-uniform specification, including results from \cite{Lin-Tian-Yu-2024} and this paper.

\begin{table}[h]
    \centering
    \small
    \begin{tabular}{p{2.9cm} p{4.5cm} p{7.4cm}}
        \toprule
        & (Periodic) specification & (Periodic) non-uniform $\tau$-specification: Examples ($\tau>0$) \\
        \midrule
        Density of periodic measures & $\overline{\mathcal{M}_f^{co}(X)}=\mathcal{M}_f(X)$ \cite{Sigmund-1970} & Examples with $\overline{\mathcal{M}_f^{co}(X)}\subsetneq\mathcal{M}_f(X)$ \cite[Theorem A]{Lin-Tian-Yu-2024} \\
        \midrule
        Existence and density of saturated sets & $G_K$ is dense in $X$ for every non-empty connected compact subset $K\subset\mathcal{M}_f(X)$ \cite{Sigmund-1974} & Examples with $G_K=\varnothing$ for some non-empty connected compact set $K\subset\mathcal{M}_f(X)$ \cite[Theorem B]{Lin-Tian-Yu-2024} \\
        \midrule
        Topological entropy of non-empty irregular sets & $h_{top}^B(f,\mathrm{IR}(f))=h_{top}(f)$ \cite{Thompson-2010} & Examples with $h_{top}^B(f,\mathrm{IR}(f))<h_{top}(f)$ (Theorem \ref{thm-Bowen-upper-1}) \\
        \midrule
        Topological entropy of transitive points & $h_{top}^B(f,\mathrm{Trans}(f))=h_{top}(f)$ \cite{Tian-2016} & Examples with $h_{top}^B(f,\mathrm{Trans}(f))<h_{top}(f)$ (Theorem \ref{thm-trans-1}) \\
        \midrule
        Growth of periodic orbits (with positive expansiveness) & $\lim\limits_{n\to\infty}\dfrac{1}{n}\log|\mathrm{Per}_n(f)|=h_{top}(f)$ \cite{Bowen-1974} & Examples with $\limsup\limits_{n\to\infty}\dfrac1n\log|\mathrm{Per}_n(f)|<h_{top}(f)$ (Theorem \ref{thm-growth-2}) \\
        \bottomrule
		& & \\
    \end{tabular}
	\caption{Some differences between specification and non-uniform specification.}
	\label{tab: variations}
\end{table}

\begin{table}[h]
    \centering
    \small
    \begin{tabular}{p{4.8cm} p{10.4cm}}
        \toprule
         & Non-uniform (periodic) $\tau$-specification: Quantitative estimates \\
        \midrule
        Topological entropy of non-empty irregular sets & $h_{top}^B(f,\mathrm{IR}(f))\geq \Big(\dfrac{1}{1+\tau}-\dfrac{\tau^2}{(1+\tau)^2}\Big)h_{top}(f)$ ($\tau<1$, Theorem \ref{thm-Bowen-lower-1}); $h_{top}^P(f,\mathrm{CI}(f))=h_{top}(f)$ (Theorem \ref{thm-packing}) \\
        \midrule
        Topological entropy of transitive points & $h_{top}^B(f,\mathrm{Trans}(f))\geq\dfrac{1}{(1+\tau)^2} h_{top}(f)$ (Theorem \ref{thm-trans-2}) \\
        \midrule
        Growth of periodic orbits (with positive expansiveness) & $\dfrac{1}{1+\tau}h_{top}(f)\leq\liminf\limits_{n\to\infty}\dfrac{1}{n}\log|\mathrm{Per}_n(f)|\leq\limsup\limits_{n\to\infty}\dfrac{1}{n}\log|\mathrm{Per}_n(f)|\leq h_{top}(f)$ (Theorem \ref{thm-growth-1}) \\
        \bottomrule
		& \\
    \end{tabular}
	\caption{Some quantitative estimates under non-uniform specification.}
	\label{tab: estimates}
\end{table}

\subsection{Multifractal analysis}

For a $\varphi\in C(X)$, define the {\it $\varphi$-irregular set} $I_\varphi(f)$ by
\begin{equation*}
	I_{\varphi}(f)=\bigg\{x\in X:\lim_{n\rightarrow\infty}\frac1n\sum_{k=0}^{n-1}\varphi(f^kx)\;\text{does not exist}\bigg\}.
\end{equation*}
Let $\mathfrak{C}(f)=\{\varphi\in C(X):I_{\varphi}(f)\neq\varnothing\}$.
The \emph{irregular set} $\mathrm{IR}(f)$ and \emph{completely-irregular set} $\mathrm{CI}(f)$ of $(X,f)$ are respectively defined by
\begin{equation*}
	\mathrm{IR}(f)=\bigcup_{\varphi\in C(X)} I_{\varphi}(f)\quad\text{and}\quad\mathrm{CI}(f)=\bigcap_{\varphi\in\mathfrak{C}(f)}I_{\varphi}(f).
\end{equation*}

Define
\begin{equation*}
	L_{\varphi}=\bigg[\inf_{\mu\in\mathcal{M}_f(X)}\int_X\varphi\mathrm{d}\mu,\sup_{\mu\in\mathcal{M}_f(X)}\int_X\varphi\mathrm{d}\mu\bigg].
\end{equation*}
For each $a\in L_{\varphi}$, define the \emph{$\varphi$-level set} by
\begin{equation*}
	R_{\varphi}(a)=\bigg\{x\in X:\lim_{n\rightarrow\infty}\frac{1}{n}\sum_{i=0}^{n-1}\varphi(f^ix)=a\bigg\}.
\end{equation*}
Then
\begin{equation*}
	X=I_{\varphi}(f)\cup\bigg(\bigcup_{a\in L_{\varphi}}R_{\varphi}(a)\bigg).
\end{equation*}
There is a substantial literature on multifractal analysis; see, for example, \cite{Barreira-Saussol-2000,Barreira-Saussol-2001,Chen-Tian-2021,Climenhaga-2013,Pfister-Sullivan-2007,Shi-Tian-Varandas-Wang-2023}.

According to Birkhoff's ergodic theorem, $\mu(\mathrm{IR}(f))=0$ holds for all $\mu\in\mathcal{M}_f(X)$.
However, various aspects of dynamical complexity have been extensively studied for $\varphi$-irregular sets, including residuality, topological entropy, topological pressure, distributional chaos, and so on.
See \cite{Li-Wu-2014}, \cite{Thompson-2012,Tian-2017}, and \cite{Li-Yorke-1975,Schweizer-Smital-1994,Chen-Tian-2021} for residuality, topological pressure, and Li--Yorke chaos and distributional chaos, respectively.
In this article, we mainly consider the Bowen topological entropy of $\varphi$-irregular sets.
Bowen \cite{Bowen-1973} introduced topological entropy for non-compact sets, and Pesin and Pitskel' \cite{Pesin-Pitskel-1984} first proved that the irregular set of the two-symbol full shift has full topological entropy.
Full topological entropy of irregular sets was later established for topologically mixing subshifts of finite type \cite{Barreira-Schmeling-1997,Barreira-Schmeling-2000,Fan-Feng-Wu-2001}, systems with the specification property \cite{Chen-Kupper-Shu-2005,Thompson-2010}, and systems with the shadowing property \cite{Dong-Oprocha-Tian-2018}.

Saturated sets and over-saturated sets can be used to study the topological entropy of irregular sets.
Indeed, it can be checked that if $|\{\int\varphi\mathrm{d}\mu:\mu\in K\}|\geq2$ for some $K\subset\mathcal{M}_f(X)$, then $G^K\subset I_\varphi(f)$.
%For dynamical systems with the non-uniform specification property, it follows from Lemma \ref{pre-saturated-lemma6} that every non-empty $\varphi$-irregular set contains a non-empty over-saturated set.
In \cite{Lin-Tian-Yu-2024}, Lin, Tian and Yu provided some basic properties for dynamical systems with the non-uniform specification property, including the residuality of non-empty over-saturated sets and the completely-irregular set.
It has been shown that, under the non-uniform specification property, every non-empty over-saturated set $G^K$ is residual in $X$, and the completely-irregular set $\mathrm{CI}(f)$ is residual in $X$.
These results answer Question \ref{Que1} in the sense of residuality.

\begin{theorem}[{\cite[Theorem 3.2 (4),(8)]{Lin-Tian-Yu-2024}}]
	Suppose that $(X,f)$ satisfies the non-uniform specification property, then
	\begin{enumerate}
		\item every non-empty over-saturated set $G^K$ is residual in $X$;
		\item $\mathfrak{C}(f)\neq\varnothing$, and $\mathrm{CI}(f)$ is residual in $X$.
	\end{enumerate}
\end{theorem}

We focus on the topological entropy of irregular sets.

\subsubsection{Lower bounds for Bowen entropy}

Let $h_{top}^B(f,E)$ denote the Bowen topological entropy of a subset $E\subset X$.
Clearly, $h_{top}^B(f,E)\leq h_{top}(f)$.
As established by Bowen in \cite{Bowen-1973}, for any ergodic measure $\nu$, the topological entropy of the set of generic points $G_{\nu}$ coincides with the measure-theoretic entropy $h_{\nu}(f)$.
Furthermore, Bowen demonstrated in \cite{Bowen-1973} that for any $t\geq0$, the topological entropy of $\mathrm{QR}(t)$ is bounded above by $t$ (see Proposition \ref{Bowen-upper-prop1}), where
\begin{equation*}
	\mathrm{QR}(t)=\{x\in X: \text{there is a}\;\mu\in V_f(x)\;\text{such that}\;h_{\mu}(f)\leq t\}.
\end{equation*}
These results require no additional assumptions.

Non-trivial lower bounds for the Bowen topological entropy of non-compact sets generally require additional assumptions.
Under the $g$-almost product property and the uniform separation property, Pfister and Sullivan \cite{Pfister-Sullivan-2007} proved that
\begin{equation*}
	h_{top}^B(f,G_K)=\inf\{h_\mu(f):\mu\in K\}
\end{equation*}
for every non-empty compact connected set $K\subset\mathcal{M}_f(X)$.
This result was extended to topological pressure \cite{Pei-Chen-2010}, non-uniformly hyperbolic systems \cite{Liang-Liao-Sun-Tian-2017}, systems with non-uniform structure \cite{Zhao-Chen-2018}, and amenable group actions \cite{Ren-Tian-Zhou-2023}.

Before presenting the main results, we introduce $\sigma$-controlled functions as defined in \cite{Hou-Lin-Tian-2023}, which will be used in the statement of our main theorems.

\begin{definition}
	Let $0<\sigma<1$.
	We say $\varphi\in C(X)$ is \emph{$\sigma$-controlled} if
	\begin{equation}\label{eq-intro-Bowen-sigma-controlled}
		\sup_{\mu\in\mathcal{M}_f(X)}\int_X\varphi\mathrm{d}\mu-\inf_{\mu\in\mathcal{M}_f(X)}\int_X\varphi\mathrm{d}\mu>\sigma\Big(\sup_{x\in X}\varphi(x)-\inf_{x\in X}\varphi(x)\Big).
	\end{equation}
	The set of all $\sigma$-controlled continuous functions is denoted by $C_{\sigma}(X)$.
\end{definition}

For any $\mu\in\mathcal{M}(X)$, it is clear that
\begin{equation*}
	\inf_{x\in X}\varphi(x)\leq\int_X\varphi\mathrm{d}\mu\leq\sup_{x\in X}\varphi(x).
\end{equation*}
Hence $C_{\sigma}(X)=\varnothing$ for $\sigma\geq1$.
On the other hand, Corollary \ref{pre-saturated-cor4} indicates that under the non-uniform specification property, we have $C_{\sigma}(X)\subset\mathfrak{C}(f)$ for all $0<\sigma<1$, and consequently
\begin{equation*}
	\mathfrak{C}(f)=\bigcup_{0<\sigma<1}C_{\sigma}(X).
\end{equation*}
For any $0<\sigma<1$ and any $\varphi\in\mathfrak{C}(f)$, define
\begin{equation*}
	\mathfrak{D}_f(\varphi,\sigma)=\Bigg\{(\mu,\nu)\in\mathcal{M}_f(X)\times\mathcal{M}_f(X):\bigg|\int_X\varphi\mathrm{d}\mu-\int_X\varphi\mathrm{d}\nu\bigg|>\sigma\Big(\sup_{x\in X}\varphi(x)-\inf_{x\in X}\varphi(x)\Big)\Bigg\}.
\end{equation*}
Then by \eqref{eq-intro-Bowen-sigma-controlled}, we have $\mathfrak{D}_f(\varphi,\sigma)\neq\varnothing$ for all $\varphi\in C_{\sigma}(X)$.
For $\sigma=0$, we define $C_0(X)$ and $\mathfrak{D}_f(\varphi,0)$ by the same inequalities.

In \cite{Thompson-2012}, Thompson proved that under the almost specification property, every non-empty $\varphi$-irregular set carries full topological entropy.
In this paper, we study how the asymptotic behavior of gap functions in dynamical systems satisfying the non-uniform specification property influences the Bowen topological entropy of irregular sets.
Consider a dynamical system with the non-uniform $\tau$-specification property for some $\tau>0$.
Intuitively, when tracing orbit segments, the gaps occupy a proportion of at least $\tau/(1+\tau)$ of the total time, and the orbit is not controlled during these gaps.
This suggests that a natural lower bound for the Bowen entropy of $\mathrm{IR}(f)$ should involve the factor $1/(1+\tau)$, which motivates the following question.

\begin{question}\label{Que2}
Given $\tau\geq 0$, let $(X,f)$ be a dynamical system satisfying the non-uniform $\tau$-specification property.
Do we have that 
\begin{equation*}
	h_{top}^B(f,\mathrm{IR}(f))\geq\frac{1}{1+\tau}h_{top}(f)?
\end{equation*}
\end{question}

Our first main results concern the $\varphi$-irregular sets for $\sigma$-controlled functions.
We provide a lower bound for the Bowen topological entropy of the $\varphi$-irregular set for a $\sigma$-controlled function $\varphi$ under the non-uniform specification property, where the bound is determined by the lower asymptotic linear growth rate of the gap function.
Furthermore, we estimate the Bowen topological entropy of the irregular set $\mathrm{IR}(f)$, which partially answers Question \ref{Que2}.

\begin{maintheorem}\label{thm-Bowen-lower-1}
	Suppose that $(X,f)$ satisfies the non-uniform specification property with gap function $M(n,\varepsilon)$.
	Assume that
	\begin{equation*}
		\tau:=\tau(M)=\sup_{\varepsilon>0}\liminf_{n\to\infty}\frac{M(n,\varepsilon)}{n}<1.
	\end{equation*}
	Then
	\begin{enumerate}
		\item for any $0<\sigma<1$ with $\tau\leq\sigma$ and any $\varphi\in C_{\sigma}(X)$, we have
		\begin{equation*}
			h_{top}^B(f,I_{\varphi}(f))\geq\frac{\sigma-\tau}{(1+\tau)\sigma}h_{top}(f)+\frac{\tau}{(1+\tau)\sigma}\sup_{(\mu,\nu)\in\mathfrak{D}_f(\varphi,\sigma)}\min\{h_{\mu}(f),h_{\nu}(f)\};
		\end{equation*}
		\item
		\begin{equation*}
			h_{top}^B(f,\mathrm{IR}(f))\geq\frac{1-\tau}{1+\tau}h_{top}(f)+\frac{\tau}{1+\tau}\sup_{\mu,\nu\in\mathcal{M}_f^{erg}(X),\mu\neq\nu}\min\{h_{\mu}(f),h_{\nu}(f)\};
		\end{equation*}
		\item
		\begin{equation*}
			h_{top}^B(f,\mathrm{IR}(f))\geq\Big(\frac{1}{1+\tau}-\frac{\tau^2}{(1+\tau)^2}\Big)h_{top}(f).
		\end{equation*}
	\end{enumerate}
	In particular, if $\tau=0$, then for any $\varphi\in\mathfrak{C}(f)$, we have
	\begin{equation*}
		h_{top}^B(f,I_{\varphi}(f))=h_{top}(f).
	\end{equation*}
\end{maintheorem}

\begin{remark}
	The case $\sigma=\tau$ in Theorem \ref{thm-Bowen-lower-1} (1) can be viewed as a $\tau$-controlled function version of Question \ref{Que2}.
	Indeed, this case plays a key role in the proof of Theorem \ref{thm-Bowen-lower-1}, and we state it independently as Theorem \ref{Bowen-thm2}.
\end{remark}

\begin{remark}
	Theorem \ref{thm-Bowen-lower-1} does not fully resolve Question \ref{Que2}.
	Specifically, our proof applies only when $\tau<1$, as no $\sigma$-controlled functions exist when $\sigma\geq1$.
	Furthermore, the factor $1/(1+\tau)$ in Question \ref{Que2} cannot in general be obtained by this proof.
	When applying Theorem \ref{thm-Bowen-lower-1} (2), we need to estimate
	\begin{equation*}
		\sup_{\mu,\nu\in\mathcal{M}_f^{erg}(X),\mu\neq\nu}\min\{h_{\mu}(f),h_{\nu}(f)\}.
	\end{equation*}
	Note that we cannot guarantee that this supremum coincides with the topological entropy $h_{top}(f)$.
	See Proposition \ref{Bowen-prop5}, Corollary \ref{Bowen-cor6} and Remark \ref{Bowen-rmk7}.
	Hence we only obtain the factor $(1/(1+\tau)-\tau^2/(1+\tau)^2)$ for the estimation in general.
	However, we can give a positive answer to Question \ref{Que2} for $\tau<1$ if $(X,f)$ satisfies one of the following conditions:
	\begin{itemize}
		\item $(X,f)$ has no measures of maximal entropy;
		\item $(X,f)$ has only one measure $\mu$ of maximal entropy but $h_\mu(f)=\sup\{h_\nu(f):\mu\neq\nu\in\mathcal{M}_f^{erg}(X)\}$;
		\item $(X,f)$ has at least two distinct ergodic measures of maximal entropy.
	\end{itemize}
\end{remark}

\begin{maincorollary}\label{cor-Bowen-lower-1}
	Suppose that $(X,f)$ satisfies the non-uniform $\tau$-specification property.
	Assume that one of the following conditions holds:
	\begin{enumerate}
		\item $(X,f)$ has no measures of maximal entropy;
		\item $(X,f)$ has only one measure $\mu$ of maximal entropy but $h_\mu(f)=\sup\{h_\nu(f):\mu\neq\nu\in\mathcal{M}_f^{erg}(X)\}$;
		\item $(X,f)$ has at least two distinct ergodic measures of maximal entropy.
	\end{enumerate}
	If $\tau<1$, then
	\begin{equation*}
		h_{top}^B(f,\mathrm{IR}(f))\geq\frac{1}{1+\tau}h_{top}(f).
	\end{equation*}
\end{maincorollary}

\begin{remark}
	In \cite{Pavlov-2016}, Pavlov constructed a subshift with non-uniform specification, admitting exactly two ergodic measures of maximal entropy.
	However, for this example, $\tau=0$.
	We will adapt this example in Theorem \ref{thm-construction}.
\end{remark}

\subsubsection{Upper bounds for Bowen entropy}

The following theorem refines the examples in \cite{Lin-Tian-Yu-2024} and is used to prove several main results.
See Figure \ref{fig:application-of-thm-construction}.

\begin{figure}[h]
    \centering
    \begin{tikzpicture}
        \node[draw, line width=0.8pt, minimum width=2.5cm, minimum height=0.8cm] (B) at (0,0) {Theorem \ref{thm-construction}};
        \node[draw, line width=0.8pt, minimum width=2.2cm, minimum height=0.8cm] (C) at (-3.5,-2) {Theorem \ref{thm-Bowen-upper-1}};
        \node[draw, line width=0.8pt, minimum width=2.2cm, minimum height=0.8cm] (D) at (3.5,-2) {Theorem \ref{thm-Bowen-upper-2}};
        \node[draw, line width=0.8pt, minimum width=1.8cm, minimum height=0.8cm] (E) at (-6.8,-4) {Theorem \ref{thm-CVP}};
        \node[draw, line width=0.8pt, minimum width=1.8cm, minimum height=0.8cm] (J) at (-4.6,-4) {Theorem \ref{thm-growth-2}};
        \node[draw, line width=0.8pt, minimum width=1.8cm, minimum height=0.8cm] (K) at (-2.4,-4) {Theorem \ref{thm-intermediate-1}};
        \node[draw, line width=0.8pt, minimum width=1.8cm, minimum height=0.8cm] (N) at (-0.2,-4) {Theorem \ref{thm-ergodic-optimization}};
        \node[draw, line width=0.8pt, minimum width=1.8cm, minimum height=0.8cm] (G) at (2.4,-4) {Theorem \ref{thm-trans-1}};
        \node[draw, line width=0.8pt, minimum width=1.8cm, minimum height=0.8cm] (L) at (4.6,-4) {Theorem \ref{thm-intermediate-2}};
        \draw[->, line width=0.8pt] (B.south) -- (C.north);
        \draw[->, line width=0.8pt] (B.south) -- (D.north);
        \draw[->, line width=0.8pt] (D.south) -- (G.north);
        \draw[->, line width=0.8pt] (C.south) -- (E.north);
        \draw[->, line width=0.8pt] (C.south) -- (J.north);
        \draw[->, line width=0.8pt] (C.south) -- (K.north);
        \draw[->, line width=0.8pt] (C.south) -- (N.north);
        \draw[->, line width=0.8pt] (D.south) -- (L.north);
    \end{tikzpicture}
	\caption{Connections and derivations among several main results.}
	\label{fig:application-of-thm-construction}
\end{figure}

\begin{maintheorem}\label{thm-construction}
	Given $\tau\in(0,\infty]$, a finite alphabet $\mathcal{A}$, and a subshift $Z\subset\mathcal{A}^{\mathbb{N}_0}$, there exists a one-sided subshift $(X,f)$ satisfying the non-uniform periodic $\tau$-specification property with gap function $M(n,\varepsilon)$ such that
	\begin{enumerate}
		\item for $\varepsilon>0$ sufficiently small, we have
		\begin{equation*}
			\lim_{n\rightarrow\infty}\frac{M(n,\varepsilon)}{n}=\tau;
		\end{equation*}
		\item $Z\subset X$;
		\item $\mathcal{M}_f^{erg}(X)$ can be written as $\mathcal{M}_f^{erg}(X)=\mathcal{M}_0\cup\mathcal{M}_f^{erg}(Z)$ such that $\overline{\mathcal{M}_0}\cap\overline{\mathcal{M}_f^{erg}(Z)}=\varnothing$;
		\item for any $x\in X\setminus\bigcup_{i\geq 0}f^{-i}Z$, we have $V_f(x)\cap\overline{\mathcal{M}_0}\neq\varnothing$;
		\item if $0<\tau<\infty$, then $(X,f)$ can be constructed such that for any $\mu\in\overline{\mathcal{M}_0}$, we have
		\begin{equation*}
			h_{\mu}(f)\leq \log 2+\frac{1}{1+\tau}\log|\mathcal{A}|;
		\end{equation*}
		if $\tau=\infty$, then for any $\delta>0$, $(X,f)$ can be constructed such that for any $\mu\in\overline{\mathcal{M}_0}$, we have
		\begin{equation*}
			h_{\mu}(f)\leq \log 2+\delta\log|\mathcal{A}|;
		\end{equation*}
	\end{enumerate}
\end{maintheorem}

\begin{remark}\label{Remark 1.7}
	In the proof of Theorem \ref{thm-construction}, the alphabet is enlarged by one marker symbol.
	The same construction works for two-sided subshifts, except that (4) fails.
\end{remark}

Since $Z$ is arbitrary, Theorem \ref{thm-construction} gives various systems with non-uniform $\tau$-specification.
In particular, when $0<\tau<\infty$, Theorem \ref{thm-construction} (5) and the variational principle give $h_{top}(X)=h_{top}(Z)$ whenever
\begin{equation*}
	h_{top}(Z)\geq\log2+\frac{\log|\mathcal{A}|}{1+\tau}.
\end{equation*}
Thus any sufficiently large topological entropy can be realized when $0<\tau<\infty$.

\begin{maincorollary}\label{cor-entropy}
	Given $0<\tau<\infty$ and $\log\Big(2\big\lceil 2^{\frac{1+\tau}{\tau}}+\frac1\tau\big\rceil^{\frac{1}{1+\tau}}\Big)<h<\infty$, there exists a dynamical system $(X,f)$ satisfying the non-uniform $\tau$-specification property with gap function $M(n,\varepsilon)$ such that
	\begin{enumerate}
		\item for $\varepsilon>0$ sufficiently small, we have
		\begin{equation*}
			\lim_{n\rightarrow\infty}\frac{M(n,\varepsilon)}{n}=\tau;
		\end{equation*}
		\item $h_{top}(f)=h$.
	\end{enumerate}
\end{maincorollary}

Such subsystems also occur in differential dynamical systems.
Consider a $C^1$ diffeomorphism $f$ on a compact $C^\infty$ Riemannian manifold $M$.
Let $\Lambda$ be an $f$-invariant set.
A non-trivial $Df$-invariant splitting $TM|_\Lambda=F\oplus E$ is called \emph{dominated}, denoted by $F\oplus_\prec E$, if there exist $C>0$ and $\lambda\in(0,1)$ such that
\begin{equation*}
	\|Df^n|_{F(x)}\|\cdot\|Df^{-n}|_{E(f^nx)}\|\leq C\lambda^n
\end{equation*}
for every $x\in\Lambda$ and $n\in\mathbb{N}$.
An $f$-invariant set $\Lambda$ is said to admit a \emph{dominated splitting} if there exists a non-trivial $Df$-invariant splitting $TM|_\Lambda=F\oplus_\prec E$.

Let $\mu\in\mathcal{M}_f^{erg}(M)$ be a \emph{hyperbolic measure}, i.e. with no zero Lyapunov exponents, and assume that it has at least one positive and one negative Lyapunov exponent.
Write its Oseledets splitting as
\begin{equation*}
	E_1\oplus\cdots\oplus E_l\oplus E_{l+1}\oplus\cdots\oplus E_{l+s},
\end{equation*}
with corresponding exponents
\begin{equation*}
	\lambda_1<\cdots<\lambda_l<0<\lambda_{l+1}<\cdots<\lambda_{l+s}.
\end{equation*}
Set $E^-=E_1\oplus\cdots\oplus E_l$ and $E^+=E_{l+1}\oplus\cdots\oplus E_{l+s}$.
We say that $\mu$ admits a \emph{dominated splitting} corresponding to its stable and unstable Oseledets subspaces if its support $S_\mu$ admits a dominated splitting $F\oplus_\prec E$ such that $F(x)=E^-(x)$ and $E(x)=E^+(x)$ for $\mu$-a.e. $x\in M$.

\begin{maincorollary}\label{Corollary C New}
	Suppose that $f$ is a $C^1$ diffeomorphism on a compact $C^\infty$ Riemannian manifold $M$ and $\mu$ is a hyperbolic ergodic measure with $h_{\mu}(f)>0$. Suppose also that one of the following conditions is satisfied:
	\begin{enumerate}
		\item $f$ is a $C^{1+\alpha}$ diffeomorphism for some $\alpha>0$;
		\item $\mu$ admits a dominated splitting corresponding to the stable/unstable subspaces of its Oseledets splitting.
	\end{enumerate}
	Then for any $0<\tau<\infty$ and $0<h<h_\mu(f)$,  there exist $N\in\mathbb{N}$ and a nonempty compact $f^N$-invariant subset $Y\subset M$ satisfying the following:
	\begin{enumerate}
		\item $f^i(Y)\cap f^j(Y)=\varnothing$ for any $0\leq i<j\leq N-1$;
		\item $(Y,f^N)$ satisfies the non-uniform $\tau$-specification property with gap function $M(n,\varepsilon)$ such that for $\varepsilon>0$ sufficiently small, we have
		\begin{equation*}
			\lim_{n\rightarrow\infty}\frac{M(n,\varepsilon)}{n}=\tau;
		\end{equation*}
		\item $h_{top}(f^N|_Y)=Nh$. 
	\end{enumerate}
	In particular, $(\Lambda,f|_\Lambda)$ satisfies the relative non-uniform $\tau$-specification property with disjoint steps and $h_{top}(f|_\Lambda)=h$, where $\Lambda=\bigcup_{i=0}^{N-1}f^i(Y)$.
\end{maincorollary}

\begin{remark}
	The relative non-uniform $\tau$-specification property with disjoint steps is defined in Section \ref{Section 2.6}.
\end{remark}

By choosing the alphabet $\mathcal{A}$ and the subshift $Z\subset\mathcal{A}^{\mathbb{N}_0}$, we obtain examples showing that the factor $1/(1+\tau)$ in Question \ref{Que2} is optimal.

\begin{maintheorem}\label{thm-Bowen-upper-1}
	Given $\tau\in(0,\infty]$ and $0<\delta<1-1/(1+\tau)$, for any $m\in\mathbb{N}$, there exists a dynamical system $(X,f)$ satisfying the non-uniform $\tau$-specification property with gap function $M(n,\varepsilon)$ such that
	\begin{enumerate}
		\item for $\varepsilon>0$ sufficiently small, we have
		\begin{equation*}
			\lim_{n\rightarrow\infty}\frac{M(n,\varepsilon)}{n}=\tau;
		\end{equation*}
		\item $(X,f)$ has exactly $m$ ergodic measures of maximal entropy $\nu_1,\cdots,\nu_m$;
		\item $\nu_1,\cdots,\nu_m$ are $m$ isolated points in $\mathcal{M}_f^{erg}(X)$.
		\item the irregular set $\mathrm{IR}(f)$ does not carry full Bowen topological entropy, and more precisely we have
		\begin{equation*}
			h_{top}^B(f,\mathrm{IR}(f))\leq\bigg(\frac{1}{1+\tau}+\delta\bigg)\cdot h_{top}(f).
		\end{equation*}
	\end{enumerate}
\end{maintheorem}

The next theorem shows that full topological entropy is the optimal upper bound for irregular sets and gives an example for which the completely-irregular set has strictly smaller Bowen topological entropy than the irregular set.

\begin{maintheorem}\label{thm-Bowen-upper-2}
	For any $\tau\in(0,\infty]$ and $0<\delta<1-1/(1+\tau)$, there exists a dynamical system $(X,f)$ satisfying the non-uniform $\tau$-specification property with gap function $M(n,\varepsilon)$ such that
	\begin{enumerate}
		\item for any $\varepsilon>0$ sufficiently small, we have
		\begin{equation*}
			\lim_{n\rightarrow\infty}\frac{M(n,\varepsilon)}{n}=\tau;
		\end{equation*}
		\item there is a $\psi\in\mathfrak{C}(f)$ with $h_{top}^B(f,I_{\psi}(f))=h_{top}(f)$, and consequently $h_{top}^B(f,\mathrm{IR}(f))=h_{top}(f)$;
		\item there is a $\varphi\in\mathfrak{C}(f)$ with
		\begin{equation*}
			h_{top}^B(f,I_{\varphi}(f))\leq\bigg(\frac{1}{1+\tau}+\delta\bigg)\cdot h_{top}(f)<h_{top}(f),
		\end{equation*}
		and consequently
		\begin{equation*}
			h_{top}^B(f,\mathrm{CI}(f))\leq\bigg(\frac{1}{1+\tau}+\delta\bigg)\cdot h_{top}(f)<h_{top}(f).
		\end{equation*}
	\end{enumerate}
\end{maintheorem}

\subsubsection{Packing entropy of over-saturated sets}

Feng and Huang \cite{Feng-Huang-2012} introduced the notion of packing topological entropy in a way resembling packing dimension.
For each subset $E\subset X$, let $h_{top}^P(f,E)$ denote the packing topological entropy of $E$.
We consider the packing topological entropy of over-saturated sets.
It was proved in \cite{Hou-Tian-Zhang-2023} that for any dynamical system with the $g$-almost product property, every non-empty connected compact convex set $K\subset\mathcal{M}_f(X)$ satisfies
\begin{equation*}
	h_{top}^P(f,G_K)=\sup\{h_{\mu}(f):\mu\in K\}.
\end{equation*}
Consequently, for every non-empty set $K\subset\mathcal{M}_f(X)$,
\begin{equation*}
	h_{top}^P(f,G^K)\geq h_{top}^P(f,G_{\overline{\mathrm{co}}(K)})\geq\sup\{h_{\mu}(f):\mu\in K\},
\end{equation*}
where $\overline{\mathrm{co}}(K)$ denotes the closed convex hull of $K$.

In contrast to Bowen topological entropy, every over-saturated set under non-uniform specification is either empty or carries full packing topological entropy.

\begin{maintheorem}\label{thm-packing}
	Suppose that $(X,f)$ satisfies the non-uniform specification property, then
	\begin{enumerate}
		\item for any non-empty over-saturated set $G^K$, we have
		\begin{equation*}
			h_{top}^P(f,G^K)=h_{top}(f).
		\end{equation*}
		\item $h_{top}^P\big(f,\mathrm{CI}(f)\big)=h_{top}(f)$.
	\end{enumerate}
\end{maintheorem}

\subsubsection{Conditional variational principle}

We say that the \emph{conditional variational principle} holds if
\begin{equation*}
	h_{top}^B(f,R_{\varphi}(a))=\sup\bigg\{h_{\mu}(f):\mu\in\mathcal{M}_f(X),\int_X\varphi\mathrm{d}\mu=a\bigg\}
\end{equation*}
for every $\varphi\in C(X)$ with $\mathrm{Int}(L_{\varphi})\neq\varnothing$ and every $a\in\mathrm{Int}(L_{\varphi})$.
Using arguments in \cite{Pfister-Sullivan-2007}, we can show that this principle holds when $\tau=0$, see Corollary \ref{CVP-cor3}.
For $\tau>0$, the examples constructed in Theorem \ref{thm-Bowen-upper-1} show that it may fail.

\begin{maintheorem}\label{thm-CVP}
	For any $\tau\in(0,\infty]$, there exists a dynamical system $(X,f)$ satisfying the non-uniform $\tau$-specification property with gap function $M(n,\varepsilon)$ such that
	\begin{enumerate}
		\item for $\varepsilon>0$ sufficiently small, we have
		\begin{equation*}
			\lim_{n\rightarrow\infty}\frac{M(n,\varepsilon)}{n}=\tau;
		\end{equation*}
		\item there exist $\psi\in C(X)$ and $0<\gamma<1$ such that $L_{\psi}=[0,1]$ and we have $R_{\psi}(a)=\varnothing$ for all $\gamma<a<1$, while
		\begin{equation*}
			\sup\bigg\{h_{\mu}(f):\mu\in\mathcal{M}_f(X),\int_X\psi\mathrm{d}\mu=a\bigg\}>0=h_{top}^B(f,R_{\psi}(a)),
		\end{equation*}
		and consequently, the conditional variational principle does not hold.
	\end{enumerate}
\end{maintheorem}

\begin{remark}
	If $0<\tau<\infty$, then one may take $\gamma=1/(1+\tau)$ in Theorem \ref{thm-CVP}.
	If $\tau=\infty$, then $\gamma$ can be chosen arbitrarily small.
\end{remark}

\subsection{Transitive points}

For $x\in X$, let $\omega(x,f)$ denote the \emph{$\omega$-limit set} of $x$, i.e. the limit set of $\{f^n(x):n\geq 0\}$.
A point $x\in X$ is said to be
\begin{itemize}
	\item a \emph{transitive point}, if $\omega(x,f)=X$;
	\item a \emph{recurrent point}, if $x\in\omega(x,f)$.
\end{itemize}
The sets of transitive points and recurrent points in $(X,f)$ are denoted by $\mathrm{Trans}(f)$ and $\mathrm{Rec}(f)$, respectively.
Clearly $\mathrm{Trans}(f)\subset\mathrm{Rec}(f)$.
According to \cite[Theorem 3.1]{Tian-2016}, $\mathrm{Rec}(f)$ always carries full Bowen topological entropy.
However, for every $\tau>0$, Theorems \ref{thm-Bowen-upper-1} and \ref{thm-Bowen-upper-2} provide systems with the non-uniform $\tau$-specification property for which $\mathrm{Trans}(f)$ does not carry full Bowen topological entropy.

\begin{maintheorem}\label{thm-trans-1}
	For any $\tau\in(0,\infty]$ and $0<\delta<1-1/(1+\tau)$, there exists a dynamical system $(X,f)$ satisfying the non-uniform $\tau$-specification property such that $\mathrm{Trans}(f)$ does not carry full topological entropy.
	More precisely we have
	\begin{equation*}
		h_{top}^B(f,\mathrm{Trans}(f))\leq\bigg(\frac{1}{1+\tau}+\delta\bigg)\cdot h_{top}(f)<h_{top}(f).
	\end{equation*}
	As a result, we have $h_{top}^B(f,\mathrm{Trans}(f))<h_{top}^B(f,\mathrm{Rec}(f))$.
\end{maintheorem}

In Theorem \ref{thm-Bowen-lower-1}, we estimate the lower bound for Bowen topological entropy of certain irregular sets.
Although the proof of Theorem \ref{thm-Bowen-lower-1} does not apply when $\tau\geq 1$, we can answer a modified version of Question \ref{Que2}, where the set $\mathrm{Trans}(f)\setminus(\bigcup_{\mu\in\mathcal{M}_f^{erg}(X)}G_{\mu})$ replaces $\mathrm{IR}(f)$ and the factor $1/(1+\tau)^2$ replaces $1/(1+\tau)$.

Note that
\begin{equation*}
	\mathrm{IR}(f)\cup\bigg(\bigcup_{\omega\in\mathcal{M}_f(X)\setminus\mathcal{M}_f^{erg}(X)}G_{\omega}\bigg)=X\setminus\bigg(\bigcup_{\omega\in\mathcal{M}_f^{erg}(X)}G_{\omega}\bigg).
\end{equation*}
According to Birkhoff's ergodic theorem, this set is also of zero measure for any invariant measure.

\begin{maintheorem}\label{thm-trans-2}
	Suppose that $(X,f)$ satisfies the non-uniform $\tau$-specification property.
	Then
	\begin{equation*}
		\max\bigg\{h_{top}^B\Big(f,\mathrm{IR}(f)\cap\mathrm{Trans}(f)\Big),h_{top}^B\Big(f,\bigcup_{\omega\in\mathcal{M}_f(X)\setminus\mathcal{M}_f^{erg}(X)}G_{\omega}\cap\mathrm{Trans}(f)\Big)\bigg\}\geq\frac{1}{(1+\tau)^2}h_{top}(f).
	\end{equation*}
\end{maintheorem}

\subsection{Growth of periodic orbits}

For any set $S$, let $|S|$ denote the cardinality of $S$.

A homeomorphism $f$ on $X$ is \emph{expansive} if there exists $\varepsilon_0>0$ such that $x\neq y$ implies $d(f^ix,f^iy)\geq\varepsilon_0$ for some $i\in\mathbb{Z}$.
A continuous map $f:X\to X$ is \emph{positively expansive} if there exists $\varepsilon_0>0$ such that $x\neq y$ implies $d(f^ix,f^iy)\geq\varepsilon_0$ for some $i\in\mathbb{N}_0$.
In either case, $\varepsilon_0$ is called an \emph{expansive constant}.
According to Bowen's classic results \cite{Bowen-1971,Bowen-1974}, when $(X,f)$ is (positively) expansive and satisfies the periodic specification property, we always have
\begin{equation*}
	\lim_{n\rightarrow\infty}\frac{1}{n}\log|\mathrm{Per}_n(f)|=h_{top}(f).
\end{equation*}
For dynamical systems satisfying the non-uniform periodic specification property, we can exhibit a lower bound for the exponential periodic growth rate, which is determined by the asymptotic behavior of the gap function.

\begin{maintheorem}\label{thm-growth-1}
	Suppose that $(X,f)$ satisfies the non-uniform periodic $\tau$-specification property.
	Then
	\begin{equation*}
		\frac{1}{1+\tau}h_{top}(f)\leq\liminf_{n\rightarrow\infty}\frac1n\log|\mathrm{Per}_n(f)|.
	\end{equation*}
	If we further assume that $(X,f)$ is positively expansive, then
	\begin{equation*}
		\frac{1}{1+\tau}h_{top}(f)\leq\liminf_{n\rightarrow\infty}\frac1n\log|\mathrm{Per}_n(f)|\leq\limsup_{n\rightarrow\infty}\frac1n\log|\mathrm{Per}_n(f)|\leq h_{top}(f).
	\end{equation*}
\end{maintheorem}

The examples in \cite[Theorem A]{Lin-Tian-Yu-2024} have exponential periodic growth rate equal to the topological entropy, showing that the latter is the optimal upper bound.
Sigmund \cite{Sigmund-1970} proved that periodic measures are dense in the invariant measures under the periodic specification property.
In \cite[Theorem 3.4]{Lin-Tian-Yu-2024}, Lin, Tian and Yu showed that this property still holds for dynamical systems satisfying the non-uniform periodic $0$-specification property.
Moreover, it was shown in \cite{Lin-Tian-Yu-2024} that the condition $\tau=0$ is optimal.
More precisely, when $\tau>0$, each subshift $(X,f)$ constructed in \cite[Theorem A]{Lin-Tian-Yu-2024} satisfies
\begin{equation*}
	\overline{\mathcal{M}_f^{co}(X)}=\overline{\mathcal{M}_f^{erg}(X)}\subsetneq\mathcal{M}_f(X).
\end{equation*}

The system in Theorem \ref{thm-Bowen-upper-1} satisfies $\overline{\mathcal{M}_f^{co}(X)}\subsetneq\overline{\mathcal{M}_f^{erg}(X)}\subsetneq\mathcal{M}_f(X)$ and
\begin{equation*}
	\limsup_{n\rightarrow\infty}\frac{1}{n}\log|\mathrm{Per}_n(f)|<h_{top}(f).
\end{equation*}
The next theorem shows that the factor $1/(1+\tau)$ in Theorem \ref{thm-growth-1} is optimal.

\begin{maintheorem}\label{thm-growth-2}
	Given $\tau\in(0,\infty]$ and $0<\delta<1-1/(1+\tau)$, there exists a positively expansive dynamical system $(X,f)$ satisfying the non-uniform periodic $\tau$-specification property such that
	\begin{enumerate}
		\item $\overline{\mathcal{M}_f^{co}(X)}\subsetneq\overline{\mathcal{M}_f^{erg}(X)}\subsetneq\mathcal{M}_f(X)$;
		\item the exponential periodic growth rate is strictly less than the topological entropy, and more precisely we have
		\begin{equation*}
			\limsup_{n\rightarrow\infty}\frac1n\log|\mathrm{Per}_n(f)|\leq\bigg(\frac{1}{1+\tau}+\delta\bigg)\cdot h_{top}(f)<h_{top}(f).
		\end{equation*}
	\end{enumerate}
\end{maintheorem}

\begin{remark}
	One can easily see that the subshift $(X,f)$ defined in Theorem \ref{thm-Bowen-upper-2} satisfies
	\begin{equation*}
		h_{top}(f)=\lim_{n\rightarrow\infty}\frac{1}{n}\log|\mathrm{Per}_n(f)|.
	\end{equation*}
	Together with Theorem \ref{thm-growth-2}, this shows that the bounds in Theorem \ref{thm-growth-1} are optimal.
\end{remark}

\subsection{Intermediate entropy and pressure}

We say that $(X,f)$ satisfies the \emph{intermediate entropy property} if every $0\leq h<h_{top}(f)$ is the entropy of some $\mu\in\mathcal{M}_f^{erg}(X)$.
Katok \cite{Katok-1980} proved that every $C^{1+\alpha}$ diffeomorphism on a compact surface with positive topological entropy has horseshoes of large entropies and hence satisfies this property.
See \cite{Burguet-2020,Guan-Sun-Wu-2017,Hou-Tian-2024,Huang-Xu-Xu-2021,Konieczny-Kupsa-Kwietniak-2018,Li-Oprocha-2018,Li-Shi-Wang-Wang-2020,Sun-2025} for further results.

Under the non-uniform $0$-specification property, the intermediate entropy property holds under either of the following additional assumptions:
\begin{enumerate}
	\item $(X,f)$ is asymptotically entropy expansive, by combining \cite[Lemma 3.5]{Lin-Tian-Yu-2024} and \cite[Theorem 1.4]{Sun-2025};
	\item $\sup_{\varepsilon>0}\lim_{n\rightarrow\infty}\frac{M(n,\varepsilon)}{n}=0$, by using arguments in proof of \cite[Theorem 1.2]{Burguet-2020}.
\end{enumerate}
When $\tau>0$, the intermediate entropy property may not hold under non-uniform $\tau$-specification.

\begin{maintheorem}\label{thm-intermediate-1}
	Given $\tau\in(0,\infty]$ and $0<\delta<1-1/(1+\tau)$, there exists a (positively expansive) dynamical system $(X,f)$ satisfying the non-uniform $\tau$-specification property such that for any $h$ with
	\begin{equation*}
		\bigg(\frac{1}{1+\tau}+\delta\bigg)\cdot h_{top}(f)<h<h_{top}(f),
	\end{equation*}
	there is no ergodic measure with $h_{\mu}(f)=h$.
	Consequently, the intermediate entropy property does not hold.
\end{maintheorem}

For $\phi\in C(X)$, let $P(f,\phi)$ denote the \emph{topological pressure}.
It follows from the variational principle for topological pressure that
\begin{equation*}
	P(f,\phi)=\sup\{P_{\mu}(f,\phi):\mu\in\mathcal{M}_f(X)\},
\end{equation*}
where
\begin{equation*}
	P_{\mu}(f,\phi):=h_{\mu}(f)+\int_X\phi\mathrm{d}\mu.
\end{equation*}
Set
\begin{equation*}
	P_{\inf}(f,\phi)=\inf\{P_{\mu}(f,\phi):\mu\in\mathcal{M}_f(X)\}.
\end{equation*}
We say that $(X,f)$ satisfies the \emph{intermediate pressure property} if every $P_{\inf}(f,\phi)<P<P(f,\phi)$ equals $P_{\mu}(f,\phi)$ for some $\mu\in\mathcal{M}_f^{erg}(X)$.
When $\tau>0$, we can provide an example with the non-uniform specification property such that the intermediate entropy property holds, while the intermediate pressure property does not.

\begin{maintheorem}\label{thm-intermediate-2}
	For any $\tau\in(0,\infty]$, there exists a (positively expansive) dynamical system $(X,f)$ satisfying the non-uniform $\tau$-specification property such that the intermediate entropy property holds, while the intermediate pressure property does not.
\end{maintheorem}

Although the intermediate entropy property may fail for $\tau>0$, under positive expansiveness every interval $(\gamma,\gamma+\beta)\subset[0,h_{top}(f))$ with $\beta>\tau h_{top}(f)/(1+\tau)$ contains the entropy of an invariant subsystem.

\begin{maintheorem}\label{thm-intermediate-3}
	Suppose that $(X,f)$ satisfies the non-uniform $\tau$-specification property.
	Assume that $f$ is positively expansive.
	If $\tau<\infty$, then for any $0\leq\gamma<h_{top}(f)$ and
	\begin{equation*}
		\frac{\tau}{1+\tau}h_{top}(f)<\beta<h_{top}(f)-\gamma,
	\end{equation*}
	there is a closed $f$-invariant set $Y\subset X$ such that
	\begin{equation*}
		\gamma<h_{top}(f|_{Y})<\gamma+\beta.
	\end{equation*}
\end{maintheorem}

\begin{remark}
	The case of $\tau=0$ follows directly from \cite[Lemma 3.5]{Lin-Tian-Yu-2024} and \cite[Theorem 1.4]{Sun-2025}.
	The case $0<\tau<\infty$ is proved in Section \ref{Sect-intermediate}.
\end{remark}

\subsection{Ergodic optimization}

For $\varphi\in C(X)$, define
\begin{equation*}
	\beta(\varphi)=\sup_{\mu\in\mathcal{M}_f(X)}\int_X\varphi\mathrm{d}\mu
\end{equation*}
and
\begin{equation*}
	\mathcal{M}_{\max}(\varphi)=\bigg\{\mu\in\mathcal{M}_f(X):\int_X\varphi\mathrm{d}\mu=\beta(\varphi)\bigg\}.
\end{equation*}
The elements of $\mathcal{M}_{\max}(\varphi)$ are the maximizing measures of $\varphi$.
We refer to \cite{Bochi-2018,Jenkinson-2019} for background on ergodic optimization.
Assume that $(X,f)$ satisfies the non-uniform periodic $0$-specification property. By \cite[Theorem 3.4]{Lin-Tian-Yu-2024}, the periodic measures are dense in $\mathcal{M}_f(X)$, and hence $\overline{\mathcal{M}_f^{erg}(X)}=\mathcal{M}_f(X)$. Since periodic measures have zero entropy, the set of zero-entropy ergodic measures is dense in $\mathcal{M}_f^{erg}(X)$.
Moreover, it follows from \cite[Theorem 1.1]{Carvalho-Condori-2023} that
\begin{equation*}
	\mathcal{Z}_0=\{\mu\in\mathcal{M}_f^{erg}(X):h_\mu(f)=0\}
\end{equation*}
is a $G_\delta$ subset of $\mathcal{M}_f^{erg}(X)$, and hence a dense $G_\delta$ subset of $\mathcal{M}_f^{erg}(X)$.
Combining \cite[Corollary 1]{Morris-2010} and \cite[Theorem 3.2]{Jenkinson-2006}, we conclude that
\begin{equation*}
	\big\{\varphi\in C(X):\mathcal{M}_{\max}(\varphi)=\{\mu\}
	\text{ for some }\mu\in\mathcal{M}_f^{erg}(X)\text{ with }h_\mu(f)=0\big\}
\end{equation*}
is a dense $G_\delta$ subset of $C(X)$.
The next theorem shows that positive asymptotic linear growth of the gap function may destroy this conclusion.

\par\medskip
\begin{maintheorem}\label{thm-ergodic-optimization}
	For any $\tau\in(0,\infty]$ and $0<\delta<1-\frac{1}{1+\tau}$, there exists a dynamical system $(X,f)$ satisfying the non-uniform periodic $\tau$-specification property with gap function $M(n,\varepsilon)$ such that
	\begin{enumerate}
		\item $(X,f)$ has a unique ergodic measure $\nu_1$ of maximal entropy, and
		\begin{equation*}
			0<h_{\nu_1}(f)=h_{top}(f)<\infty;
		\end{equation*}
		\item there exists a non-empty open subset $U$ of $C(X)$ such that, for any $\varphi\in U$, we have
		\begin{equation*}
			\sup\bigg\{\int\varphi\mathrm{d}\nu:\nu\in\overline{\mathcal{M}_f^{erg}(X)\setminus\{\nu_1\}}\bigg\}
			\leq-1+\frac{1}{1+\tau}+\frac{\delta}{2}< -\frac{\delta}{2}<\int\varphi\mathrm{d}\nu_1.
		\end{equation*}
		In particular, the set
		\begin{equation*}
			\big\{\varphi\in C(X):h_\mu(f)<h_{top}(f)
			\text{ for some }\mu\in\mathcal{M}_{\max}(\varphi)\big\}
		\end{equation*}
		is not dense in $C(X)$.
	\end{enumerate}
\end{maintheorem}

\textbf{Organization of this paper.}
In Section \ref{Sect-pre}, we will introduce some basic notions and preliminary results.
In Section \ref{Sect-Bowen-lower}, we will study the lower bound estimate for Bowen topological entropy of irregular sets, and show Theorem \ref{thm-Bowen-lower-1} and Corollary \ref{cor-Bowen-lower-1}.
In Section \ref{Sect-Bowen-upper}, we will construct dynamical systems with non-uniform specification in proof of Theorem \ref{thm-construction}, show Corollary \ref{cor-entropy} and Corollary \ref{Corollary C New}, and use these examples to study the upper bound estimate for Bowen topological entropy of irregular sets, and show Theorem \ref{thm-Bowen-upper-1} and Theorem \ref{thm-Bowen-upper-2}.
In Section \ref{Sect-packing}, we will study the packing topological entropy for over-saturated sets and the completely-irregular set, and show Theorem \ref{thm-packing}.
In Section \ref{Sect-CVP}, we will study the conditional variational principle and show Theorem \ref{thm-CVP}.
In Section \ref{Sect-trans}, we will estimate the entropy for transitive points and show Theorem \ref{thm-trans-1} and Theorem \ref{thm-trans-2}.
In Section \ref{Sect-growth}, we will study the exponential periodic growth rate, and show Theorem \ref{thm-growth-1} and Theorem \ref{thm-growth-2}.
In Section \ref{Sect-intermediate}, we will study the intermediate entropy property and intermediate pressure property, and show Theorem \ref{thm-intermediate-1}, Theorem \ref{thm-intermediate-2} and Theorem \ref{thm-intermediate-3}.
In Section \ref{Sect-optimization}, we will study the ergodic optimization and prove Theorem \ref{thm-ergodic-optimization}.

\section{Preliminaries}\label{Sect-pre}

\subsection{The first Wasserstein metric}

Let $\rho$ be the first Wasserstein metric on $\mathcal{M}(X)$, which metricizes the weak$^\ast$ topology on $\mathcal{M}(X)$, see \cite{Villani-2009} for more information.
According to \cite[Page 95, (6.3)]{Villani-2009}, for any $\mu,\nu\in\mathcal{M}(X)$, $\rho(\mu,\nu)$ can be represented by
\begin{equation}\label{eq-pre-metric-1}
	\rho(\mu,\nu)=\sup_{\varphi\in\mathsf{Lip}^1(X)}\left|\int_X\varphi\mathrm{d}\mu-\int_X\varphi\mathrm{d}\nu\right|,
\end{equation}
where $\mathsf{Lip}^1(X)$ is the space of real-valued Lipschitz functions on $X$ with Lipschitz constant at most $1$.
Then for any $x,y\in X$, we have
\begin{equation}\label{eq-pre-metric-2}
	\rho(\delta_x,\delta_y)=d(x,y).
\end{equation}

The following proposition can be easily checked by using \eqref{eq-pre-metric-1}, hence we omit the proof.

\begin{proposition}\label{pre-metric-prop1}
	Let $\mu,\mu_1,\cdots,\mu_n,\nu_1,\cdots,\nu_n\in\mathcal{M}(X)$ and $0\leq s_1,\cdots,s_n\leq 1$ with $s_1+\cdots+s_n=1$, then
	\begin{equation}\label{eq-pre-metric-prop1-1}
		\rho\bigg(\mu,\sum_{i=1}^ns_i\mu_i\bigg)\leq\sum_{i=1}^ns_i\rho(\mu,\mu_i)
	\end{equation}
	and
	\begin{equation}\label{eq-pre-metric-prop1-2}
		\rho\bigg(\sum_{i=1}^ns_i\mu_i,\sum_{i=1}^ns_i\nu_i\bigg)\leq\sum_{i=1}^ns_i\rho(\mu_i,\nu_i).
	\end{equation}
\end{proposition}

\begin{lemma}\cite[Lemma 2.2]{Lin-Tian-Yu-2024}\label{pre-metric-lemma2}
	Given $0<\varepsilon,\delta\leq 1$.
	Let $\{x_i\}_{i=0}^{n-1},\{y_i\}_{i=0}^{n-1}\subset X$.
	If $|\{i\in[0,n-1]:d(x_i,y_i)\leq\varepsilon\}|\geq(1-\delta)n$, then
	\begin{equation}\label{eq-pre-metric-lemma2}
		\rho\bigg(\frac1n\sum_{i=0}^{n-1}\delta_{x_i},\frac1n\sum_{i=0}^{n-1}\delta_{y_i}\bigg)\leq\varepsilon+\delta\cdot\mathrm{diam}(X),
	\end{equation}
	where $\mathrm{diam}(X)=\sup_{x,y\in X}d(x,y)$ is the diameter of $(X,d)$.
\end{lemma}

\subsection{Over-saturated sets}

In this subsection, we will provide some basic properties of over-saturated sets for dynamical systems with the non-uniform specification property.

%lemma1
\begin{lemma}\label{pre-saturated-lemma1}
	For any non-empty subset $K\subset\mathcal{M}_f(X)$, there exists a sequence $\{\nu_j\}_{j=1}^\infty$ contained in $K$ such that
	\begin{equation}\label{eq-pre-saturated-lemma1}
		\overline{\{\nu_j:j\geq n\}}=\overline{K},\quad\forall n\in\mathbb{N}.
	\end{equation}
	Moreover, if a sequence of invariant measures $\{\nu_j\}_{j=1}^{\infty}$ satisfies \eqref{eq-pre-saturated-lemma1}, then
	\begin{equation*}
		G^K=\bigcap_{j\geq1}G^{\nu_j}.
	\end{equation*}
\end{lemma}

\begin{proof}
	From the compactness of $\overline{K}$, we can take a sequence $\{\nu_j\}_{j=1}^{\infty}\subset K$ and an increasing sequence of non-negative integers $\{m_k\}_{k=0}^{\infty}$ with $m_0=0$ such that the collection $\{B(\nu_{j},2^{-k}):j=m_{k-1}+1,\cdots,m_k\}$ covers $\overline{K}$ for all $k\in\mathbb{N}$.
	Then
	\begin{equation*}
		\overline{\{\nu_j:j\geq n\}}=\overline{K},\quad\forall n\in\mathbb{N}.
	\end{equation*}

	Since $\nu_k\in K$, we have $G^K\subset G^{\nu_k}$ for any $k\in\mathbb{N}$.
	It follows that $G^K\subset\bigcap_{j\geq 1} G^{\nu_j}$.
	Now we show that $G^K\supset\bigcap_{j\geq 1} G^{\nu_j}$.
	Fix $x\in\bigcap_{j\geq 1} G^{\nu_j}$.
	For any $j\in\mathbb{N}$, we can take an increasing sequence of integers $\{n_j(l)\}_{l=1}^{\infty}$ such that
	\begin{equation}\label{eq-pre-saturated-lemma1-proof-1}
		\lim_{l\rightarrow\infty}\delta_x^{n_j(l)}=\nu_j.
	\end{equation}
	Arbitrarily given $\mu\in K$.
	For any $\varepsilon>0$, we can find a subsequence $\{\nu_{j_i}\}_{i=1}^{\infty}\subset\{\nu_j\}_{j=1}^\infty$ such that
	\begin{equation}\label{eq-pre-saturated-lemma1-proof-2}
		\rho(\nu_{j_i},\mu)<2^{-i}\varepsilon,\quad\forall i\in\mathbb{N}.
	\end{equation}
	By \eqref{eq-pre-saturated-lemma1-proof-1}, we can choose an increasing sequence of integers $\{l_i\}_{i=1}^{\infty}$ such that
	\begin{equation}\label{eq-pre-saturated-lemma1-proof-3}
		\rho(\delta_x^{r_i},\nu_{j_i})<2^{-i}\varepsilon,\quad\forall i\in\mathbb{N},
	\end{equation}
	where $r_i=n_{j_i}(l_i)$.
	Combining \eqref{eq-pre-saturated-lemma1-proof-2} and \eqref{eq-pre-saturated-lemma1-proof-3}, we obtain that
	\begin{equation*}
		\rho(\delta_x^{r_i},\mu)\leq\rho(\delta_x^{r_i},\nu_{j_i})+\rho(\nu_{j_i},\mu)< 2^{-i+1}\varepsilon,
	\end{equation*}
	which implies that $\delta_x^{r_i}\rightarrow\mu$ and hence $x\in G^{\mu}$.
	By the arbitrariness of $\mu$, we have $x\in G^K$.
\end{proof}

%lemma2
\begin{lemma}[{\cite[Theorem 3.2 (4)]{Lin-Tian-Yu-2024}}]\label{pre-saturated-lemma2}
	Suppose that $(X,f)$ satisfies the non-uniform specification property.
	Let $K\subset\mathcal{M}_f(X)$ be a non-empty subset.
	If $G^K\neq\varnothing$, then $G^K$ is residual in $X$.
\end{lemma}

%prop3
\begin{proposition}\label{pre-saturated-prop3}
	Suppose that $(X,f)$ satisfies the non-uniform specification property.
	Let $K\subset\mathcal{M}_f(X)$ be a non-empty subset and $\{\nu_j\}_{j=1}^\infty$ be a sequence contained in $K$ satisfying
	\begin{equation*}
		\overline{\{\nu_j:j\geq n\}}=\overline{K},\quad\forall n\in\mathbb{N}.
	\end{equation*}
	Then $G^K\neq\varnothing$ if and only if $G^{\nu_j}\neq\varnothing$ for all $j\in\mathbb{N}$.
\end{proposition}

\begin{proof}
	The necessity follows from Lemma \ref{pre-saturated-lemma1}.
	For the sufficiency, by Lemma \ref{pre-saturated-lemma2}, every $G^{\nu_j}$ is residual in $X$.
	Using Lemma \ref{pre-saturated-lemma1}, we conclude that $G^K=\bigcap_{j\geq1}G^{\nu_j}$ is residual in $X$.
\end{proof}

%cor4
\begin{corollary}\label{pre-saturated-cor4}
	Suppose that $(X,f)$ satisfies the non-uniform specification property.
	Then
	\begin{enumerate}
		\item $G^{\mathcal{M}_f^{erg}(X)}=G^{\overline{\mathcal{M}_f^{erg}(X)}}\neq\varnothing$;
		\item for any $0<\sigma<1$, we have $C_{\sigma}(X)\subset\mathfrak{C}(f)$.
	\end{enumerate}
\end{corollary}

\begin{proof}
	Statement (1) follows directly from Proposition \ref{pre-saturated-prop3}.
	We proceed to show (2).
	It suffices to show that for any $\varphi\in C(X)$, if
	\begin{equation}\label{eq-pre-saturated-cor4-proof-1}
		\sup_{\mu\in\mathcal{M}_f(X)}\int_X\varphi\mathrm{d}\mu>\inf_{\mu\in\mathcal{M}_f(X)}\int_X\varphi\mathrm{d}\mu,
	\end{equation}
	then $\varphi\in\mathfrak{C}(f)$.
	Note that \eqref{eq-pre-saturated-cor4-proof-1} is equivalent to
	\begin{equation*}
		\sup_{\mu\in\mathcal{M}_f^{erg}(X)}\int_X\varphi\mathrm{d}\mu>\inf_{\mu\in\mathcal{M}_f^{erg}(X)}\int_X\varphi\mathrm{d}\mu.
	\end{equation*}
	Hence we can take ergodic measures $\mu$ and $\nu$ such that
	\begin{equation}\label{eq-pre-saturated-cor4-proof-2}
		\int_X\varphi\mathrm{d}\mu>\int_X\varphi\mathrm{d}\nu.
	\end{equation}
	By statement (1), we can fix an $x\in G^{\mathcal{M}_f^{erg}(X)}$.
	Then we can take increasing sequences $\{n_k\}_{k=1}^{\infty}$ and $\{m_k\}_{k=1}^{\infty}$ such that $\delta_x^{n_k}\rightarrow\mu$ and $\delta_x^{m_k}\rightarrow\nu$ as $k\rightarrow\infty$.
	From \eqref{eq-pre-saturated-cor4-proof-2}, it follows that
	\begin{equation*}
		\lim_{k\rightarrow\infty}\frac{1}{n_k}\sum_{i=0}^{n_k-1}\varphi(f^ix)=\int_X\varphi\mathrm{d}\mu>\int_X\varphi\mathrm{d}\nu=\lim_{k\rightarrow\infty}\frac{1}{m_k}\sum_{i=0}^{m_k-1}\varphi(f^ix),
	\end{equation*}
	which implies that $\varphi\in\mathfrak{C}(f)$.
	This completes the proof of Corollary \ref{pre-saturated-cor4}.
\end{proof}

%prop5
\begin{proposition}\label{pre-saturated-prop5}
	Suppose that $(X,f)$ satisfies the non-uniform specification property.
	Let
	\begin{equation}\label{eq-pre-saturated-prop5}
		\mathcal{K}=\bigcup_{x\in X}V_f(x).
	\end{equation}
	Then
	\begin{enumerate}
		\item $G^K\neq\varnothing$ if and only if $K\in\mathcal{P}_*(\mathcal{K})$, where $\mathcal{P}_*(\mathcal{K})$ is the collection of non-empty subsets of $\mathcal{K}$;
		\item if $G^K\neq\varnothing$, then $G^{\mathcal{K}}\subset G^K$;
		\item $G_{\mathcal{K}}=G^{\mathcal{K}}\neq\varnothing$, and consequently $\mathcal{K}$ is a non-empty compact connected set.
	\end{enumerate}
\end{proposition}

\begin{proof}
	If $G^K\neq\varnothing$, then $K\subset V_f(x)\subset\mathcal{K}$ for any $x\in G^K$.
	For the sufficiency, using Proposition \ref{pre-saturated-prop3}, we obtain that $G^{\mathcal{K}}\neq\varnothing$, which implies that for any $K\in\mathcal{P}_*(\mathcal{K})$, $G^K\neq\varnothing$ since $G^K\supset G^{\mathcal{K}}$.
	This completes the proof of statement (1).
	The statement (2) follows immediately from (1).

	Now we show statement (3).
	Clearly $G_{\mathcal{K}}\subset G^{\mathcal{K}}$.
	On the other hand, for any $y\in G^{\mathcal{K}}$, we have
	\begin{equation*}
		\mathcal{K}\subset V_f(y)\subset\bigcup_{x\in X}V_f(x)=\mathcal{K},
	\end{equation*}
	which implies that $V_f(y)=\mathcal{K}$.
	We conclude that $G_{\mathcal{K}}=G^{\mathcal{K}}\neq\varnothing$.
	It follows from \cite[Proposition 3.8]{DGS1976} that $\mathcal{K}$ is a non-empty compact connected set since $\mathcal{K}=V_f(x)$ for any $x\in G_{\mathcal{K}}$.
\end{proof}

%lemma6
\begin{lemma}\label{pre-saturated-lemma6}
	Suppose that $(X,f)$ satisfies the non-uniform specification property.
	Then
	\begin{enumerate}
		\item for any $\varphi\in\mathfrak{C}(f)$, there exists $K\subset\mathcal{M}_f(X)$ such that $G^K\neq\varnothing$ and $G^K\subset I_{\varphi}(f)$;
		\item $G^{\mathcal{K}}\subset\mathrm{CI}(f)$.
	\end{enumerate}
\end{lemma}

\begin{proof}
	By the definition of $\mathfrak{C}(f)$, we have $I_{\varphi}(f)\neq\varnothing$.
	Hence it follows from the ergodic decomposition theorem that there exist two different ergodic measures $\mu,\nu\in\mathcal{M}_f^{erg}(X)$ such that
	\begin{equation}\label{eq-pre-saturated-lemma6-proof-1}
		\int_X\varphi\mathrm{d}\mu<\int_X\varphi\mathrm{d}\nu.
	\end{equation}
	Let $K=\{\mu,\nu\}$.
	Then by Proposition \ref{pre-saturated-prop5}, we have $G^K\neq\varnothing$.
	We aim to show that $G^K\subset I_{\varphi}(f)$.
	Fix any $x\in G^K$.
	Since $K\subset V_f(x)$, we can choose sequences $\{m_k\}_{k=1}^{\infty}$ and $\{n_k\}_{k=1}^{\infty}$ of $\mathbb{N}$ such that
	\begin{equation}\label{eq-pre-saturated-lemma6-proof-2}
		\lim_{k\rightarrow\infty}\delta_x^{m_k}=\mu\quad\text{and}\quad\lim_{k\rightarrow\infty}\delta_x^{n_k}=\nu.
	\end{equation}
	Combining \eqref{eq-pre-saturated-lemma6-proof-1} and \eqref{eq-pre-saturated-lemma6-proof-2}, we obtain that
	\begin{align*}
		\lim_{k\rightarrow\infty}\frac{1}{m_k}\sum_{i=0}^{m_k-1}\varphi(f^ix)=\lim_{k\rightarrow\infty}\int_X\varphi\mathrm{d}\delta_x^{m_k}&=\int_X\varphi\mathrm{d}\mu \\
		&<\int_X\varphi\mathrm{d}\nu=\lim_{k\rightarrow\infty}\int_X\varphi\mathrm{d}\delta_x^{n_k}=\lim_{k\rightarrow\infty}\frac{1}{n_k}\sum_{i=0}^{n_k-1}\varphi(f^ix).
	\end{align*}
	Therefore $x\in I_{\varphi}(f)$.
	As a result, $G^K\subset I_{\varphi}(f)$, which completes the proof of statement (1).

	Now we show the statement (2).
	By Proposition \ref{pre-saturated-prop5}, we have $G^{\mathcal{K}}\subset G^K$.
	For any $\varphi\in\mathfrak{C}(f)$, it follows from statement (1) that $G^{\mathcal{K}}\subset I_{\varphi}(f)$.
	Hence $G^{\mathcal{K}}\subset\mathrm{CI}(f)$.
\end{proof}

\subsection{Separated sets}

For $x,y\in X$, $\varepsilon>0$ and $n\in\mathbb{N}$, define
\begin{equation*}
	d_n(x,y)=\max\{d(f^ix,f^iy):0\leq i\leq n-1\}.
\end{equation*}
Let $B_n(x,\varepsilon):=\{y\in X:d_n(x,y)<\varepsilon\}$ and $\overline{B_n(x,\varepsilon)}:=\{y\in X:d_n(x,y)\leq\varepsilon\}$ denote the open and closed Bowen balls, respectively.

For $\varepsilon>0$ and $n\in\mathbb{N}$, two points $x,y\in X$ are called \emph{$(n,\varepsilon)$-separated} if
\begin{equation*}
	|\{0\leq i\leq n-1:d(f^ix,f^iy)>\varepsilon\}|\geq 1.
\end{equation*}
A subset $E\subset X$ is called \emph{$(n,\varepsilon)$-separated} if each pair of distinct points in $E$ are $(n,\varepsilon)$-separated.
It is known that
\begin{equation*}
	\begin{split}
		h_{top}(f)&=\lim_{\varepsilon\rightarrow0}\limsup_{n\rightarrow\infty}\frac1n\log s_n(\varepsilon)\\
		&=\lim_{\varepsilon\rightarrow0}\liminf_{n\rightarrow\infty}\frac1n\log s_n(\varepsilon),
	\end{split}
\end{equation*}
where $s_n(\varepsilon)$ denotes the maximal cardinality of $(n,\varepsilon)$-separated sets.
It is known that, if $(X,f)$ is positively expansive with expansive constant $\varepsilon_0$, then for any $0<\varepsilon<\varepsilon_0$, we have
\begin{equation*}
	h_{top}(f)=\lim_{n\rightarrow\infty}\frac{1}{n}\log s_n(\varepsilon).
\end{equation*}

For $\delta,\varepsilon>0$ and $n\in\mathbb{N}$, a pair of points $x,y\in X$ is called \emph{$(\delta,n,\varepsilon)$-separated} if
\begin{equation*}
	|\{0\leq i\leq n-1:d(f^ix,f^iy)>\varepsilon\}|\geq\delta n.
\end{equation*}
A subset $E\subset X$ is called \emph{$(\delta,n,\varepsilon)$-separated} if each pair of distinct points in $E$ are $(\delta,n,\varepsilon)$-separated.
Clearly, each $(\delta,n,\varepsilon)$-separated set is an $(n,\varepsilon)$-separated set.
Let $F$ be a neighborhood of $\nu\in\mathcal{M}_f(X)$ in $\mathcal{M}(X)$.
Define
\begin{equation*}
	X_{n,F}=\bigg\{x\in X:\frac1n\sum_{i=0}^{n-1}\delta_{f^ix}\in F\bigg\}.
\end{equation*}

%prop1
\begin{proposition}[{\cite[Proposition 2.1]{Pfister-Sullivan-2005}}]\label{pre-separation-prop1}
	Let $\mu\in\mathcal{M}_f^{erg}(X)$ be an ergodic measure and $0<\gamma<h_{\mu}(f)$.
	Then there exist $\delta^*>0$ and $\varepsilon^*>0$ such that for any neighborhood $F$ of $\mu$ in $\mathcal{M}(X)$, there exists $n^*_F\in\mathbb{N}$ with the following property:
	for any $n\geq n^*_F$, there is a $(\delta^*,n,\varepsilon^*)$-separated set $\Gamma^{\mu}_n\subset X_{n,F}$ satisfying
	\begin{equation*}
		|\Gamma^{\mu}_n|\geq e^{n\gamma}.
	\end{equation*}
\end{proposition}

For $\mu\in\mathcal{M}_f(X)$, define
\begin{equation*}
	\underline{s}(\mu,\delta,\varepsilon)=\inf_{F\ni\mu}\liminf_{n\rightarrow\infty}\frac1n\log N(F;\delta,n,\varepsilon),
\end{equation*}
where the infimum is taken over any base of neighborhoods of $\mu$ in $\mathcal{M}(X)$ and $N(F;\delta,n,\varepsilon)$ denotes the maximal cardinality of $(\delta,n,\varepsilon)$-separated subset of $X_{n,F}$.
Proposition \ref{pre-separation-prop1} indicates that for any ergodic measure $\mu\in\mathcal{M}_f^{erg}(X)$ and any $0<\gamma<h_{\mu}(f)$, there exist $\delta^*>0$ and $\varepsilon^*>0$ such that
\begin{equation*}
	\underline{s}(\mu,\delta^*,\varepsilon^*)\geq\gamma.
\end{equation*}
For invariant measures, we have the following lemma.

%lemma2
\begin{lemma}[{\cite[Lemma 6.2]{Pfister-Sullivan-2007}}]\label{pre-separation-lemma2}
	Suppose that $\mu\in\mathcal{M}_f(X)$ is an invariant measure and $0<\gamma<h_{\mu}(f)$.
	Then there exist $\delta^*>0$, $\varepsilon^*>0$ such that for any $j\in\mathbb{N}$, there is a finite convex combination of ergodic measures with rational coefficients
	\begin{equation*}
		\mu_j=\sum_{u=1}^{m_j}p_{j,u}\mu_{j,u}
	\end{equation*}
	satisfying $\mu_j\rightarrow\mu$ as $j\rightarrow\infty$ and
	\begin{equation*}
		\sum_{u=1}^{m_j}p_{j,u}\cdot\underline{s}(\mu_{j,u},\delta^*,\varepsilon^*)>\gamma.
	\end{equation*}
\end{lemma}

\subsection{Topological entropies for subsets}

\subsubsection{Upper capacity topological entropy}

Fix a subset $E\subset X$.
For $\varepsilon>0$, let $s_n(\varepsilon,E)$ denote the maximal cardinality of $(n,\varepsilon)$-separated sets contained in $E$.
Define the \emph{upper capacity topological entropy} of $E$ by
\begin{equation*}
	h_{top}^{UC}(f,E)=\lim_{\varepsilon\to0}\limsup_{n\to\infty}\frac1n\log s_n(\varepsilon,E).
\end{equation*}
It is known that if $f$ is positively expansive with expansive constant $\varepsilon_0$, then for any $0<\delta<\varepsilon_0$,
\begin{equation}\label{eq-pre-entropy-UC-1-expansive}
	h_{top}^{UC}(f,E)=\limsup_{n\to\infty}\frac{1}{n}\log s_n(\delta,E).
\end{equation}

\subsubsection{Bowen topological entropy}

Fix a subset $F\subset X$.
For $\varepsilon>0$ and $N\in\mathbb{N}$, let $\mathcal{F}_N(F,\varepsilon)$ be the collection of all finite or countable families $\{B_{n_k}(x_k,\varepsilon)\}$ such that $x_k\in F$, $n_k\geq N$ and $\{B_{n_k}(x_k,\varepsilon)\}$ covers $F$.
For $t\geq 0$, let
\begin{equation*}
	C(F;t,N,\varepsilon,f)=\inf_{\mathcal{C}\in\mathcal{F}_N(F,\varepsilon)}\sum_{B_n(x,\varepsilon)\in\mathcal{C}}e^{-tn}
\end{equation*}
and
\begin{equation*}
	C(F;t,\varepsilon,f)=\lim_{N\rightarrow\infty}C(F;t,N,\varepsilon,f).
\end{equation*}
It is easy to see that the quantity $C(F;s,\varepsilon,f)$ has a critical value of parameter $s$ jumping from $\infty$ to $0$,
which is denoted by
\begin{equation*}
	h_{top}^B(F;\varepsilon,f)=\inf\{s\geq 0:C(F;s,\varepsilon,f)=0\}=\sup\{s\geq 0:C(F;s,\varepsilon,f)=\infty\}.
\end{equation*}
Note that $h_{top}^B(F;\varepsilon,f)$ increases as $\varepsilon$ decreases.
Define the \emph{Bowen topological entropy} of $F$ by
\begin{equation*}
	h_{top}^B(f,F)=\lim_{\varepsilon\rightarrow0}h_{top}^B(F;\varepsilon,f).
\end{equation*}

\subsubsection{Packing topological entropy}
Fix a subset $E\subset X$.
For $\varepsilon>0$ and $N\in\mathbb{N}$, let $\mathcal{G}_N(E,\varepsilon)$ be the collection of all finite or countable pairwise disjoint families $\Big\{\overline{B_{n_k}(x_k,\varepsilon)}\Big\}$ such that $x_k\in E$ and $n_k\geq N$.
For $t\geq 0$, let
\begin{equation*}
	Q(E;t,N,\varepsilon,f)=\sup_{\mathcal{D}\in\mathcal{G}_N(E,\varepsilon)}\sum_{\overline{B_n(x,\varepsilon)}\in\mathcal{D}}e^{-tn}
\end{equation*}
and
\begin{equation*}
	Q(E;t,\varepsilon,f)=\lim_{N\rightarrow\infty}Q(E;t,N,\varepsilon,f).
\end{equation*}
Define
\begin{equation*}
	P(E;t,\varepsilon,f)=\inf\bigg\{\sum_{j=1}^{\infty}Q(E_j;t,\varepsilon,f):\bigcup_{j=1}^{\infty}E_j\supset E\bigg\}.
\end{equation*}
It is easy to see that the quantity $P(E;s,\varepsilon,f)$ has a critical value of parameter $s$ jumping from $\infty$ to $0$,
which is denoted by
\begin{equation*}
	h_{top}^P(E;\varepsilon,f)=\inf\{s\geq 0:P(E;s,\varepsilon,f)=0\}=\sup\{s\geq 0:P(E;s,\varepsilon,f)=\infty\}.
\end{equation*}
Note that $h_{top}^P(E;\varepsilon,f)$ increases as $\varepsilon$ decreases.
Define the \emph{packing topological entropy} of $E$ by
\begin{equation*}
	h_{top}^P(f,E)=\lim_{\varepsilon\rightarrow0}h_{top}^P(E;\varepsilon,f).
\end{equation*}

\subsubsection{Some basic properties}

%prop1
\begin{proposition}[{\cite[Proposition 2.1]{Feng-Huang-2012}}]\label{pre-entropy-prop1}
	Let $(X,f)$ be a dynamical system, then the following statements hold:
	\begin{enumerate}
		\item for $E\subset F\subset X$, we have
		\begin{equation*}
			h_{top}^{UC}(f,E)\leq h_{top}^{UC}(f,F),\quad h_{top}^B(f,E)\leq h_{top}^B(f,F)\quad\text{and}\quad h_{top}^P(f,E)\leq h_{top}^P(f,F);
		\end{equation*}
		\item for $E\subset\bigcup_{k=1}^{\infty}E_k$, we have
		\begin{equation*}
			h_{top}^B(f,E)\leq \sup_{k\geq1}h_{top}^B(f,E_k)\quad\text{and}\quad h_{top}^P(f,E)\leq\sup_{k\geq1}h_{top}^P(f,E_k);
		\end{equation*}
		\item for any $E\subset X$, we have
		\begin{equation*}
			h_{top}^B(f,E)\leq h_{top}^P(f,E)\leq h_{top}^{UC}(f,E);
		\end{equation*}
		\item if $E\subset X$ is compact and $f$-invariant, then
		\begin{equation*}
			h_{top}^B(f,E)=h_{top}^P(f,E)=h_{top}^{UC}(f,E)=h_{top}(f|_E).
		\end{equation*}
	\end{enumerate}
\end{proposition}

\subsubsection{Variational principles}

In \cite{Feng-Huang-2012}, Feng and Huang introduced the notion of measure-theoretical lower and upper entropies to study the variational principles for topological entropies on subsets.
\begin{definition}
	Let $\mu\in\mathcal{M}(X)$.
	The \emph{measure-theoretical lower entropy} and \emph{measure-theoretical upper entropy} of $\mu$ are defined by
	\begin{equation*}
		\underline{h}_{\mu}(f)=\int_X\underline{h}_{\mu}(f,x)\mathrm{d}\mu(x)
		\quad\text{and}\quad
		\overline{h}_{\mu}(f)=\int_X\overline{h}_{\mu}(f,x)\mathrm{d}\mu(x),
	\end{equation*}
	respectively, where
	\begin{equation*}
		\begin{split}
			\underline{h}_{\mu}(f,x)&=\lim_{\varepsilon\rightarrow0}\liminf_{n\rightarrow\infty}-\frac1n\log\mu(B_n(x,\varepsilon)),\\
			\overline{h}_{\mu}(f,x)&=\lim_{\varepsilon\rightarrow0}\limsup_{n\rightarrow\infty}-\frac1n\log\mu(B_n(x,\varepsilon)).
		\end{split}
	\end{equation*}
\end{definition}

%lemma2
\begin{lemma}[Variational principles, {\cite[Theorem 1.2 and 1.3]{Feng-Huang-2012}}]\label{pre-entropy-lemma2}
	Let $(X,f)$ be a dynamical system.
	If $E\subset X$ is non-empty and compact, then
	\begin{equation*}
		\begin{split}
			h_{top}^B(f,E)&=\sup\{\underline{h}_{\mu}(f):\mu\in\mathcal{M}(X),\mu(E)=1\},\\
			h_{top}^P(f,E)&=\sup\{\overline{h}_{\mu}(f):\mu\in\mathcal{M}(X),\mu(E)=1\}.
		\end{split}
	\end{equation*}
\end{lemma}

\subsection{Symbolic dynamics}
Let $\mathcal{A}$ be a finite alphabet and $\mathcal{A}^{\mathbb{N}_0}:=\{x_0x_1x_2\cdots:x_i\in\mathcal{A}\;\text{for}\;i\in\mathbb{N}_0\}$.
Each element in $\mathcal{A}^{\mathbb{N}_0}$ is also called a sequence over $\mathcal{A}$.
The one-sided full shift over $\mathcal{A}$ is defined by $(\mathcal{A}^{\mathbb{N}_0},\sigma)$, where
\begin{equation*}
	\sigma(x_0x_1x_2\cdots)=x_1x_2\cdots,\quad\forall x_0x_1x_2\cdots\in \mathcal{A}^{\mathbb{N}_0}.
\end{equation*}
Without causing ambiguity, we also use $f$ to denote $\sigma$ sometimes. The metric on $(\mathcal{A}^{\mathbb{N}_0},f)$ is defined as follows: for $x=(x_0x_1x_2\cdots)$ and $y=(y_0y_1y_2\cdots)$,
\begin{equation*}
	d(x,y)=\sum_{k\geq0}\frac{1}{2^k}d'(x_k,y_k),
\end{equation*}
where $d'$ is the discrete metric on $\mathcal{A}$.

Given a finite alphabet $\mathcal{A}$.
Each member $w$ of $\mathcal{A}^{\{i,i+1,\cdots,j\}}$ for some $0\leq i\leq j$ is called a \emph{word} over $\mathcal{A}$, and its length $j-i+1$ is denoted by $|w|$.
Let $\mathcal{A}^\ast:=\bigcup_{0\leq i\leq j}\mathcal{A}^{\{i,i+1,\cdots,j\}}$ denote the set of all words over $\mathcal{A}$.
For $n\in\mathbb{N}$, denote $\mathcal{A}^n:=\mathcal{A}^{\{0,\cdots,n-1\}}$.
A word $w$ is said to be \emph{empty} if $w\in\mathcal{A}^{\varnothing}$ has length $0$.

Let $X\subset\mathcal{A}^{\mathbb{N}_0}$ be a non-empty closed $f$-invariant set, then $(X,f)$ is a subsystem of $(\mathcal{A}^{\mathbb{N}_0},f)$, called a \emph{subshift}.
A standard way to define a subshift is by using forbidden words.
Given a subset $\mathcal{F}\subset \mathcal{A}^\ast$, we can define
\begin{equation*}
	X=X(\mathcal{F}):=\{x\in \mathcal{A}^{\mathbb{N}_0}:x_ix_{i+1}\cdots x_j\notin\mathcal{F}\;\text{for all}\;i,j\in\mathbb{N}_0\;\text{with}\;0\leq i\leq j\}.
\end{equation*}
Then $(X,f)$ is a subshift of $(\mathcal{A}^{\mathbb{N}_0},f)$.
The elements of $\mathcal{F}$ are called the \emph{forbidden words}.
Given $w\in\mathcal{A}^\ast$ and $n\in\mathbb{N}_{0}\cup\{\infty\}$, we define $w^n:=\underbrace{ww\cdots w}_\text{$n$ items}$.

The \emph{language} of a subshift $X$, denoted by $\mathcal{L}(X)$, is defined as the set of all words which appear in points of $X$.
For any $n\in\mathbb{N}$, denote $\mathcal{L}_n(X):=\mathcal{L}(X)\cap \mathcal{A}^n$.
For $w\in\mathcal{L}_n(X)$, the \emph{cylinder set} $[w]$ is defined as the set of all $x\in X$ with $x_0x_1\cdots x_{n-1}=w$.
Note that every cylinder set $[w]$ is both open and closed, hence its characteristic function $\chi_{[w]}$ is continuous.
It is known that the topological entropy of a subshift $(X,f)$ can be computed as
\begin{equation*}
	h_{top}(f)=h_{top}(X):=\lim_{n\rightarrow\infty}\frac{1}{n}\log|\mathcal{L}_n(X)|.
\end{equation*}
We may use $h_{top}(X)$ to denote the topological entropy of a subshift $(X,f)$, rather than $h_{top}(f)$, without ambiguity.

Given a subshift $(X,f)$, an ergodic measure $\mu\in\mathcal{M}_f^{erg}(X)$, a word $w\in\mathcal{L}(X)$, an integer $n\in\mathbb{N}$, and $\varepsilon>0$, let $\mathcal{W}_{n,\varepsilon,w}(X)$ be the set of all $u\in\mathcal{L}_n(X)$ having between $n(\mu([w])-\varepsilon)$ and $n(\mu([w])+\varepsilon)$ occurrences of $w$.
The following lemma, which is just an application of Birkhoff's ergodic theorem, is used to estimate the measure-theoretic entropy of ergodic measures of subshifts, where the proof (of a more general version) can be found as \cite[Lemma 4.8]{Pavlov-2011}.

%lemma1
\begin{lemma}[{\cite[Lemma 2.12]{Pavlov-2016}}]\label{pre-symbolic-lemma1}
	Let $(X,f)$ be a subshift and $\mu\in\mathcal{M}_f^{erg}(X)$ an ergodic measure.
	Given $w\in\mathcal{L}(X)$ and $0<\varepsilon<1$, for any $n\in\mathbb{N}$, let $\mathcal{W}_{n,\varepsilon,w}(X)$ be as above.
	Then
	\begin{equation*}
		h_{\mu}(f)\leq\liminf_{n\rightarrow\infty}\frac{1}{n}\log|\mathcal{W}_{n,\varepsilon,w}(X)|.
	\end{equation*}
\end{lemma}

Recall that a dynamical system $(X,f)$ is said to be \emph{strictly ergodic} if the following two conditions hold:
\begin{itemize}
	\item $(X,f)$ is \emph{minimal}, i.e. $X$ has no non-trivial closed $f$-invariant subset;
	\item $(X,f)$ is \emph{uniquely ergodic}, i.e. $\mathcal{M}_f^{erg}(X)$ is a singleton set.
\end{itemize}

Hahn and Katznelson constructed strictly ergodic systems with arbitrarily large entropy in \cite{Hahn-Katznelson-1967}.
We use the following version from \cite{Damanik-2007}.
%lemma2
\begin{lemma}[{\cite[Theorem 2.22]{Damanik-2007}}]\label{pre-symbolic-lemma2}
	For any $\delta>0$, there are alphabets $\mathcal{A}^{(j)}$ and sequences $z^{(j)}$ over $\mathcal{A}^{(j)}$ for all $j\in\mathbb{N}$ such that
	\begin{enumerate}
		\item $|\mathcal{A}^{(j)}|\rightarrow\infty$ as $j\rightarrow\infty$;
		\item for any $j\in\mathbb{N}$, the subshift $(Z^{(j)},f)$ is strictly ergodic, where $Z^{(j)}=\omega(z^{(j)},f)$ is the $\omega$-limit set of $z^{(j)}$;
		\item for any $j\in\mathbb{N}$, we have
		\begin{equation*}
			-\delta+\log|\mathcal{A}^{(j)}|<h_{top}(Z^{(j)})<\log|\mathcal{A}^{(j)}|.
		\end{equation*}
	\end{enumerate}
\end{lemma}

\subsection{Relative non-uniform specification property}\label{Section 2.6}
For $n\in\mathbb{N}$, a collection $\mathcal{D}=\{D_0,D_1,\cdots,D_{n-1}\}$ is said to be a \emph{regular periodic decomposition} of $(X,f)$ if the following holds:
\begin{itemize}
	\item $D_i=\overline{int(D_i)}$ for any $0\leq i\leq n-1$;
	\item $D_i\cap D_j$ is nowhere dense for any $0\leq i\neq j\leq n-1$;
	\item $f(D_i)\subset D_{i+1\text{ } (\text{mod } n)}$ for any $0\leq i\leq n-1$;
	\item $\bigcup_{i=0}^{n-1}D_i=X$.
\end{itemize}
If further $D_i\cap D_j=\varnothing$ for any $0\leq i\neq j\leq n-1$, then we say that $\mathcal{D}$ is a \emph{regular periodic decomposition with disjoint steps}.

\begin{lemma}\cite[Lemma 2.1]{Banks1997}\label{Lem-regular-periodic}
	Suppose that $(X,f)$ is topologically transitive and has a regular periodic decomposition $\mathcal{D}=\{D_0,D_1,\cdots,D_{n-1}\}$, then $D_0=f^n(D_0)$ and $D_j=f^j(D_0)$ for any $0\leq j\leq n-1$.
\end{lemma}
\begin{definition}
	A dynamical system $(X,f)$ is said to satisfy the \emph{relative non-uniform ($\tau$-)specification property (with disjoint steps)} if there exists a regular periodic decomposition $$\mathcal{D}=\{D_0,D_1,\cdots,D_{n-1}\}$$ (with disjoint steps) such that $(D_i,f^n)$ has non-uniform ($\tau$-)specification property for any $0\leq i\leq n-1$.
\end{definition}
\begin{remark}
	\begin{enumerate}
		\item It is clear that if $f$ is surjective and has the relative non-uniform specification property, then $f$ is topologically transitive. Hence, by Lemma \ref{Lem-regular-periodic}, $D_0=f^n(D_0)$ and $D_j=f^j(D_0)$ for any $0\leq j\leq n-1$.
		\item Given $N\in\mathbb{N}$ and $\varepsilon>0$, choose $0<\delta(\varepsilon)\leq\varepsilon$ such that for any $x,y\in X$ with $d(x,y)\leq\delta(\varepsilon)$, one has $d(f^ix,f^iy)\leq\varepsilon$ for any $0\leq i\leq N-1$. Recall that the gap function $M(n,\varepsilon)$ is non-increasing with $\varepsilon$, hence, we always have $$\sup_{\varepsilon>0}\lim_{n\rightarrow\infty}\frac{M(n,\delta(\varepsilon))}{n}=\sup_{\varepsilon>0}\lim_{n\rightarrow\infty}\frac{M(n,\varepsilon)}{n}\text{ and }\sup_{\varepsilon>0}\liminf_{n\rightarrow\infty}\frac{M(n,\delta(\varepsilon))}{n}=\sup_{\varepsilon>0}\liminf_{n\rightarrow\infty}\frac{M(n,\varepsilon)}{n}$$
		If $(D,f^N)$ has non-uniform ($\tau$-)specification property with gap function $M(n,\varepsilon)$, then for any $0\leq i\leq N-1$, $(f^i(D),f^N)$  has non-uniform ($\tau$-)specification property with gap function $M(n,\delta(\varepsilon))$. Hence, we can always assume that $(f^i(D),f^N)$  has non-uniform ($\tau$-)specification property with the same gap function.
	\end{enumerate}
\end{remark}

\section{Bowen topological entropy: lower bound estimate}\label{Sect-Bowen-lower}

In this section, we will prove Theorem \ref{thm-Bowen-lower-1} and Corollary \ref{cor-Bowen-lower-1}.
The following lemma is used to estimate the lower bound for Bowen topological entropy.

\begin{lemma}\label{Bowen-lemma1}
	Let $\{n_i\}_{i=1}^{\infty}$ and $\{t_i\}_{i=1}^{\infty}$ be two sequences of positive integers, and let $\{M_i\}_{i=1}^{\infty}$ be a sequence of non-negative integers.
	Define an array of non-negative integers $\{s(j,k):j\in\mathbb{N},0\leq k\leq n_j-1\}$ by
	\begin{equation*}
		s(j,k)=
		\begin{cases}
			kt_1 & \text{if}\quad j=1,\\
			\sum_{i=1}^{j-1}(n_it_i+M_i)+kt_j & \text{if}\quad j\geq2.
		\end{cases}
	\end{equation*}
	Given $\alpha,\beta>0$.
	Assume that
	\begin{equation}\label{eq-Bowen-lemma1-1}
		\liminf_{j\rightarrow\infty}\frac{\sum_{i=1}^jn_it_i}{s(j+1,1)-1}\geq\alpha.
	\end{equation}
	Let $\{\Omega_j(i):j\in\mathbb{N},0\leq i\leq n_j-1\}$ be the index sets such that for any $j\in\mathbb{N}$, we have
	\begin{equation}\label{eq-Bowen-lemma1-2}
		\exp(\beta\cdot t_j)\leq|\Omega_j(i)|<\infty.
	\end{equation}
	Let $\varepsilon>0$.
	Suppose that
	\begin{equation*}
		Y=\bigg\{x_{\xi}\in X:\xi=(\xi_1^0,\cdots,\xi_1^{n_1-1},\xi_2^0,\cdots,\xi_2^{n_2-1},\cdots)\in\prod_{j=1}^{\infty}\prod_{i=0}^{n_j-1}\Omega_j(i)\bigg\}\subset X
	\end{equation*}
	satisfies the following condition: $f^{s(j,k)}x_{\xi}$ and $f^{s(j,k)}x_{\zeta}$ are $(t_j,\varepsilon)$-separated whenever $\xi_j^k\neq\zeta_j^k$ for some $j\in\mathbb{N}$ and $0\leq k\leq n_j-1$.
	Then
	\begin{equation*}
		h_{top}^B(f,\overline{Y})\geq\alpha\beta.
	\end{equation*}
\end{lemma}

\begin{proof}
	Fix $\zeta=(\zeta_1^0,\cdots,\zeta_1^{n_1-1},\zeta_2^0,\cdots,\zeta_2^{n_2-1},\cdots)\in\prod_{j=1}^{\infty}\prod_{i=0}^{n_j-1}\Omega_j(i)$.
	For any $l\in\mathbb{N}$, let
	\begin{equation*}
		\Delta_l=\bigg\{\xi\in\prod_{j=1}^{\infty}\prod_{i=0}^{n_j-1}\Omega_j(i):\xi_j^i=\zeta_j^i,\forall j\geq l+1,0\leq i\leq n_j-1\bigg\}
	\end{equation*}
	and
	\begin{equation*}
		\mu_l=\frac{1}{|\Delta_l|}\sum_{\xi\in\Delta_l}\delta_{x_{\xi}}=\frac{1}{\prod_{j=1}^l\prod_{i=0}^{n_j-1}|\Omega_j(i)|}\sum_{\xi\in\Delta_l}\delta_{x_{\xi}}\in\mathcal{M}(X).
	\end{equation*}
	Then $\mu_l(\overline{Y})=1$.
	Let $\mu\in\mathcal{M}(X)$ be a limit point of the sequence $\{\mu_l\}_{l=1}^{\infty}$ and $\{l_r\}_{r=1}^{\infty}$ be an increasing subsequence of positive integers such that $\mu_{l_r}\rightarrow\mu$ weakly as $r\rightarrow\infty$.
	Then
	\begin{equation*}
		\mu(\overline{Y})\geq\limsup_{r\rightarrow\infty}\mu_{l_r}(\overline{Y})=1.
	\end{equation*}
	Using Lemma \ref{pre-entropy-lemma2}, we obtain that
	\begin{equation}\label{eq-Bowen-lemma1-proof-1}
		h_{top}^B(f,\overline{Y})\geq\underline{h}_{\mu}(f)=\int_{\overline{Y}}\underline{h}_{\mu}(f,x)\mathrm{d}\mu(x)\geq\int_{\overline{Y}}\liminf_{n\rightarrow\infty}-\frac{1}{n}\log\mu(B_n(x,\varepsilon/2))\mathrm{d}\mu(x).
	\end{equation}
	Fix $x\in\overline{Y}$ and $n\in\mathbb{N}$, then there exists $y\in Y$ such that $d_n(x,y)<\varepsilon/2$, and hence $B_n(x,\varepsilon/2)\subset B_n(y,\varepsilon)$.
	Thus for any $r\in\mathbb{N}$, we have
	\begin{equation}\label{eq-Bowen-lemma1-proof-2}
		\mu_{l_r}(B_n(x,\varepsilon/2))\leq\mu_{l_r}(B_n(y,\varepsilon))=\mu_{l_r}(B_n(y,\varepsilon)\cap Y).
	\end{equation}

	Now we estimate the upper bound of $\mu_{l_r}(B_n(x,\varepsilon/2))$ for $n\geq s(2,0)$ and sufficiently large $l_r$.
	Write $y=x_{\theta}$.
	For any $x_{\xi}\in Y$, if $\xi_j^k\neq\theta_j^k$, then by the assumptions, $f^{s(j,k)}x_{\xi}$ and $f^{s(j,k)}x_{\theta}$ are $(t_j,\varepsilon)$-separated, and hence $x_{\xi}\notin B_n(y,\varepsilon)$ for all $n\geq s(j+1,0)-M_j$.

	{\bf Case 1.}
	If $s(j+1,0)-M_j\leq n\leq s(j+1,1)-1$ for some $j\geq2$, then $x_{\xi}\notin B_n(y,\varepsilon)$ whenever $\xi_i^k\neq\theta_i^k$ for some $1\leq i\leq j$ and $0\leq k\leq n_i-1$.
	Hence, using \eqref{eq-Bowen-lemma1-2} and \eqref{eq-Bowen-lemma1-proof-2}, for any $l_r\geq j+1$, we obtain that
	\begin{equation*}
		\mu_{l_r}(B_n(x,\varepsilon/2))\leq\frac{\prod_{i=j+1}^{l_r}\prod_{k=0}^{n_i-1}|\Omega_i(k)|}{\prod_{i=1}^{l_r}\prod_{k=0}^{n_i-1}|\Omega_i(k)|}=\frac{1}{\prod_{i=1}^j\prod_{k=0}^{n_i-1}|\Omega_i(k)|}\leq\exp\bigg(-\beta\sum_{i=1}^jn_it_i\bigg).
	\end{equation*}
	Therefore
	\begin{equation}\label{eq-Bowen-lemma1-proof-3}
		\begin{split}
			-\frac{1}{n}\log\mu(B_n(x,\varepsilon/2))&\geq-\frac{1}{n}\liminf_{r\rightarrow\infty}\log\mu_{l_r}(B_n(x,\varepsilon/2))\\
			&\geq-\frac{1}{n}\bigg(-\beta\sum_{i=1}^jn_it_i\bigg)\geq\frac{\sum_{i=1}^jn_it_i}{s(j+1,1)-1}\beta.
		\end{split}		
	\end{equation}

	{\bf Case 2.}
	Otherwise, there exist $j\in\mathbb{N}$ and $1\leq k\leq n_{j+1}-1$ such that
	\begin{equation*}
		s(j+1,k)\leq n\leq s(j+1,k)+t_{j+1}-1.
	\end{equation*}
	Then for $x_{\xi}$ with $\xi_i^u\neq\theta_i^u$ for some $1\leq i\leq j$ and $0\leq u\leq n_i-1$ or $i=j+1$ and $0\leq u\leq k-1$, we have $x_{\xi}\notin B_n(y,\varepsilon)$.
	Hence, using \eqref{eq-Bowen-lemma1-2} and \eqref{eq-Bowen-lemma1-proof-2}, for any $l_r\geq j+2$, we have
	\begin{equation*}
		\begin{split}
			\mu_{l_r}(B_n(x,\varepsilon/2))&\leq\frac{\Big(\prod_{s=k}^{n_{j+1}-1}|\Omega_{j+1}(s)|\Big)\cdot\Big(\prod_{i=j+2}^{l_r}\prod_{s=0}^{n_i-1}|\Omega_i(s)|\Big)}{\prod_{i=1}^{l_r}\prod_{s=0}^{n_i-1}|\Omega_i(s)|}\\
			&=\frac{1}{\Big(\prod_{s=0}^{k-1}|\Omega_{j+1}(s)|\Big)\cdot\Big(\prod_{i=1}^j\prod_{s=0}^{n_i-1}|\Omega_i(s)|\Big)}\leq\exp\Bigg(-\beta\bigg(k\cdot t_{j+1}+\sum_{i=1}^jn_it_i\bigg)\Bigg).
		\end{split}
	\end{equation*}
	Therefore
	\begin{equation}\label{eq-Bowen-lemma1-proof-4}
		\begin{split}
			-\frac{1}{n}\log\mu(B_n(x,\varepsilon/2))&\geq-\frac{1}{n}\liminf_{r\rightarrow\infty}\log\mu_{l_r}(B_n(x,\varepsilon/2))\\
			&\geq-\frac{1}{n}(-\beta\bigg(k\cdot t_{j+1}+\sum_{i=1}^jn_it_i\bigg)\Bigg)\\
			&\geq\frac{kt_{j+1}+\sum_{i=1}^jn_it_i}{s(j+1,k)+t_{j+1}-1}\beta=\frac{kt_{j+1}+\sum_{i=1}^jn_it_i}{kt_{j+1}+s(j+1,1)-1}\beta\geq\frac{\sum_{i=1}^jn_it_i}{s(j+1,1)-1}\beta.
		\end{split}
	\end{equation}
	Combining \eqref{eq-Bowen-lemma1-1}, \eqref{eq-Bowen-lemma1-proof-3} and \eqref{eq-Bowen-lemma1-proof-4}, we obtain that
	\begin{equation}\label{eq-Bowen-lemma1-proof-5}
		\liminf_{n\to\infty}-\frac{1}{n}\log\mu(B_n(x,\varepsilon/2))\geq\liminf_{j\rightarrow\infty}\frac{\sum_{i=1}^jn_it_i}{s(j+1,1)-1}\beta\geq\alpha\beta.
	\end{equation}
	Combining \eqref{eq-Bowen-lemma1-proof-1} and \eqref{eq-Bowen-lemma1-proof-5}, by the arbitrariness of $x\in\overline{Y}$, we conclude that
	\begin{equation*}
		h_{top}^B(f,\overline{Y})\geq\alpha\beta.
	\end{equation*}
\end{proof}

\subsection{The critical case}

The case $\sigma=\tau$ of statement (1) is crucial to the proof of Theorem \ref{thm-Bowen-lower-1}, and we state it as the following theorem.

%thm2
\begin{theorem}\label{Bowen-thm2}
	Suppose that $(X,f)$ satisfies the non-uniform specification property with gap function $M(n,\varepsilon)$.
	Assume that
	\begin{equation*}
		\tau:=\tau(M)=\sup_{\varepsilon>0}\liminf_{n\to\infty}\frac{M(n,\varepsilon)}{n}<1.
	\end{equation*}
	Then for any $\varphi\in C_{\tau}(X)$, we have
	\begin{equation*}
		h_{top}^B(f,I_{\varphi}(f))\geq\frac{1}{1+\tau}\sup_{(\mu,\nu)\in\mathfrak{D}_f(\varphi,\tau)}\min\{h_{\mu}(f),h_{\nu}(f)\}.
	\end{equation*}
\end{theorem}

\begin{proof}
%	Let $M(n,\varepsilon)\in\mathscr{G}(X,f)$ be a gap function satisfying
%	\begin{equation*}
%		\tau=\tau(M)=\sup_{\varepsilon>0}\liminf_{n\to\infty}\frac{M(n,\varepsilon)}{n}.
%	\end{equation*}
	Given $(\mu,\nu)\in\mathfrak{D}_f(\varphi,\tau)$.
	Our goal is to prove that
	\begin{equation*}
		h_{top}^B(f,I_{\varphi}(f))\geq\frac{1}{1+\tau}\min\{h_{\mu}(f),h_{\nu}(f)\}.
	\end{equation*}
	Without loss of generality, we assume that
	\begin{equation*}
		\int_X\varphi\mathrm{d}\mu<\int_X\varphi\mathrm{d}\nu.
	\end{equation*}
	Arbitrarily given $0<\gamma<\min\{h_{\mu}(f),h_{\nu}(f)\}$.
	Take $\kappa>0$ sufficiently small such that
	\begin{equation}\label{eq-Bowen-thm2-proof-1}
		\gamma+2\kappa<\min\{h_{\mu}(f),h_{\nu}(f)\}.
	\end{equation}

	Replacing $\varphi$ with $\varphi+C$ for some $C>0$, we may assume that $\inf_{x\in X}\varphi(x)>0$.
	Then
	\begin{equation*}
		\int_X\varphi\mathrm{d}\nu>\int_X\varphi\mathrm{d}\mu\geq\mu(X)\inf_{x\in X}\varphi(x)>0.
	\end{equation*}
	Fix $\eta>0$ sufficiently small such that
	\begin{equation}\label{eq-Bowen-thm2-proof-2}
		\int_X\varphi\mathrm{d}\nu-\int_X\varphi\mathrm{d}\mu>\tau\Big(\sup_{x\in X}\varphi(x)-\inf_{x\in X}\varphi(x)\Big)+10\eta.
	\end{equation}
	By Lemma \ref{pre-separation-lemma2}, we can take positive numbers $\delta^*>0$ and $\varepsilon^*>0$, positive integers $p_1,\cdots,p_m,q_1,\cdots,q_n$, and ergodic measures $\mu_1,\cdots,\mu_m,\nu_1,\cdots,\nu_n$ such that
	\begin{equation}\label{eq-Bowen-thm2-proof-3}
		\Bigg|\sum_{u=1}^m\frac{p_u}{p}\int_X\varphi\mathrm{d}\mu_u-\int_X\varphi\mathrm{d}\mu\Bigg|<\eta,\quad\Bigg|\sum_{u=1}^n\frac{q_u}{q}\int_X\varphi\mathrm{d}\nu_u-\int_X\varphi\mathrm{d}\nu\Bigg|<\eta,
	\end{equation}
	and
	\begin{equation}\label{eq-Bowen-thm2-proof-4}
		\min\left\{\sum_{u=1}^m\frac{p_u}{p}\underline{s}(\mu_u,\delta^*,\varepsilon^*),\sum_{u=1}^n\frac{q_u}{q}\underline{s}(\nu_u,\delta^*,\varepsilon^*)\right\}>\gamma+\kappa,
	\end{equation}
	where $p=p_1+\cdots+p_m$ and $q=q_1+\cdots+q_n$.
	Define
	\begin{itemize}
		\item $P_0=Q_0=0$;
		\item $P_u=p_1+\cdots+p_u$ for $1\leq u\leq m$;
		\item $Q_u=q_1+\cdots+q_u$ for $1\leq u\leq n$.
	\end{itemize}
	Then there exists $N^*\in\mathbb{N}$ such that for any $l\geq N^*$, we have
	\begin{equation}\label{eq-Bowen-thm2-proof-5}
		\begin{split}
			N\Big(\mathfrak{M}_u;\delta^*,l,\varepsilon^*\Big)&\geq\exp \big(l(\underline{s}(\mu_u,\delta^*,\varepsilon^*)-\kappa)\big)\quad\text{for}\quad u=1,\cdots,m,\\
			N\Big(\mathfrak{N}_u;\delta^*,l,\varepsilon^*\Big)&\geq\exp \big(l(\underline{s}(\nu_u,\delta^*,\varepsilon^*)-\kappa)\big)\quad\text{for}\quad u=1,\cdots,n,
		\end{split}
	\end{equation}
	where
	\begin{align*}
		\mathfrak{M}_u&=\left\{\lambda\in\mathcal{M}(X):\left|\int_X\varphi\mathrm{d}\lambda-\int_X\varphi\mathrm{d}\mu_u\right|<\eta\right\}\quad\text{for}\quad u=1,\cdots,m,\\
		\mathfrak{N}_u&=\left\{\lambda\in\mathcal{M}(X):\left|\int_X\varphi\mathrm{d}\lambda-\int_X\varphi\mathrm{d}\nu_u\right|<\eta\right\}\quad\text{for}\quad u=1,\cdots,n.
	\end{align*}
	By \eqref{eq-Bowen-thm2-proof-5}, for any $l\geq N^*$, we can take $(\delta^*,l,\varepsilon^*)$-separated sets $\Gamma_l^{\mu_u}\subset X_{l,\mathfrak{M}_u}$ and $\Gamma_l^{\nu_u}\subset X_{l,\mathfrak{N}_u}$ satisfying
	\begin{equation}\label{eq-Bowen-thm2-proof-6}
		|\Gamma_l^{\mu_u}|\geq\exp\big(l(\underline{s}(\mu_u,\delta^*,\varepsilon^*)-\kappa)\big)\quad\text{and}\quad|\Gamma_l^{\nu_u}|\geq\exp\big(l(\underline{s}(\nu_u,\delta^*,\varepsilon^*)-\kappa)\big),
	\end{equation}
	respectively.

	Consider a function $\Theta:[0,1/5]\rightarrow\mathbb{R}$ defined as
	\begin{align*}
		\Theta(\theta)=&\frac{1}{1+\tau+2\theta}\bigg(\int_X\varphi\mathrm{d}\nu-3\eta-(\tau+2\theta)\sup_{x\in X}\varphi(x)\bigg)\\
		&-\frac{1}{1+\tau-2\theta}\bigg(\int_X\varphi\mathrm{d}\mu+3\eta-(\tau-2\theta)\inf_{x\in X}\varphi(x)\bigg)-\theta\bigg(\int_X\varphi\mathrm{d}\nu-3\eta+\sup_{x\in X}\varphi(x)\bigg).
	\end{align*}
	Clearly $\Theta$ is continuous.
	It follows from \eqref{eq-Bowen-thm2-proof-2} that $\Theta(0)>4\eta/(1+\tau)>2\eta$.
	Hence we can take $\theta^*>0$ sufficiently small such that $\theta^*<\min\{1/N^*,\eta\}$ and $\Theta(\theta)>\eta$ for any $0<\theta<\theta^*$.
	In other words, for any $0<\theta<\theta^*$, we have
	\begin{equation}\label{eq-Bowen-thm2-proof-7}
		\begin{split}
			\Bigg(\bigg(\frac{1}{1+\tau+2\theta}-\theta\bigg)&\bigg(\int_X\varphi\mathrm{d}\nu-3\eta\bigg)+\frac{\tau-2\theta}{1+\tau-2\theta}\inf_{x\in X}\varphi(x)\Bigg)\\
			-&\Bigg(\frac{1}{1+\tau-2\theta}\bigg(\int_X\varphi\mathrm{d}\mu+3\eta\bigg)+\bigg(\frac{\tau+2\theta}{1+\tau+2\theta}+\theta\bigg)\sup_{x\in X}\varphi(x)\Bigg)>\eta.
		\end{split}
	\end{equation}
	Arbitrarily given $0<\theta<\theta^*$.
	Since $M(n,\varepsilon)$ is non-decreasing as $\varepsilon$ decreases, we have
	\begin{equation*}
		\tau=\sup_{\varepsilon>0}\liminf_{n\rightarrow\infty}\frac{M(n,\varepsilon)}{n}=\lim_{\varepsilon\rightarrow0}\liminf_{n\rightarrow\infty}\frac{M(n,\varepsilon)}{n}.
	\end{equation*}
	Therefore we can take $\varepsilon_0>0$ sufficiently small such that
	\begin{itemize}
		\item $\varepsilon_0<\varepsilon^*/3$;
		\item if two points $x,y\in X$ satisfy $d(x,y)\leq\varepsilon_0$, then $|\varphi(x)-\varphi(y)|\leq\eta$;
		\item for any $0<\varepsilon\leq\varepsilon_0$, one has
		\begin{equation}\label{eq-Bowen-thm2-proof-8}
			\left|\tau-\liminf_{n\rightarrow\infty}\frac{M(n,\varepsilon)}{n}\right|<\theta.
		\end{equation}
	\end{itemize}

	Now we define sequences of positive integers $\{n_i\}_{i=1}^{\infty}$ and $\{t_i\}_{i=1}^{\infty}$.
	By \eqref{eq-Bowen-thm2-proof-8}, we can take $l_0\in\mathbb{N}$ sufficiently large such that $l_0>\max\{N^*,1/\theta\}$ and
	\begin{equation}\label{eq-Bowen-thm2-proof-9}
		\left|\tau-\frac{M(l_0,\varepsilon_0)}{l_0}\right|<2\theta.
	\end{equation}
	Let $m_0=M(l_0,\varepsilon_0)$.
	For $i\in\mathbb{N}$, define
	\begin{equation*}
		t_i=
		\begin{cases}
			p(l_0+m_0) & \text{if}\;i\;\text{is odd},\\
			q(l_0+m_0) & \text{if}\;i\;\text{is even},\\
		\end{cases}
		\quad\text{and}\quad n_i=\big((p+q)(l_0+m_0)\big)^{4i}.
	\end{equation*}
	Then for any $i\in\mathbb{N}$, we have
	\begin{equation}\label{eq-Bowen-thm2-proof-10}
		\frac{t_{i+1}}{n_it_i}<\frac{1}{(p+q)(l_0+m_0)}<\theta
	\end{equation}
	and
	\begin{equation}\label{eq-Bowen-thm2-proof-11}
		\frac{(p+q)(l_0+m_0)\sum_{j=1}^in_jt_j}{n_{i+1}t_{i+1}}<\frac{\sum_{j=1}^{4i+2}\big((p+q)(l_0+m_0)\big)^j}{\big((p+q)(l_0+m_0)\big)^{4i+4}}<\frac{1}{(p+q)(l_0+m_0)}<\theta.
	\end{equation}
	
	Now we construct the index sets $\{\Omega_i\}_{i=1}^{\infty}$ satisfying \eqref{eq-Bowen-lemma1-2} in Lemma \ref{Bowen-lemma1}.
	For $i\in\mathbb{N}$, define
	\begin{equation*}
		\Omega_i=
		\begin{cases}
			\prod_{u=1}^m\big(\Gamma^{\mu_u}_{l_0}\big)^{p_u} & \text{if}\;i\;\text{is odd},\\
			\prod_{u=1}^n\big(\Gamma^{\nu_u}_{l_0}\big)^{q_u} & \text{if}\;i\;\text{is even},
		\end{cases}
	\end{equation*}
	where $\Gamma_{l_0}^{\mu_u}\subset X_{l_0,\mathfrak{M}_u}$ and $\Gamma^{\nu_u}_{l_0}\subset X_{l_0,\mathfrak{N}_u}$ are $(\delta^*,l_0,\varepsilon^*)$-separated sets satisfying \eqref{eq-Bowen-thm2-proof-6}.
	Then by \eqref{eq-Bowen-thm2-proof-4}, for any $k\in\mathbb{N}$, we have $|\Omega_k|<\infty$ and
	\begin{align*}
		|\Omega_{2k-1}|&=\prod_{u=1}^m\big|\Gamma^{\mu_u}_{l_0}\big|^{p_u}\geq\prod_{u=1}^m\exp\big(p_ul_0(\underline{s}(\mu_u,\delta^*,\varepsilon^*)-\kappa)\big)\geq\exp(pl_0\gamma)\geq\exp\frac{t_{2k-1}\gamma}{1+\tau+2\theta},\\
		|\Omega_{2k}|&=\prod_{u=1}^n\big|\Gamma^{\nu_u}_{l_0}\big|^{q_u}\geq\prod_{u=1}^n\exp\big(q_ul_0(\underline{s}(\nu_u,\delta^*,\varepsilon^*)-\kappa)\big)\geq\exp(ql_0\gamma)\geq\exp\frac{t_{2k}\gamma}{1+\tau+2\theta}.
	\end{align*}
	Hence for any $i\in\mathbb{N}$, we have
	\begin{equation}\label{eq-Bowen-thm2-proof-12}
		\exp\frac{\gamma\cdot t_i}{1+\tau+2\theta}\leq|\Omega_i|<\infty.
	\end{equation}
	Therefore $\{\Omega_i\}_{i=1}^{\infty}$ satisfies \eqref{eq-Bowen-lemma1-2} in Lemma \ref{Bowen-lemma1} by taking $\beta=\gamma/(1+\tau+2\theta)$.
	
	For $j\in\mathbb{N}$ and $0\leq k\leq n_j-1$, define
	\begin{equation*}
		s(j,k)=
		\begin{cases}
			kt_1 & \text{if}\quad j=1,\\
			\sum_{i=1}^{j-1}n_it_i+kt_j & \text{if}\quad j\geq2.
		\end{cases}
	\end{equation*}
	For $j\in\mathbb{N}$, $0\leq k\leq n_j-1$, and $0\leq r\leq p-1$ if $j$ is odd or $0\leq r\leq q-1$ if $j$ is even, let
	\begin{equation*}
		\begin{cases}
			a(j,k,r)=s(j,k)+r(l_0+m_0),\\
			b(j,k,r)=s(j,k)+r(l_0+m_0)+l_0-1,\\
			c_j=a(j+1,0,0)=\sum_{i=1}^jn_it_i.
		\end{cases}
	\end{equation*}
	Define
	\begin{equation*}
		\Xi_j=\prod_{i=1}^j\Omega_i^{n_i}
		\quad\text{and}\quad
		\Xi=\prod_{i=1}^{\infty}\Omega_i^{n_i}.
	\end{equation*}
	Fix $\xi=(\xi_1^0,\cdots,\xi_1^{n_1-1},\xi_2^0,\cdots,\xi_2^{n_2-1},\cdots)\in\Xi$.
	For $j\in\mathbb{N}$ and $0\leq k\leq n_j-1$, we write
	\begin{equation*}
		\xi_j^k=
		\begin{cases}
			(x_{j,k}^0,\cdots,x_{j,k}^{p-1})\in\Omega_j=\prod_{u=1}^m\big(\Gamma^{\mu_u}_{l_0}\big)^{p_u} & \text{if}\; j\;\text{is odd},\\
			(x_{j,k}^0,\cdots,x_{j,k}^{q-1})\in\Omega_j=\prod_{u=1}^n\big(\Gamma^{\nu_u}_{l_0}\big)^{q_u} & \text{if}\; j\;\text{is even},
		\end{cases}
	\end{equation*}
	where $x_{j,k}^r\in\Gamma^{\mu_u}_{l_0}$ for $P_{u-1}\leq r\leq P_u-1$ if $j$ is odd and $x_{j,k}^r\in\Gamma^{\nu_u}_{l_0}$ for $Q_{u-1}\leq r\leq Q_u-1$ if $j$ is even.

	Given $\xi\in\Xi$, we define a non-empty compact set $E_j^{\xi}$ inductively for each $j\in\mathbb{N}$, starting with
	\begin{equation*}
		E_1^{\xi}=\bigcap_{k=0}^{n_1-1}\bigcap_{r=0}^{p-1}\bigcap_{i=a(1,k,r)}^{b(1,k,r)}f^{-i}\Big(\overline{B(f^{i-a(1,k,r)}x_{1,k}^r,\varepsilon_0)}\Big)
	\end{equation*}
	and
	\begin{equation*}
		E_2^{\xi}=E_1^{\xi}\cap\Bigg(\bigcap_{k=0}^{n_2-1}\bigcap_{r=0}^{q-1}\bigcap_{i=a(2,k,r)}^{b(2,k,r)}f^{-i}\Big(\overline{B(f^{i-a(2,k,r)}x_{2,k}^r,\varepsilon_0)}\Big)\Bigg).
	\end{equation*}
	Since $a(j,k,r+1)-b(j,k,r)=m_0+1=M(l_0,\varepsilon_0)+1=M(b(j,k,r)-a(j,k,r)+1,\varepsilon_0)+1$, it follows from the non-uniform specification property that $E_1^{\xi}$ and $E_2^{\xi}$ are non-empty closed sets.
	
	Assume that $E_l^{\xi}$ has been defined for every $1\leq l\leq 2j$.
	Define $E_{2j+1}^{\xi}$ and $E_{2j+2}^{\xi}$ as
	\begin{equation*}
		E_{2j+1}^{\xi}=\bigg(\bigcap_{l=1}^{2j}E_l^{\xi}\bigg)\cap\Bigg(\bigcap_{k=0}^{n_{2j+1}-1}\bigcap_{r=0}^{p-1}\bigcap_{i=a(2j+1,k,r)}^{b(2j+1,k,r)}f^{-i}\Big(\overline{B(f^{i-a(2j+1,k,r)}x_{2j+1,k}^r,\varepsilon_0)}\Big)\Bigg)
	\end{equation*}
	and
	\begin{equation*}
		E_{2j+2}^{\xi}=\bigg(\bigcap_{l=1}^{2j+1}E_l^{\xi}\bigg)\cap\Bigg(\bigcap_{k=0}^{n_{2j+2}-1}\bigcap_{r=0}^{q-1}\bigcap_{i=a(2j+2,k,r)}^{b(2j+2,k,r)}f^{-i}\Big(\overline{B(f^{i-a(2j+2,k,r)}x_{2j+2,k}^r,\varepsilon_0)}\Big)\Bigg).
	\end{equation*}
	From the non-uniform specification property, it follows that $E_{2j+1}^{\xi}$ and $E_{2j+2}^{\xi}$ are non-empty closed sets.
	Similarly, for $l\in\mathbb{N}$, $1\leq j\leq l$ and $\zeta\in\Xi_l$, we can define a non-empty closed set $E_j^{\zeta}$.
	Note that if $\xi\in\Xi$ satisfies $\xi_j^k=\zeta_j^k$ for all $1\leq j\leq l$ and $0\leq k\leq n_j-1$, then $E_l^{\xi}=E_l^{\zeta}$.

	For $j\in\mathbb{N}$, define
	\begin{equation*}
		F_j=\bigcup_{\zeta\in\Xi_j}E_j^{\zeta}\quad\text{and}\quad F=\bigcap_{j\geq1}F_j.
	\end{equation*}
	Then $F_j$ and $F$ are non-empty and compact.
	We claim that $F\subset I_{\varphi}(f)$.
	For any $x\in F$, we will show that
	\begin{equation*}
		\limsup_{j\rightarrow\infty}\frac{1}{c_{2j-1}}\sum_{i=0}^{c_{2j-1}-1}\varphi(f^ix)<\liminf_{j\rightarrow\infty}\frac{1}{c_{2j}}\sum_{i=0}^{c_{2j}-1}\varphi(f^ix).
	\end{equation*}
	Since $x\in F$, there exists $\xi\in\Xi$ such that for any $j\in\mathbb{N}$, we have
	\begin{equation}\label{eq-Bowen-thm2-proof-13}
		\begin{cases}
			d(f^ix,f^{i-a(2j-1,k,r)}x_{2j-1,k}^r)\leq\varepsilon_0 & \text{for}\;0\leq k\leq n_{2j-1}-1,\;0\leq r\leq p-1\\&\quad\text{and}\;a(2j-1,k,r)\leq i\leq b(2j-1,k,r),\\
			d(f^ix,f^{i-a(2j,k,r)}x_{2j,k}^r)\leq\varepsilon_0 & \text{for}\;0\leq k\leq n_{2j}-1\,;0\leq r\leq q-1\\&\quad\text{and}\;a(2j,k,r)\leq i\leq b(2j,k,r),
		\end{cases}
	\end{equation}
	where $x_{2j-1,k}^r\in\Gamma_{l_0}^{\mu_u}$ for $P_{u-1}\leq r\leq P_u-1$ and $x_{2j,k}^r\in\Gamma_{l_0}^{\nu_u}$ for $Q_{u-1}\leq r\leq Q_u-1$ are points satisfying $\xi_{2j-1}^k=(x_{2j-1,k}^0,\cdots,x_{2j-1,k}^{p-1})$ and $\xi_{2j}^k=(x_{2j,k}^0,\cdots,x_{2j,k}^{q-1})$, respectively.
	Note that
	\begin{equation*}
		\Omega_{2j-1}=\prod_{u=1}^m\big(\Gamma^{\mu_u}_{l_0}\big)^{p_u}\subset\prod_{u=1}^m\big(X_{l_0,\mathfrak{M}_u}\big)^{p_u}.
	\end{equation*}
	It follows that for any $0\leq k\leq n_{2j-1}-1$, any $1\leq u\leq m$ and any $P_{u-1}\leq r\leq P_u-1$, we have $\delta_{x_{2j-1,k}^r}^{l_0}\in\mathfrak{M}_u$.
	Hence
	\begin{equation*}
		\Bigg|\int_X\varphi\mathrm{d}\mu_u-\frac{1}{l_0}\sum_{i=0}^{l_0-1}\varphi(f^ix_{2j-1,k}^r)\Bigg|=\Bigg|\int_X\varphi\mathrm{d}\mu_u-\int_X\varphi\mathrm{d}\delta_{x_{2j-1,k}^r}^{l_0}\Bigg|<\eta.
	\end{equation*}
	This implies that
	\begin{equation}\label{eq-Bowen-thm2-proof-14}
		\frac{1}{l_0}\sum_{i=0}^{l_0-1}\varphi(f^ix_{2j-1,k}^r)<\int_X\varphi\mathrm{d}\mu_u+\eta.
	\end{equation}
	Combining \eqref{eq-Bowen-thm2-proof-3} and \eqref{eq-Bowen-thm2-proof-14}, we obtain that
	\begin{equation}\label{eq-Bowen-thm2-proof-15}
		\frac{1}{pl_0}\sum_{r=0}^{p-1}\sum_{i=0}^{l_0-1}\varphi(f^ix_{2j-1,k}^r)<\sum_{u=1}^m\frac{p_u}{p}\bigg(\int_X\varphi\mathrm{d}\mu_u+\eta\bigg)<\int_X\varphi\mathrm{d}\mu+2\eta.
	\end{equation}
	Similarly we have
	\begin{equation}\label{eq-Bowen-thm2-proof-16}
		\frac{1}{ql_0}\sum_{r=0}^{q-1}\sum_{i=0}^{l_0-1}\varphi(f^ix_{2j,k}^r)>\int_X\varphi\mathrm{d}\nu-2\eta.
	\end{equation}
	Since $|\varphi(y)-\varphi(z)|\leq\eta$ whenever $d(y,z)\leq\varepsilon_0$, it follows from \eqref{eq-Bowen-thm2-proof-13} that
	\begin{equation*}
		\Bigg|\frac{1}{pn_{2j-1}l_0}\sum_{k=0}^{n_{2j-1}-1}\sum_{r=0}^{p-1}\sum_{i=0}^{l_0-1}\varphi(f^ix_{2j-1,k}^r)-\frac{1}{pn_{2j-1}l_0}\sum_{k=0}^{n_{2j-1}-1}\sum_{r=0}^{p-1}\sum_{i=a(2j-1,k,r)}^{b(2j-1,k,r)}\varphi(f^ix)\Bigg|\leq\eta
	\end{equation*}
	and
	\begin{equation*}
		\left|\frac{1}{qn_{2j}l_0}\sum_{k=0}^{n_{2j}-1}\sum_{r=0}^{q-1}\sum_{i=0}^{l_0-1}\varphi(f^ix_{2j,k}^r)-\frac{1}{qn_{2j}l_0}\sum_{k=0}^{n_{2j}-1}\sum_{r=0}^{q-1}\sum_{i=a(2j,k,r)}^{b(2j,k,r)}\varphi(f^ix)\right|\leq\eta.
	\end{equation*}
	Hence
	\begin{equation}\label{eq-Bowen-thm2-proof-17}
		\frac{1}{pn_{2j-1}l_0}\sum_{k=0}^{n_{2j-1}-1}\sum_{r=0}^{p-1}\sum_{i=a(2j-1,k,r)}^{b(2j-1,k,r)}\varphi(f^ix)\leq\frac{1}{pn_{2j-1}l_0}\sum_{k=0}^{n_{2j-1}-1}\sum_{r=0}^{p-1}\sum_{i=0}^{l_0-1}\varphi(f^ix_{2j-1,k}^r)+\eta
	\end{equation}
	and
	\begin{equation}\label{eq-Bowen-thm2-proof-18}
		\frac{1}{qn_{2j}l_0}\sum_{k=0}^{n_{2j}-1}\sum_{r=0}^{q-1}\sum_{i=a(2j,k,r)}^{b(2j,k,r)}\varphi(f^ix)\geq\frac{1}{qn_{2j}l_0}\sum_{k=0}^{n_{2j}-1}\sum_{r=0}^{q-1}\sum_{i=0}^{l_0-1}\varphi(f^ix_{2j,k}^r)-\eta.
	\end{equation}
	Combining \eqref{eq-Bowen-thm2-proof-15}, \eqref{eq-Bowen-thm2-proof-16}, \eqref{eq-Bowen-thm2-proof-17} and \eqref{eq-Bowen-thm2-proof-18}, we obtain that
	\begin{equation}\label{eq-Bowen-thm2-proof-19}
		\frac{1}{pn_{2j-1}l_0}\sum_{k=0}^{n_{2j-1}-1}\sum_{r=0}^{p-1}\sum_{i=a(2j-1,k,r)}^{b(2j-1,k,r)}\varphi(f^ix)<\int_X\varphi\mathrm{d}\mu+3\eta
	\end{equation}
	and
	\begin{equation}\label{eq-Bowen-thm2-proof-20}
		\frac{1}{qn_{2j}l_0}\sum_{k=0}^{n_{2j}-1}\sum_{r=0}^{q-1}\sum_{i=a(2j,k,r)}^{b(2j,k,r)}\varphi(f^ix)>\int_X\varphi\mathrm{d}\nu-3\eta.
	\end{equation}
	Let
	\begin{align*}
		\Lambda_{2j-1}&=\mathbb{Z}\cap\Bigg([0,c_{2j-1}-1]\setminus\bigg(\bigcup_{k=0}^{n_{2j-1}-1}\bigcup_{r=0}^{p-1}[a(2j-1,k,r),b(2j-1,k,r)]\bigg)\Bigg),\\
		\Lambda_{2j}&=\mathbb{Z}\cap\Bigg([0,c_{2j}-1]\setminus\bigg(\bigcup_{k=0}^{n_{2j}-1}\bigcup_{r=0}^{q-1}[a(2j,k,r),b(2j,k,r)]\bigg)\Bigg),
	\end{align*}
	then $|\Lambda_{2j-1}|=c_{2j-1}-pn_{2j-1}l_0$ and $|\Lambda_{2j}|=c_{2j}-qn_{2j}l_0$.
	Using \eqref{eq-Bowen-thm2-proof-19}, for any $j\in\mathbb{N}$, we obtain that
	\begin{align*}
		\frac{1}{c_{2j-1}}\sum_{i=0}^{c_{2j-1}-1}\varphi(f^ix)
		&=\frac{1}{c_{2j-1}}\Bigg(\sum_{k=0}^{n_{2j-1}-1}\sum_{r=0}^{p-1}\sum_{i=a(2j-1,k,r)}^{b(2j-1,k,r)}\varphi(f^ix)+\sum_{i\in\Lambda_{2j-1}}\varphi(f^ix)\Bigg)\\
		&=\frac{pn_{2j-1}l_0}{c_{2j-1}}\frac{1}{pn_{2j-1}l_0}\sum_{k=0}^{n_{2j-1}-1}\sum_{r=0}^{p-1}\sum_{i=a(2j-1,k,r)}^{b(2j-1,k,r)}\varphi(f^ix)+\frac{1}{c_{2j-1}}\sum_{i\in\Lambda_{2j-1}}\varphi(f^ix)\\
		&<\frac{pn_{2j-1}l_0}{c_{2j-1}}\bigg(\int_X\varphi\mathrm{d}\mu+3\eta\bigg)+\frac{c_{2j-1}-pn_{2j-1}l_0}{c_{2j-1}}\sup_{x\in X}\varphi(x).
	\end{align*}
	By \eqref{eq-Bowen-thm2-proof-11}, we have
	\begin{equation*}
		0<\frac{l_0}{l_0+m_0}-\frac{pn_{2j-1}l_0}{c_{2j-1}}=\frac{pl_0}{t_{2j-1}}-\frac{pn_{2j-1}l_0}{c_{2j-1}}=\frac{pl_0\sum_{i=1}^{2j-2}n_it_i}{c_{2j-1}t_{2j-1}}<\theta,
	\end{equation*}
	and hence
	\begin{equation*}
		0<\frac{l_0}{l_0+m_0}-\frac{pn_{2j-1}l_0}{c_{2j-1}}<\theta.
	\end{equation*}
	Therefore
	\begin{equation}\label{eq-Bowen-thm2-proof-21}
		\begin{split}
			\frac{1}{c_{2j-1}}\sum_{i=0}^{c_{2j-1}-1}\varphi(f^ix)
			&<\frac{l_0}{l_0+m_0}\bigg(\int_X\varphi\mathrm{d}\mu+3\eta\bigg)+\bigg(\frac{m_0}{l_0+m_0}+\theta\bigg)\sup_{x\in X}\varphi(x)\\
			&\leq\frac{1}{1+\tau-2\theta}\bigg(\int_X\varphi\mathrm{d}\mu+3\eta\bigg)+\bigg(\frac{\tau+2\theta}{1+\tau+2\theta}+\theta\bigg)\sup_{x\in X}\varphi(x).
		\end{split}
	\end{equation}
	Similarly we have
	\begin{equation*}
		\frac{1}{c_{2j}}\sum_{i=0}^{c_{2j}-1}\varphi(f^ix)>\frac{qn_{2j}l_0}{c_{2j}}\left(\int_X\varphi\mathrm{d}\nu-3\eta\right)+\frac{c_{2j}-qn_{2j}l_0}{c_{2j}}\inf_{x\in X}\varphi(x)
	\end{equation*}
	and
	\begin{equation*}
		0<\frac{l_0}{l_0+m_0}-\frac{qn_{2j}l_0}{c_{2j}}<\theta.
	\end{equation*}
	Therefore
	\begin{equation}\label{eq-Bowen-thm2-proof-22}
		\frac{1}{c_{2j}}\sum_{i=0}^{c_{2j}-1}\varphi(f^ix)>\bigg(\frac{1}{1+\tau+2\theta}-\theta\bigg)\bigg(\int_X\varphi\mathrm{d}\nu-3\eta\bigg)+\frac{\tau-2\theta}{1+\tau-2\theta}\inf_{x\in X}\varphi(x).
	\end{equation}
	Combining \eqref{eq-Bowen-thm2-proof-7}, \eqref{eq-Bowen-thm2-proof-21} and \eqref{eq-Bowen-thm2-proof-22}, we conclude that
	\begin{equation*}
		\liminf_{j\rightarrow\infty}\frac{1}{c_{2j}}\sum_{i=0}^{c_{2j}-1}\varphi(f^ix)-\limsup_{j\rightarrow\infty}\frac{1}{c_{2j-1}}\sum_{i=0}^{c_{2j-1}-1}\varphi(f^ix)>\eta,
	\end{equation*}
	which implies that $x\in I_{\varphi}(f)$.
	Therefore $F\subset I_{\varphi}(f)$.

	For $\xi\in\Xi$, define
	\begin{equation*}
		E^{\xi}=\bigcap_{j\geq1}E_j^{\xi}.
	\end{equation*}
	Since $E^{\xi}_k\supset E^{\xi}_{k+1}$, the set $E^{\xi}$ is non-empty and closed.
	Clearly $E_j^{\xi}\subset F_j$ and $E^{\xi}\subset F$ for all $\xi\in\Xi$ and $j\in\mathbb{N}$.
	Arbitrarily take $x_{\xi}\in E^{\xi}$.
	Let
	\begin{equation*}
		Y=\bigg\{x_{\xi}\in X:\xi=(\xi_1^0,\cdots,\xi_1^{n_1-1},\xi_2^0,\cdots,\xi_2^{n_2-1},\cdots)\in\prod_{i=1}^{\infty}\Omega_i^{n_i}\bigg\}.
	\end{equation*}
	Then $Y\subset F$.
	To apply Lemma \ref{Bowen-lemma1}, we need to verify that for any $\xi,\zeta\in\Xi$, if $\xi_j^k\neq\zeta_j^k$ for some $j\in\mathbb{N}$ and $0\leq k\leq n_j-1$, then $f^{s(j,k)}x_{\xi}$ and $f^{s(j,k)}x_{\zeta}$ are $(t_j,\varepsilon_0)$-separated.
	For simplicity, we only consider the case where $j$ is odd, as the case for even $j$ follows analogously.
	We write $\xi_j^k=(x_{j,k}^0,\cdots,x_{j,k}^{p-1})$ and $\zeta_j^k=(z_{j,k}^0,\cdots,z_{j,k}^{p-1})$ with $x_{j,k}^r,z_{j,k}^r\in\Gamma_{l_0}^{\mu_u}$ for $u=1,\cdots,m$ and $P_{u-1}\leq r\leq P_u-1$.
	Since $\xi_j^k\neq\zeta_j^k$, there is an $r$ such that $x_{j,k}^r\neq z_{j,k}^r$.
	From the definition of $E_j^{\xi}$, it follows that for any $0\leq i\leq l_0-1$, we have
	\begin{equation}\label{eq-Bowen-thm2-proof-23}
		\begin{split}
			d(f^{a(j,k,r)+i}x_{\xi},f^{a(j,k,r)+i}x_{\zeta})
			&\geq d(f^ix_{j,k}^r,f^iz_{j,k}^r)-d(f^{a(j,k,r)+i}x_{\xi},f^ix_{j,k}^r)-d(f^{a(j,k,r)+i}x_{\zeta},f^iz_{j,k}^r)\\
			&\geq d(f^ix_{j,k}^r,f^iz_{j,k}^r)-2\varepsilon_0.
		\end{split}
	\end{equation}
	Take $1\leq u\leq m$ such that $P_{u-1}\leq r\leq P_u-1$.
	Since $\Gamma_{l_0}^{\mu_u}$ is a $(\delta^*,l_0,\varepsilon^*)$-separated set and $x_{j,k}^r\neq z_{j,k}^r$, there exists $i\in\{0,\cdots,l_0-1\}$ such that
	\begin{equation}\label{eq-Bowen-thm2-proof-24}
		d(f^ix_{j,k}^r,f^iz_{j,k}^r)>\varepsilon^*\geq3\varepsilon_0.
	\end{equation}
	Let $\iota=r(l_0+m_0)+i$, then $0\leq\iota\leq t_j-1$.
	Combining \eqref{eq-Bowen-thm2-proof-23} and \eqref{eq-Bowen-thm2-proof-24}, we conclude that
	\begin{equation*}
		d(f^{s(j,k)+\iota}x_{\xi},f^{s(j,k)+\iota}x_{\zeta})=d(f^{a(j,k,r)+i}x_{\xi},f^{a(j,k,r)+i}x_{\zeta})>\varepsilon_0.
	\end{equation*}
	As a result, $f^{s(j,k)}x_{\xi}$ and $f^{s(j,k)}x_{\zeta}$ are $(t_j,\varepsilon_0)$-separated.

	From \eqref{eq-Bowen-thm2-proof-10}, it follows that
	\begin{equation}\label{eq-Bowen-thm2-proof-25}
		\liminf_{j\rightarrow\infty}\frac{\sum_{i=1}^jn_it_i}{s(j+1,1)-1}=\liminf_{j\rightarrow\infty}\frac{\sum_{i=1}^jn_it_i}{\sum_{i=1}^jn_it_i+t_{j+1}-1}\geq\frac{1}{1+\theta}.
	\end{equation}
	Combining \eqref{eq-Bowen-thm2-proof-12} and \eqref{eq-Bowen-thm2-proof-25}, by Lemma \ref{Bowen-lemma1}, we obtain that
	\begin{equation*}
		h_{top}^B(f,\overline{Y})\geq\frac{\gamma}{(1+\theta)(1+\tau+2\theta)}.
	\end{equation*}
	Since $Y\subset F\subset I_{\varphi}(f)$ and $F$ is compact, it follows that
	\begin{equation*}
		h_{top}^B(f,I_{\varphi}(f))\geq h_{top}^B(f,F)\geq h_{top}^B(f,\overline{Y})\geq\frac{\gamma}{(1+\theta)(1+\tau+2\theta)}.
	\end{equation*}
	By the arbitrariness of $\gamma\in(0,\min\{h_{\mu}(f),h_{\nu}(f)\})$ and $\theta\in(0,\theta^*)$, we conclude that
	\begin{equation*}
		h_{top}^B(f,I_{\varphi}(f))\geq\frac{1}{1+\tau}\min\{h_{\mu}(f),h_{\nu}(f)\}.
	\end{equation*}
	Therefore
	\begin{equation*}
		h_{top}^B(f,I_{\varphi}(f))\geq\frac{1}{1+\tau}\sup_{(\mu,\nu)\in\mathfrak{D}_f(\varphi,\tau)}\min\{h_{\mu}(f),h_{\nu}(f)\}.
	\end{equation*}
	This completes the proof of Theorem \ref{Bowen-thm2}.
\end{proof}

\subsection{Lower bound for Bowen topological entropy of the irregular set}

Now we proceed to show Theorem \ref{thm-Bowen-lower-1}.
Fix a gap function $M(n,\varepsilon)\in\mathscr{G}(X,f)$ satisfying
\begin{equation*}
	\tau=\tau(M)=\sup_{\varepsilon>0}\liminf_{n\to\infty}\frac{M(n,\varepsilon)}{n}.
\end{equation*}

\subsubsection{Proof of Theorem \ref{thm-Bowen-lower-1} (1)}

\begin{proof}[{\bf Proof of Theorem \ref{thm-Bowen-lower-1} (1)}]
	We have shown in Theorem \ref{Bowen-thm2} that the statement (1) of Theorem \ref{thm-Bowen-lower-1} holds for the case $\sigma=\tau$.
	Assume that $0<\tau<\sigma<1$.
	Given $(\mu,\nu)\in\mathfrak{D}_f(\varphi,\sigma)$.
	We aim to show that
	\begin{equation*}
		h_{top}^B(f,I_{\varphi}(f))\geq\frac{\sigma-\tau}{(1+\tau)\sigma}h_{top}(f)+\frac{\tau}{(1+\tau)\sigma}\min\{h_{\mu}(f),h_{\nu}(f)\}.
	\end{equation*}
	Without loss of generality, we assume that
	\begin{equation*}
		\int_X\varphi\mathrm{d}\mu<\int_X\varphi\mathrm{d}\nu.
	\end{equation*}
	Arbitrarily given $0<\gamma<h_{top}(f)$.
	From the traditional variational principle, we can take $\lambda\in\mathcal{M}_f(X)$ such that $h_{\lambda}(f)>\gamma$.
	Let $\tilde{\mu}=(\tau/\sigma)\mu+(1-\tau/\sigma)\lambda$ and $\tilde{\nu}=(\tau/\sigma)\nu+(1-\tau/\sigma)\lambda$.
	Then
	\begin{align*}
		\min\{h_{\tilde{\mu}}(f),h_{\tilde{\nu}}(f)\}&=\frac{\sigma-\tau}{\sigma}h_{\lambda}(f)+\frac{\tau}{\sigma}\min\{h_{\mu}(f),h_{\nu}(f)\}\\
		&>\frac{(\sigma-\tau)}{\sigma}\gamma+\frac{\tau}{\sigma}\min\{h_{\mu}(f),h_{\nu}(f)\}.
	\end{align*}
	Since $(\mu,\nu)\in\mathfrak{D}_f(\varphi,\sigma)$, we have
	\begin{equation*}
		\int_X\varphi\mathrm{d}\tilde{\nu}-\int_X\varphi\mathrm{d}\tilde{\mu}=\frac{\tau}{\sigma}\bigg(\int_X\varphi\mathrm{d}\nu-\int_X\varphi\mathrm{d}\mu\bigg)>\tau\Big(\sup_{x\in X}\varphi(x)-\inf_{x\in X}\varphi(x)\Big).
	\end{equation*}
	Hence $(\tilde{\mu},\tilde{\nu})\in\mathfrak{D}_f(\varphi,\tau)$.
	It follows from Theorem \ref{Bowen-thm2} that
	\begin{align*}
		h_{top}^B(f,I_{\varphi}(f))&\geq\frac{1}{1+\tau}\min\{h_{\tilde{\mu}}(f),h_{\tilde{\nu}}(f)\}\\
		&>\frac{(\sigma-\tau)}{(1+\tau)\sigma}\gamma+\frac{\tau}{(1+\tau)\sigma}\min\{h_{\mu}(f),h_{\nu}(f)\}.
	\end{align*}
	By the arbitrariness of $\gamma\in(0,h_{top}(f))$, we conclude that
	\begin{equation*}
		h_{top}^B(f,I_{\varphi}(f))\geq\frac{\sigma-\tau}{(1+\tau)\sigma}h_{top}(f)+\frac{\tau}{(1+\tau)\sigma}\min\{h_{\mu}(f),h_{\nu}(f)\}.
	\end{equation*}
	Therefore
	\begin{equation}\label{eq-Bowen-thmC-proof-1}
		h_{top}^B(f,I_{\varphi}(f))\geq\frac{\sigma-\tau}{(1+\tau)\sigma}h_{top}(f)+\frac{\tau}{(1+\tau)\sigma}\sup_{(\mu,\nu)\in\mathfrak{D}_f(\varphi,\sigma)}\min\{h_{\mu}(f),h_{\nu}(f)\}.
	\end{equation}
	This completes the proof of statement (1).

	Now we assume that $\tau=0$.
	For any $\varphi\in\mathfrak{C}(f)$, there is a $0<\sigma<1$ such that $\varphi\in C_{\sigma}(X)$.
	For any $0<\kappa<\sigma$, let $M_{\kappa}(n,\varepsilon):=M(n,\varepsilon)+\lceil\kappa n\rceil$.
	Then $M_{\kappa}$ is an admissible gap function with $\tau(M_{\kappa})=\kappa$.
	By \eqref{eq-Bowen-thmC-proof-1}, we have
	\begin{equation*}
		h_{top}^B(f,I_{\varphi}(f))\geq\frac{\sigma-\kappa}{(1+\kappa)\sigma}h_{top}(f).
	\end{equation*}
	By the arbitrariness of $\kappa$, we conclude that
	\begin{equation*}
		h_{top}^B(f,I_{\varphi}(f))=h_{top}(f).
	\end{equation*}
	This completes the proof.
\end{proof}

\subsubsection{Proof of Theorem \ref{thm-Bowen-lower-1} (2)}

Before proving Theorem \ref{thm-Bowen-lower-1} (2), we provided some preliminary results.

%lemma3
\begin{lemma}[cf. the proof of {\cite[Lemma 3.4 (2)]{Hou-Lin-Tian-2023}}]\label{Bowen-lemma3}
	Let $(X,f)$ be a dynamical system and $0<\sigma<1$.
	Suppose that $(X,f)$ is not uniquely ergodic.
	Then for any two ergodic measures $\mu,\nu$ with $\mu\neq\nu$, there is a $\varphi\in C_{\sigma}(X)$ such that $(\mu,\nu)\in\mathfrak{D}_f(\varphi,\sigma)$.
\end{lemma}

By \cite[Theorem 3.2 (8)]{Lin-Tian-Yu-2024}, the non-uniform specification property implies $\mathfrak{C}(f)\neq\varnothing$.
Hence $(X,f)$ is not uniquely ergodic.
As a consequence of Lemma \ref{Bowen-lemma3}, we obtain the following result.

%cor4
\begin{corollary}\label{Bowen-cor4}
	If $(X,f)$ satisfies the non-uniform specification property, then for any $0<\sigma<1$ and any two distinct ergodic measures $\mu,\nu\in\mathcal{M}_f^{erg}(X)$, there exists $\varphi\in C_{\sigma}(X)$ such that $(\mu,\nu)\in\mathfrak{D}_f(\varphi,\sigma)$.
\end{corollary}

\begin{proof}[{\bf Proof of Theorem \ref{thm-Bowen-lower-1} (2)}]
	We proceed to show statement (2) of Theorem \ref{thm-Bowen-lower-1}.
	Given $\mu,\nu\in\mathcal{M}_f^{erg}(X)$ with $\mu\neq\nu$, we aim to show that
	\begin{equation*}
		h_{top}^B(f,\mathrm{IR}(f))\geq\frac{1-\tau}{1+\tau}h_{top}(f)+\frac{\tau}{1+\tau}\min\{h_{\mu}(f),h_{\nu}(f)\}.
	\end{equation*}
	Arbitrarily given $\tau<\sigma<1$.
	By Corollary \ref{Bowen-cor4}, there is a $\varphi\in C_{\sigma}(X)$ such that $(\mu,\nu)\in\mathfrak{D}_f(\varphi,\sigma)$.
	From Theorem \ref{thm-Bowen-lower-1}, it follows that
	\begin{equation*}
		h_{top}^B(f,\mathrm{IR}(f))\geq h_{top}^B(f,I_{\varphi}(f))\geq\frac{\sigma-\tau}{(1+\tau)\sigma}h_{top}(f)+\frac{\tau}{(1+\tau)\sigma}\min\{h_{\mu}(f),h_{\nu}(f)\}.
	\end{equation*}
	By taking the limit as $\sigma\rightarrow1$, we obtain that
	\begin{equation*}
		h_{top}^B(f,\mathrm{IR}(f))\geq\frac{1-\tau}{1+\tau}h_{top}(f)+\frac{\tau}{1+\tau}\min\{h_{\mu}(f),h_{\nu}(f)\}.
	\end{equation*}
	Therefore
	\begin{equation*}
		h_{top}^B(f,\mathrm{IR}(f))\geq\frac{1-\tau}{1+\tau}h_{top}(f)+\frac{\tau}{1+\tau}\sup_{\mu,\nu\in\mathcal{M}_f^{erg}(X),\mu\neq\nu}\min\{h_{\mu}(f),h_{\nu}(f)\}.
	\end{equation*}
	This completes the proof of statement (2).
\end{proof}

\subsubsection{Proof of Theorem \ref{thm-Bowen-lower-1} (3) and Corollary \ref{cor-Bowen-lower-1}}

To prove Theorem \ref{thm-Bowen-lower-1} (3), we need to estimate
\begin{equation*}
	\sup_{\mu,\nu\in\mathcal{M}_f^{erg}(X),\mu\neq\nu}\min\{h_{\mu}(f),h_{\nu}(f)\}.
\end{equation*}

%prop5
\begin{proposition}\label{Bowen-prop5}
	Suppose that $(X,f)$ satisfies the non-uniform specification property with gap function $M(n,\varepsilon)$.
	If $\tau:=\tau(M)<\infty$, then for any $\gamma<h_{top}(f)/(1+\tau)$, there exist at least two ergodic measures whose measure-theoretic entropies are greater than $\gamma$.
\end{proposition}

\begin{proof}
	From the traditional variational principle of topological entropy, the ergodic decomposition and the affine character of the entropy, it follows that, if there is no measure of maximal entropy, then we can take $\mu,\nu\in\mathcal{M}_f^{erg}(X)$ such that $\gamma<h_{\nu}(f)<h_{\mu}(f)<h_{top}(f)$, which is what we desire.

	Now we fix an ergodic measure $\mu\in\mathcal{M}_f^{erg}(X)$ with $h_{\mu}(f)=h_{top}(f)$.
	Take $q\in\mathbb{N}$ sufficiently large and $\theta^*>0$ sufficiently small such that
	\begin{equation}\label{eq-Bowen-prop5-proof-1}
		\gamma<\frac{q-1}{q(1+\tau+2\theta^*)}h_{top}(f).
	\end{equation}
	Take $0<\sigma<1$ such that
	\begin{equation}\label{eq-Bowen-prop5-proof-2}
		\sigma>q\tau(1-\sigma).
	\end{equation}

	Given $\nu_0\in\mathcal{M}_f^{erg}(X)\setminus\{\mu\}$.
	From Corollary \ref{Bowen-cor4}, we can choose $\varphi\in C_{\sigma}(X)$ such that $\inf_{x\in X}\varphi(x)=0$ and
	\begin{equation*}
		\int_X\varphi\mathrm{d}\mu-\int_X\varphi\mathrm{d}\nu_0>\sigma\Big(\sup_{x\in X}\varphi(x)-\inf_{x\in X}\varphi(x)\Big)=\sigma\sup_{x\in X}\varphi(x)=\sigma\|\varphi\|.
	\end{equation*}
	Then
	\begin{equation}\label{eq-Bowen-prop5-proof-3}
		\sigma\|\varphi\|\leq\int_X\varphi\mathrm{d}\mu\leq\|\varphi\|\quad\text{and}\quad 0\leq\int_X\varphi\mathrm{d}\nu_0\leq(1-\sigma)\|\varphi\|.
	\end{equation}
	Let $\mu_0=(1-1/q)\mu+(1/q)\nu_0$, then
	\begin{equation}\label{eq-Bowen-prop5-proof-4}
		\int_X\varphi\mathrm{d}\mu-\int_X\varphi\mathrm{d}\mu_0=\frac{1}{q}\bigg(\int_X\varphi\mathrm{d}\mu-\int_X\varphi\mathrm{d}\nu_0\bigg)>\frac{\sigma\|\varphi\|}{q}.
	\end{equation}

	Consider a function $\Theta:[0,1/5]\rightarrow\mathbb{R}$ defined by
	\begin{equation*}
		\Theta(\theta)=\int_X\varphi\mathrm{d}\mu-\bigg(\frac{1}{1+\tau-2\theta}\int_X\varphi\mathrm{d}\mu_0+\frac{\tau+2\theta}{1+\tau+2\theta}\|\varphi\|\bigg).
	\end{equation*}
	Clearly $\Theta$ is continuous.
	Combining \eqref{eq-Bowen-prop5-proof-2}, \eqref{eq-Bowen-prop5-proof-3} and \eqref{eq-Bowen-prop5-proof-4}, we obtain that
	\begin{align*}
		\Theta(0)&=\int_X\varphi\mathrm{d}\mu-\frac{1}{1+\tau}\bigg(\int_X\varphi\mathrm{d}\mu_0+\tau\|\varphi\|\bigg)\\
		&=\frac{1}{1+\tau}\Bigg(\tau\int_X\varphi\mathrm{d}\mu+\bigg(\int_X\varphi\mathrm{d}\mu-\int_X\varphi\mathrm{d}\mu_0\bigg)-\tau\|\varphi\|\Bigg)\\
		&>\frac{1}{1+\tau}\bigg(\tau\int_X\varphi\mathrm{d}\mu-\Big(\tau-\frac{\sigma}{q}\Big)\|\varphi\|\bigg)>\frac{\big(\sigma-q\tau(1-\sigma)\big)\|\varphi\|}{q(1+\tau)}>0.
	\end{align*}
	Take $\eta>0$ sufficiently small such that $\Theta(0)>6\eta$.
	Then there is a $0<\theta_0<\theta^*/2$ such that for any $0<\theta<\theta_0$, we have $\Theta(\theta)>5\eta$.
	In other words, for any $0<\theta<\theta_0$, we have
	\begin{equation}\label{eq-Bowen-prop5-proof-5}
		\frac{1}{1+\tau-2\theta}\int_X\varphi\mathrm{d}\mu_0+\frac{\tau+2\theta}{1+\tau+2\theta}\|\varphi\|<\int_X\varphi\mathrm{d}\mu-5\eta.
	\end{equation}
	Let
	\begin{equation*}
		\mathfrak{M}=\Bigg\{\lambda\in\mathcal{M}(X):\bigg|\int_X\varphi\mathrm{d}\lambda-\int_X\varphi\mathrm{d}\mu\bigg|<\eta\Bigg\}.
	\end{equation*}
	Then $\mathfrak{M}$ is an open neighborhood in $\mathcal{M}(X)$ of $\mu$.
	By \eqref{eq-Bowen-prop5-proof-1}, it follows from Proposition \ref{pre-separation-prop1} that there exist $\delta^*>0$, $\varepsilon^*>0$ and $N^*\in\mathbb{N}$ such that for any $n\geq N^*$, there is a $(\delta^*,n,\varepsilon^*)$-separated set $\Gamma_n\subset X_{n,\mathfrak{M}}$ satisfying
	\begin{equation}\label{eq-Bowen-prop5-proof-6}
		|\Gamma_n|\geq\exp\frac{nq(1+\tau+2\theta^*)\gamma}{q-1}.
	\end{equation}
	Arbitrarily given $0<\theta<\theta_0$.
	Since $M(n,\varepsilon)$ is non-decreasing as $\varepsilon$ decreases, we have
	\begin{equation*}
		\tau=\sup_{\varepsilon>0}\liminf_{n\rightarrow\infty}\frac{M(n,\varepsilon)}{n}=\lim_{\varepsilon\rightarrow0}\liminf_{n\rightarrow\infty}\frac{M(n,\varepsilon)}{n}.
	\end{equation*}
	Therefore we can take $\varepsilon_0>0$ sufficiently small such that
	\begin{itemize}
		\item $\varepsilon_0<\varepsilon^*/3$;
		\item if two points $x,y\in X$ satisfy $d(x,y)\leq\varepsilon_0$, then $|\varphi(x)-\varphi(y)|\leq\eta$;
		\item for any $0<\varepsilon\leq\varepsilon_0$, one has
		\begin{equation}\label{eq-Bowen-prop5-proof-7}
			\left|\tau-\liminf_{n\rightarrow\infty}\frac{M(n,\varepsilon)}{n}\right|<\theta.
		\end{equation}
	\end{itemize}
	Let
	\begin{equation*}
		\mathfrak{N}=\Bigg\{\lambda\in\mathcal{M}(X):\bigg|\int_X\varphi\mathrm{d}\lambda-\int_X\varphi\mathrm{d}\nu_0\bigg|<\eta\Bigg\}.
	\end{equation*}
	Then $\mathfrak{N}$ is an open neighborhood in $\mathcal{M}(X)$ of $\nu_0$.
	Since $\nu_0$ is ergodic, it follows from Birkhoff's ergodic theorem that we can take $n^*\geq N^*$ sufficiently large such that $X_{n,\mathfrak{N}}\neq\varnothing$ holds for all $n\geq n^*$.
	By \eqref{eq-Bowen-prop5-proof-7}, we can take $n_0\geq n^*$ such that
	\begin{equation*}
		\bigg|\tau-\frac{M(n_0,\varepsilon_0)}{n_0}\bigg|<2\theta.
	\end{equation*}
	Let $m_0=M(n_0,\varepsilon_0)$, then
	\begin{equation}\label{eq-Bowen-prop5-proof-8}
		(\tau-2\theta)n_0\leq m_0\leq (\tau+2\theta)n_0.
	\end{equation}
	For $j\in\mathbb{N}_0$, let $a(j)=j(n_0+m_0)$ and $b(j)=j(n_0+m_0)+n_0-1$.

	Let $\Gamma=\Gamma_{n_0}$ be a $(\delta^*,n_0,\varepsilon^*)$-separated set contained in $X_{n_0,\mathfrak{M}}$ that satisfies \eqref{eq-Bowen-prop5-proof-6} and fix $y\in X_{n_0,\mathfrak{N}}$.
	Define $\Omega=\{y\}\times\Gamma^{q-1}$.
	For any $\xi=(\xi_0,\xi_1,\cdots)\in\Omega^{\mathbb{N}_0}$, we write $\xi_i=(x_i^0,\cdots,x_i^{q-1})\in\Omega$ for each $i\in\mathbb{N}_0$ and define
	\begin{equation*}
		E^{\xi}_k=\bigcap_{j=0}^{k-1}\bigcap_{s=0}^{q-1}\bigcap_{i=a(jq+s)}^{b(jq+s)}f^{-i}\Big(\overline{B(f^{i-a(jq+s)}x_j^s,\varepsilon_0)}\Big).
	\end{equation*}
	From the non-uniform specification property, it follows that each $E^{\xi}_k$ is a non-empty and closed set.
	Similarly, for any $l\in\mathbb{N}$, we can define $E^{\zeta}_k$ for $\zeta\in\Omega^l$ and $1\leq k\leq l$.
	For $l\in\mathbb{N}$, define
	\begin{equation*}
		F_l=\bigcup_{\zeta\in\Omega^l}E^{\zeta}_l
	\end{equation*}
	and $F=\bigcap_{l\geq1}F_l$.
	Then $F_l$ and $F$ are non-empty closed sets.
	Let $E^{\xi}=\bigcap_{k\geq1}E^{\xi}_k$, then $E^{\xi}$ is non-empty and closed.
	It is easy to see that $F=\bigcup_{\xi\in\Omega^{\mathbb{N}_0}}E^{\xi}$.
	Note that if $x\in E^{\xi}$, then for any $i\in\mathbb{N}$, we have $f^{iq(n_0+m_0)}x\in E^{\sigma^i\xi}$, where $\sigma:\Omega^{\mathbb{N}_0}\rightarrow\Omega^{\mathbb{N}_0}$ is the shift map defined by $(\sigma\xi)_k=\xi_{k+1}$.
	As a consequence, $F$ is a non-empty closed $f^{q(n_0+m_0)}$-invariant set.

	Now we find a closed $f$-invariant set $Y\subset X$ satisfying
	\begin{equation*}
		h_{top}(f|_Y)>\gamma\quad\text{and}\quad\mu\notin\mathcal{M}_f^{erg}(Y).
	\end{equation*}
	Choose $L\in\mathbb{N}$ sufficiently large such that $L>\eta^{-1}\|\varphi\|$.
	Let $r_0=Lq(n_0+m_0)$.
	Define
	\begin{equation*}
		\mathfrak{L}=\Bigg\{\lambda\in\mathcal{M}(X):\int_X\varphi\mathrm{d}\lambda\leq\int_X\varphi\mathrm{d}\mu-2\eta\Bigg\}
	\end{equation*}
	and
	\begin{equation*}
		Y=\big\{x\in X:f^jx\in X_{r_0,\mathfrak{L}},\forall j\in\mathbb{N}_0\big\}.
	\end{equation*}
	Clearly $Y$ is a closed $f$-invariant set.
	We claim that $F\subset Y$.
	Since $F$ is $f^{q(n_0+m_0)}$-invariant, it suffices to show that for any $x\in F$ and any $0\leq j\leq q(n_0+m_0)-1$, we have
	\begin{equation*}
		\frac{1}{r_0}\sum_{i=j}^{j+r_0-1}\varphi(f^ix)<\int_X\varphi\mathrm{d}\mu-2\eta.
	\end{equation*}
	Given $x\in F$ and $0\leq j\leq q(n_0+m_0)-1$.
	Fix $\xi\in\Omega^{\mathbb{N}_0}$ such that $x$ lies in $E^{\xi}$.
	Write $\xi_i=(x_i^0,\cdots,x_i^{q-1})$, where $x_i^0=y$ and $x_i^1,\cdots,x_i^{q-1}\in\Gamma$.
	By the assumption that $\inf_{x\in X}\varphi(x)=0$, we have
	\begin{equation}\label{eq-Bowen-prop5-proof-9}
		\begin{split}
			\frac{1}{r_0}\sum_{i=j}^{j+r_0-1}\varphi(f^ix)&\leq\frac{1}{r_0}\sum_{i=0}^{r_0+q(n_0+m_0)-1}\varphi(f^ix)\\
			&\leq\frac{1}{r_0}\sum_{i=0}^{r_0-1}\varphi(f^ix)+\frac{1}{r_0}\sum_{i=r_0}^{r_0+q(n_0+m_0)-1}\varphi(f^ix)\\
			&\leq\frac{1}{r_0}\sum_{i=0}^{r_0-1}\varphi(f^ix)+\frac{\|\varphi\|}{L}<\frac{1}{r_0}\sum_{i=0}^{r_0-1}\varphi(f^ix)+\eta.
		\end{split}
	\end{equation}
	By the definition of $E^{\xi}$, since $|\varphi(x_1)-\varphi(x_2)|\leq\eta$ whenever $d(x_1,x_2)\leq\varepsilon_0$, for any $0\leq s\leq L-1$, any $0\leq t\leq q-1$ and any $a(sq+t)\leq i\leq b(sq+t)$, we have
	\begin{equation}\label{eq-Bowen-prop5-proof-10}
		\varphi(f^ix)\leq\varphi(f^{i-a(sq+t)}x_s^t)+\eta.
	\end{equation}
	Let
	\begin{equation*}
		\Lambda=\mathbb{Z}\cap\bigg([0,r_0-1]\setminus\bigcup_{s=0}^{L-1}\bigcup_{t=0}^{q-1}[a(sq+t),b(sq+t)]\bigg),
	\end{equation*}
	then $|\Lambda|=Lqm_0$.
	It follows from \eqref{eq-Bowen-prop5-proof-10} that
	\begin{equation}\label{eq-Bowen-prop5-proof-11}
		\begin{split}
			\frac{1}{r_0}\sum_{i=0}^{r_0-1}\varphi(f^ix)
			&=\frac{1}{Lq(n_0+m_0)}\Bigg(\sum_{s=0}^{L-1}\sum_{t=0}^{q-1}\sum_{i=a(sq+t)}^{b(sq+t)}\varphi(f^ix)+\sum_{i\in\Lambda}\varphi(f^ix)\Bigg)\\
			&\leq\frac{1}{Lq(n_0+m_0)}\Bigg(\sum_{s=0}^{L-1}\sum_{t=0}^{q-1}\sum_{i=a(sq+t)}^{b(sq+t)}\big(\varphi(f^{i-a(sq+t)}x_s^t)+\eta\big)+Lqm_0\|\varphi\|\Bigg)\\
			&<\frac{1}{Lq(n_0+m_0)}\sum_{s=0}^{L-1}\sum_{t=0}^{q-1}\sum_{i=0}^{n_0-1}\varphi(f^ix_s^t)+\frac{m_0}{n_0+m_0}\|\varphi\|+\eta.
		\end{split}
	\end{equation}
	For any $0\leq s\leq L-1$, since $x_s^t\in\Gamma\subset X_{n_0,\mathfrak{M}}$ for $1\leq t\leq q-1$ and $x_s^0=y\in X_{n_0,\mathfrak{N}}$, we have
	\begin{equation}\label{eq-Bowen-prop5-proof-12}
		\frac{1}{qn_0}\sum_{t=0}^{q-1}\sum_{i=0}^{n_0-1}\varphi(f^ix_s^t)<\frac{q-1}{q}\int_X\varphi\mathrm{d}\mu+\frac{1}{q}\int_X\varphi\mathrm{d}\nu_0+\eta=\int_X\varphi\mathrm{d}\mu_0+\eta.
	\end{equation}
	Combining \eqref{eq-Bowen-prop5-proof-5}, \eqref{eq-Bowen-prop5-proof-8}, \eqref{eq-Bowen-prop5-proof-9}, \eqref{eq-Bowen-prop5-proof-11} and \eqref{eq-Bowen-prop5-proof-12}, we conclude that
	\begin{equation*}
		\begin{split}
			\frac{1}{r_0}\sum_{i=j}^{j+r_0-1}\varphi(f^ix)&<\frac{n_0}{n_0+m_0}\int_X\varphi\mathrm{d}\mu_0+\frac{m_0}{n_0+m_0}\|\varphi\|+3\eta\\
			&\leq\frac{1}{1+\tau-2\theta}\int_X\varphi\mathrm{d}\mu_0+\frac{\tau+2\theta}{1+\tau+2\theta}\|\varphi\|+3\eta<\int_X\varphi\mathrm{d}\mu-2\eta,
		\end{split}
	\end{equation*}
	which implies that $f^jx\in X_{r_0,\mathfrak{L}}$.
	Therefore $F\subset Y$.

	Now we show that $\mu\notin\mathcal{M}_f^{erg}(Y)$.
	Otherwise we can choose a point $x\in Y$ such that $\delta_x^n\rightarrow\mu$ as $n\rightarrow\infty$.
	Then we can take $N\in\mathbb{N}$ sufficiently large such that for any $k\geq N$, we have
	\begin{equation}\label{eq-Bowen-prop5-proof-13}
		\frac1k\sum_{i=0}^{k-1}\varphi(f^ix)>\int_X\varphi\mathrm{d}\mu-\eta.
	\end{equation}
	Since $x\in Y$, for any $j\in\mathbb{N}$, we have
	\begin{equation*}
		\frac{1}{r_0}\sum_{i=(j-1)r_0}^{jr_0-1}\varphi(f^ix)\leq\int_X\varphi\mathrm{d}\mu-2\eta.
	\end{equation*}
	Therefore
	\begin{equation*}
		\frac{1}{2Nr_0}\sum_{i=0}^{2Nr_0-1}\varphi(f^ix)=\frac{1}{2N}\sum_{j=1}^{2N}\bigg(\frac{1}{r_0}\sum_{i=(j-1)r_0}^{jr_0-1}\varphi(f^ix)\bigg)\leq\int_X\varphi\mathrm{d}\mu-2\eta,
	\end{equation*}
	which contradicts with \eqref{eq-Bowen-prop5-proof-13}.

	It remains to show that $h_{top}(f|_Y)\geq\gamma$.
	For $l\in\mathbb{N}$, let $s_l(\kappa,Y)$ denote the maximal cardinality of $(l,\kappa)$-separated sets in $Y$.
	We will show that
	\begin{equation*}
		s_{lq(n_0+m_0)}(\varepsilon_0,Y)\geq\exp(ln_0q(1+\tau+2\theta^*)\gamma).
	\end{equation*}
	For any $\xi=(\xi_0,\ldots,\xi_{l-1})\in\Omega^l$, choose $\tilde{\xi}\in\Omega^{\mathbb{N}_0}$ satisfying $\xi=(\tilde{\xi}_0,\ldots,\tilde{\xi}_{l-1})$, and fix $x_{\xi}\in E^{\tilde{\xi}}\subset F\subset Y$.
	Then $x_{\xi}\in E^{\tilde{\xi}}_l=E^{\xi}_l$.
	For $\xi,\zeta\in\Omega^l$, if $\xi_s\neq\zeta_s$ for some $0\leq s\leq l-1$, then we write $\xi_s=(x_s^0,\cdots,x_s^{q-1})$ and $\zeta_s=(z_s^0,\cdots,z_s^{q-1})$, and assume that $x_s^t\neq z_s^t$ for some $0\leq t\leq q-1$.
	Since $x_s^0=z_s^0=y$, we have $1\leq t\leq q-1$ and $x_s^t,z_s^t\in\Gamma$.
	From the definition of $E^{\xi}_l$, it follows that for any $a(sq+t)\leq i\leq b(sq+t)$, we have
	\begin{equation}\label{eq-Bowen-prop5-proof-14}
		\begin{split}
			d(f^ix_{\xi},f^ix_{\zeta})&\geq d(f^{i-a(sq+t)}x_s^t,f^{i-a(sq+t)}z_s^t)-\Big(d(f^{i-a(sq+t)}x_s^t,f^ix_{\xi})+d(f^ix_{\zeta},f^{i-a(sq+t)}z_s^t)\Big)\\
			&\geq d(f^{i-a(sq+t)}x_s^t,f^{i-a(sq+t)}z_s^t)-2\varepsilon_0.
		\end{split}
	\end{equation}
	Since $\Gamma$ is $(\delta^*,n_0,\varepsilon^*)$-separated, there is an $i\in\{a(sq+t),\cdots,b(sq+t)\}$ such that
	\begin{equation}\label{eq-Bowen-prop5-proof-15}
		d(f^{i-a(sq+t)}x_s^t,f^{i-a(sq+t)}z_s^t)>\varepsilon^*>3\varepsilon_0.
	\end{equation}
	Combining \eqref{eq-Bowen-prop5-proof-14} and \eqref{eq-Bowen-prop5-proof-15}, we obtain that
	\begin{equation*}
		d(f^ix_{\xi},f^ix_{\zeta})>\varepsilon_0.
	\end{equation*}
	Hence $\{x_{\xi}:\xi\in\Omega^l\}$ is an $(lq(m_0+n_0),\varepsilon_0)$-separated set in $Y$.
	It follows from \eqref{eq-Bowen-prop5-proof-6} that
	\begin{equation}\label{eq-Bowen-prop5-proof-16}
		s_{lq(n_0+m_0)}(\varepsilon_0,Y)\geq|\Omega^l|=|\Gamma|^{(q-1)l}\geq\exp(ln_0q(1+\tau+2\theta^*)\gamma).
	\end{equation}
	Combining \eqref{eq-Bowen-prop5-proof-8} and \eqref{eq-Bowen-prop5-proof-16}, we conclude that
	\begin{align*}
		h_{top}(f|_Y)&\geq\limsup_{l\rightarrow\infty}\frac{1}{lq(n_0+m_0)}\log s_{lq(n_0+m_0)}(\varepsilon_0,Y)\\
		&\geq\frac{n_0(1+\tau+2\theta^*)}{n_0+m_0}\gamma\geq\frac{1+\tau+2\theta^*}{1+\tau+2\theta}\gamma>\gamma.
	\end{align*}

	Now we consider the subsystem $(Y,f|_Y)$.
	By the variational principle, there is a $\nu\in\mathcal{M}_f^{erg}(Y)$ such that $h_{\nu}(f|_{Y})>\gamma$.
	Since $\mu\notin\mathcal{M}_f(Y)$, we have $\nu\neq\mu$.
	View $\nu$ as a Borel probability measure on $X$ by setting $\nu(E)=\nu(E\cap Y)$ for every $E\in\mathcal{B}(X)$.
	Then $\nu\in\mathcal{M}_f^{erg}(X)$ and $h_{\nu}(f)>\gamma$.
	This completes the proof of Proposition \ref{Bowen-prop5}.
\end{proof}

%cor6
\begin{corollary}\label{Bowen-cor6}
	Suppose that $(X,f)$ satisfies the non-uniform specification property with gap function $M(n,\varepsilon)$.
	Then
	\begin{equation*}
		\sup_{\mu,\nu\in\mathcal{M}_f^{erg}(X),\mu\neq\nu}\min\{h_{\mu}(f),h_{\nu}(f)\}\geq\frac{1}{1+\tau(M)}h_{top}(f).
	\end{equation*}
\end{corollary}

\begin{remark}\label{Bowen-rmk7}
	The case $m=1$ in Theorem \ref{thm-Bowen-upper-1}, together with Theorem \ref{thm-intermediate-1}, shows that the factor $1/(1+\tau)$ in Proposition \ref{Bowen-prop5} and Corollary \ref{Bowen-cor6} is optimal.
\end{remark}

\begin{proof}[{\bf Proof of Theorem \ref{thm-Bowen-lower-1} (3)}]
	We proceed to show statement (3) of Theorem \ref{thm-Bowen-lower-1}.
	From statement (2) and Corollary \ref{Bowen-cor6}, it follows that
	\begin{equation*}
		h_{top}^B(f,\mathrm{IR}(f))\geq\frac{1-\tau}{1+\tau}h_{top}(f)+\frac{\tau}{(1+\tau)^2}h_{top}(f)=\Big(\frac{1}{1+\tau}-\frac{\tau^2}{(1+\tau)^2}\Big)h_{top}(f).
	\end{equation*}
	This completes the proof of Theorem \ref{thm-Bowen-lower-1}.
\end{proof}

\begin{proof}[{\bf Proof of Corollary \ref{cor-Bowen-lower-1}}]
	From Theorem \ref{thm-Bowen-lower-1} (2), it follows that
	\begin{equation*}
		h_{top}^B(f,\mathrm{IR}(f))\geq\frac{1-\tau}{1+\tau}h_{top}(f)+\frac{\tau}{1+\tau}\sup_{\mu,\nu\in\mathcal{M}_f^{erg}(X),\mu\neq\nu}\min\{h_{\mu}(f),h_{\nu}(f)\}.
	\end{equation*}
	By \cite[Theorem 3.2 (8)]{Lin-Tian-Yu-2024}, $(X,f)$ is not uniquely ergodic.
	By the assumption and the variational principle, for any $0<\gamma<h_{top}(f)$, there are two distinct ergodic measures $\mu,\nu\in\mathcal{M}_f^{erg}(X)$ with $h_{\mu}(f)\geq h_{\nu}(f)>\gamma$.
	Therefore
	\begin{equation*}
		h_{top}^B(f,\mathrm{IR}(f))\geq\frac{1-\tau}{1+\tau}h_{top}(f)+\frac{\tau}{1+\tau}\min\{h_{\mu}(f),h_{\nu}(f)\}\geq\frac{1}{1+\tau}\gamma.
	\end{equation*}
	By the arbitrariness of $\gamma$, we conclude that
	\begin{equation*}
		h_{top}^B(f,\mathrm{IR}(f))\geq\frac{1}{1+\tau}h_{top}(f).
	\end{equation*}
	This completes the proof of Corollary \ref{cor-Bowen-lower-1}.
\end{proof}

\section{Examples with non-uniform specification: upper bound of Bowen topological entropy}\label{Sect-Bowen-upper}
In this section, we will construct dynamical systems that satisfy the non-uniform specification property in proof of Theorem \ref{thm-construction}, and use these examples to prove Theorem \ref{thm-Bowen-upper-1} and Theorem \ref{thm-Bowen-upper-2}.

\subsection{Construction of examples in symbolic dynamics}

The idea of the construction in Theorem \ref{thm-construction} is from the examples in the proof of \cite[Theorem 1.1]{Pavlov-2016} and \cite[Theorem A]{Lin-Tian-Yu-2024}.
%We prove Theorem \ref{thm-construction} with some technical details omitted.

\begin{proof}[{\bf Proof of Theorem \ref{thm-construction}}]
	Assume that $\mathcal{A}=\{1,2,\cdots,A\}$.
	Let $\widetilde{\mathcal{A}}=\{0\}\cup\mathcal{A}=\{0,1,2,\cdots,A\}$.

	First, we show the case $0<\tau<\infty$.
	Consider the set of forbidden words $\mathcal{F}_{\tau}\subset\widetilde{\mathcal{A}}^*$, consisting of:
	\begin{itemize}
		\item all words $w\in\mathcal{A}^*\setminus\mathcal{L}(Z)$;
		\item all words of the form $x_1\cdots x_s0^tx_{s+1}$, where $s,t\in\mathbb{N}$, $x_i\neq0$, and $t<\tau s$.
	\end{itemize}
	Define
	\begin{equation*}
		X=X(\mathcal{F}_{\tau})=\{x\in\widetilde{\mathcal{A}}^{\mathbb{N}_0}:x_ix_{i+1}\cdots x_l\notin\mathcal{F}_{\tau},\;\forall 0\leq i\leq l\}.
	\end{equation*}
	Then $(X,f)$ is a subshift of $(\widetilde{\mathcal{A}}^{\mathbb{N}_0},f)$.
	A typical point in $X$ is of the form $0^{t_0}w^{(1)}0^{t_1}w^{(2)}0^{t_2}\cdots$, where $w^{(i)}\in\mathcal{L}(Z)$ and $t_i\geq\lceil\tau|w^{(i)}|\rceil$ for $i\geq1$.
	Clearly $Z$ is a subshift of $X$.
	This shows (2).

	Let
	\begin{equation*}
		M(n,\varepsilon)=
		\begin{cases}
			1 & \text{if}\;\varepsilon\geq2,\\
			\big\lceil\tau\cdot(n+\lceil-\log_2\varepsilon\rceil)\big\rceil+\lceil-\log_2\varepsilon\rceil+1 & \text{if}\;0<\varepsilon<2.
		\end{cases}
	\end{equation*}
	Clearly, the function $M$ is non-decreasing in $n$ and non-increasing in $\varepsilon$, and
	\begin{equation*}
		\lim_{n\to\infty}\frac{M(n,\varepsilon)}n=\tau,\quad\forall\,\varepsilon\in(0,2).
	\end{equation*}
	
	Now we verify the non-uniform periodic specification property.
	It suffices to consider $0<\varepsilon<2$.
	Given $0<\varepsilon<2$, write $r=\lceil-\log_2\varepsilon\rceil$.
	Given $k\geq 2$ and $k$ points $x^{(1)},\ldots,x^{(k)}\in X$, write
	\begin{equation*}
		x^{(i)}=x^{(i)}_0x^{(i)}_1\cdots\in X,\quad\forall\,1\leq i\leq k.
	\end{equation*}
	Let $\{a_i\}_{i=1}^k$, $\{b_i\}_{i=1}^k$ be two sequences of non-negative integers with $a_1\leq b_1<a_2\leq b_2\cdots<a_k\leq b_k$ and
	\begin{equation*}
		a_{i+1}-b_i\geq M(b_i-a_i+1,\varepsilon),\quad\forall\,1\leq i\leq k-1.
	\end{equation*}
	Fix a positive integer $p$ with $p\geq b_k-a_1+M(b_k-a_k+1,\varepsilon)$.
	Write $n_i=b_i-a_i+1$.
	Let $t_i:=a_{i+1}-b_i-r-1$ for $1\leq i\leq k-1$ and $t_k:=p-b_k+a_1-r-1$.
	Then
	\begin{equation*}
		t_i\geq\lceil\tau(n_i+r)\rceil\;\text{for}\;1\leq i\leq k-1\quad\text{and}\quad t_k\geq\lceil\tau(n_k+r)\rceil.
	\end{equation*}
	Let $w^{(i)}:=x^{(i)}_0x^{(i)}_1\cdots x^{(i)}_{n_i+r-1}\in\mathcal{L}_{n_i+r}(X)$ for each $1\leq i\leq k$ and
	\begin{equation*}
		z'=(w^{(1)}0^{t_1}w^{(2)}0^{t_2}\cdots w^{(k)}0^{t_k})^{\infty}.
	\end{equation*}
	Then $z'\in X$ and $f^p(z')=z'$.
	Let $q_0\in\{0,\ldots,p-1\}$ be the unique integer such that $q_0+a_1\equiv0\pmod p$, and set $z=f^{q_0}(z')$.
	Then
	\begin{equation*}
		d(f^{n-a_i}x^{(i)},f^nz)\leq\sum_{s=r+1}^{\infty}2^{-s}=2^{-r}\leq\varepsilon,\quad\text{for}\quad a_i\leq n\leq b_i,1\leq i\leq k.
	\end{equation*}
	Hence $(X,f)$ satisfies the non-uniform periodic specification property with gap function $M(n,\varepsilon)$.

	It follows from the definition of $M(n,\varepsilon)$ and $(X,f)$ that $\tau_*(X,f)=\tau$.
	Indeed, consider another gap function $M'\in\mathscr{G}(X,f)$ and fix $\varepsilon_0=1/4$.
	Given $n\in\mathbb{N}$, choose two points $x=x_0x_1\cdots\in Z$ and $y=0y_1y_2\cdots$, where $f(y)=y_1y_2\cdots\in Z$.
	By non-uniform specification, there exists $z\in X$ satisfying
	\begin{equation*}
		d(f^kx,f^kz)\leq\varepsilon_0,\;\text{for}\; 0\leq k\leq n-1,\quad\text{and}\quad d(y,f^Tz)\leq\varepsilon_0,
	\end{equation*}
	where $T=n-1+M'(n,\varepsilon_0)$. Since $\varepsilon_0<1$, the first set of inequalities implies that $z_0,\ldots,z_{n-1}$ are all non-zero. Since $\varepsilon_0<1/2$, the second inequality gives $z_T=0$ and $z_{T+1}=y_1\neq0$.
	Let $s\in[n,T]$ be the first zero coordinate of $z$, and let $t\geq1$ be the length of the zero block beginning at $s$. Then $s+t\leq T+1$. By the definition of the forbidden words, $t\geq\lceil\tau s\rceil$. Consequently,
	\begin{equation*}
		M'(n,\varepsilon_0)=T-n+1\geq s+t-n\geq t
		\geq\lceil\tau s\rceil\geq\lceil\tau n\rceil.
	\end{equation*}
	This implies that
	\begin{equation*}
		\tau(M')=\sup_{\varepsilon>0}\liminf_{n\rightarrow\infty}\frac{M'(n,\varepsilon)}{n}
		\geq\liminf_{n\rightarrow\infty}\frac{M'(n,\varepsilon_0)}{n}\geq\tau.
	\end{equation*}
	Hence $\tau_*(X,f)=\tau$.
	This completes the proof of the non-uniform periodic $\tau$-specification property.
	
	We proceed to show (3).
	Let $\mathcal{M}_0=\mathcal{M}_f^{erg}(X)\setminus\mathcal{M}_f^{erg}(Z)$.
%	Similar to the proof of \cite[Proposition 5.1]{Lin-Tian-Yu-2024}, one can show that
%	\begin{equation*}
%		\overline{\mathcal{M}_0}\subset\{\mu\in\mathcal{M}_f(X):G_{\mu}\neq\varnothing\}.
%	\end{equation*}
	Let $\chi_{[0]}$ denote the characteristic function of the cylinder set $[0]$.
	Then $\chi_{[0]}$ is continuous.
	We claim that
	\begin{equation}\label{eq-thm-construction-proof-1}
		\int_X\chi_{[0]}\mathrm{d}\mu\geq\frac{\tau}{1+\tau}>0,\quad\forall\mu\in\overline{\mathcal{M}_0}.
	\end{equation}
	Given $\mu\in \mathcal{M}_0$, let $x\in G_{\mu}$ and write $x=x_0x_1\cdots$.
	If $\mu=\delta_{0^{\infty}}$, then $\int_X\chi_{[0]}\mathrm{d}\mu=1$, and hence \eqref{eq-thm-construction-proof-1} holds.
	Now we assume that $\mu\neq\delta_{0^{\infty}}$.
	Since $\mu\notin\mathcal{M}_f(Z)$, by the construction of forbidden words, $x$ must be of the form
	\begin{equation*}
		x=0^{t_0}w^{(1)}0^{t_1}\cdots w^{(k)}0^{t_k}\cdots,
	\end{equation*}
	where, for each $i\in\mathbb{N}$, $w^{(i)}$ is a word over $Z$ and $t_i$ is a positive integer with $t_i\geq\tau\cdot|w^{(i)}|$.
	Hence
	\begin{equation*}
		\int_X\chi_{[0]}\mathrm{d}\mu=\lim_{k\to\infty}\frac{\sum_{i=0}^kt_i}{t_0+\sum_{i=1}^k(|w^{(i)}|+t_i)}\geq\frac{\tau}{1+\tau}.
	\end{equation*}
	This shows \eqref{eq-thm-construction-proof-1}.
	Note that
	\begin{equation}\label{eq-thm-construction-proof-2}
		\int_X\chi_{[0]}\mathrm{d}\mu=0,\quad\forall\mu\in\mathcal{M}_f(Z).
	\end{equation}
	Combining \eqref{eq-thm-construction-proof-1} and \eqref{eq-thm-construction-proof-2}, we conclude that $\overline{\mathcal{M}_0}$ and $\overline{\mathcal{M}_f^{erg}(Z)}$ are two disjoint closed sets.
	This completes the proof of (3).

	Now we show (4).
	Clearly, for any $x=x_0x_1\cdots\in \bigcup_{j\geq0}f^{-j}Z$, there are only finite many $0$s in $x_0,x_1,\cdots$.
	Thus each $x\in X\setminus \bigcup_{j\geq0}f^{-j}Z$ is either of the form $x=w0^{\infty}$ for some $w\in\mathcal{L}(X)$ or
	\begin{equation*}
		x=0^{t_0}w^{(1)}0^{t_1}w^{(2)}0^{t_2}\cdots,
	\end{equation*}
	where $w^{(i)}\in\mathcal{L}(Z)$, $|w^{(i)}|\geq1$, and $t_i\geq\lceil\tau|w^{(i)}|\rceil$ for $i\geq1$.
	The first case follows immediately from the fact that $\delta_{0^\infty}\in V_f(x)\cap\mathcal{M}_0$.
	In the second case, let $\mu^{(i)}$ be the periodic measure supported on the orbit of periodic point $(w^{(1)}0^{t_1}\cdots w^{(i)}0^{t_i})^{\infty}$, then $\mu^{(i)}\in\mathcal{M}_0$.
	Let $\mu$ be an accumulation point of the sequence $\{\mu^{(i)}\}_{i\geq1}$, then $\mu\in V_f(x)\cap\overline{\mathcal{M}_0}$.
	This completes the proof of (4).

	Finally, we show (5).
	For any ergodic measure $\mu\in\mathcal{M}_0$, we claim that
	\begin{equation}\label{eq-thm-construction-proof-3}
		h_{\mu}(f)\leq\log2+\frac{1}{1+\tau}\log A.
	\end{equation}
	By \eqref{eq-thm-construction-proof-1}, we have $\beta:=\mu([0])\geq\tau/(1+\tau)$.
	For any $0<\kappa<1$ and $n\in\mathbb{N}$, let $\mathcal{W}_{n,\kappa\beta,0}(X)\subset\mathcal{L}_n(X)$ be the set of $n$-letter words with between $(1-\kappa)n\beta$ and $(1+\kappa)n\beta$ occurrences of $0$.
	By Lemma \ref{pre-symbolic-lemma1}, we have
	\begin{equation}\label{eq-thm-construction-proof-4}
		h_{\mu}(f)\leq\liminf_{n\rightarrow\infty}\frac{\log|\mathcal{W}_{n,\kappa\beta,0}(X)|}{n}.
	\end{equation}

	Now we estimate $|\mathcal{W}_{n,\kappa\beta,0}(X)|$.
	Clearly $\mathcal{W}_{n,\kappa\beta,0}(X)\subset\mathcal{W}_{n,\kappa\beta,0}(\widetilde{\mathcal{A}}^{\mathbb{N}_0})$.
	For any $w\in\mathcal{W}_{n,\kappa\beta,0}(\widetilde{\mathcal{A}}^{\mathbb{N}_0})$, there exists $0\leq k\leq n(1-(1-\kappa)\beta)$ such that $w_0,\cdots,w_{n-1}$ have exactly $k$ non-zero letters.
	Note that the positions of $k$ non-zero letters can be chosen in $\binom{n}{k}$ ways, and each non-zero position has $A$ choices for letters.
	Since $\beta\geq\tau/(1+\tau)$, we obtain that
	\begin{equation}\label{eq-thm-construction-proof-5}
		\begin{split}
			|\mathcal{W}_{n,\kappa\beta,0}(X)|\leq|\mathcal{W}_{n,\kappa\beta,0}(\widetilde{\mathcal{A}}^{\mathbb{N}_0})|&\leq\sum_{k=0}^{\lfloor n(1-(1-\kappa)\beta)\rfloor}\binom{n}{k}\cdot A^k\\
			&\leq 2^n A^{n(1-(1-\kappa)\beta)}\leq\Big(2A^{\frac{1+\kappa\tau}{1+\tau}}\Big)^n.
		\end{split}
	\end{equation}
	Combining \eqref{eq-thm-construction-proof-4} and \eqref{eq-thm-construction-proof-5}, we obtain that
	\begin{equation*}
		h_{\mu}(f)\leq\liminf_{n\rightarrow\infty}\frac{\log|\mathcal{W}_{n,\kappa\beta,0}(X)|}{n}\leq\log2+\frac{1+\kappa\tau}{1+\tau}\log A.
	\end{equation*}
	By the arbitrariness of $\kappa$, we conclude that \eqref{eq-thm-construction-proof-3} holds for all $\mu\in\mathcal{M}_0$.

	Now we show that \eqref{eq-thm-construction-proof-3} holds for all $\mu\in\mathcal{M}_f(X)$ satisfying
	\begin{equation}\label{eq-thm-construction-proof-6}
		\int_X\chi_{[0]}\mathrm{d}\mu\geq\frac{\tau}{1+\tau}.
	\end{equation}
	From the ergodic decomposition, there exist $0\leq\alpha\leq 1$ and two probability measures $\pi_0$ and $\pi_Z$ on $\mathcal{M}_0$ and $\mathcal{M}_f^{erg}(Z)$, respectively, such that
	\begin{equation}\label{eq-thm-construction-proof-7}
		\int_X\chi_{[0]}\mathrm{d}\mu=\alpha\int_{\mathcal{M}_0}\int_X\chi_{[0]}\mathrm{d}\omega\mathrm{d}\pi_0(\omega)+(1-\alpha)\int_{\mathcal{M}_f^{erg}(Z)}\int_X\chi_{[0]}\mathrm{d}\omega\mathrm{d}\pi_Z(\omega).
	\end{equation}
	Using \eqref{eq-thm-construction-proof-2}, \eqref{eq-thm-construction-proof-6} and \eqref{eq-thm-construction-proof-7}, we obtain that
	\begin{equation}\label{eq-thm-construction-proof-8}
		\alpha\int_{\mathcal{M}_0}\int_X\chi_{[0]}\mathrm{d}\omega\mathrm{d}\pi_0(\omega)=\int_X\chi_{[0]}\mathrm{d}\mu\geq\frac{\tau}{1+\tau}.
	\end{equation}
	From affine character of the entropy, it follows that
	\begin{equation}\label{eq-thm-construction-proof-9}
		\begin{split}
			h_{\mu}(f)&=\alpha\int_{\mathcal{M}_0}h_{\omega}(f)\mathrm{d}\pi_0(\omega)+(1-\alpha)\int_{\mathcal{M}_f^{erg}(Z)}h_{\omega}(f)\mathrm{d}\pi_Z(\omega)\\
			&\leq\alpha\int_{\mathcal{M}_0}h_{\omega}(f)\mathrm{d}\pi_0(\omega)+(1-\alpha)\cdot h_{top}(Z)\leq\alpha\int_{\mathcal{M}_0}h_{\omega}(f)\mathrm{d}\pi_0(\omega)+(1-\alpha)\cdot\log A.
		\end{split}
	\end{equation}
	For each $\omega\in\mathcal{M}_0$, let $\beta_{\omega}=\int_X\chi_{[0]}\mathrm{d}\omega$.
	Then by \eqref{eq-thm-construction-proof-8}, we have
	\begin{equation}\label{eq-thm-construction-proof-10}
		\alpha\int_{\mathcal{M}_0}\beta_{\omega}\mathrm{d}\pi_0(\omega)\geq\frac{\tau}{1+\tau}.
	\end{equation}
	Similar to \eqref{eq-thm-construction-proof-5}, for $\kappa>0$ sufficiently small, we have
	\begin{equation*}
		|\mathcal{W}_{n,\kappa\beta_{\omega},0}(X)|\leq\Big(2A^{1-(1-\kappa)\beta_{\omega}}\Big)^n.
	\end{equation*}
	Applying Lemma \ref{pre-symbolic-lemma1}, we obtain that $h_{\omega}(f)\leq\log2+(1-(1-\kappa)\beta_{\omega})\cdot \log A$.
	By the arbitrariness of $\kappa$, we obtain that
	\begin{equation}\label{eq-thm-construction-proof-11}
		h_{\omega}(f)\leq\log2+(1-\beta_{\omega})\cdot\log A.
	\end{equation}
	Combining \eqref{eq-thm-construction-proof-9}, \eqref{eq-thm-construction-proof-10} and \eqref{eq-thm-construction-proof-11}, we conclude that
	\begin{equation}
		\begin{split}
			h_{\mu}(f)&\leq(1-\alpha)\cdot\log A+\alpha\bigg(\log2+\log A\cdot\int_{\mathcal{M}_0}(1-\beta_{\omega})\mathrm{d}\pi_0(\omega)\bigg)\\
			&\leq\log2+(1-\alpha)\cdot\log A+\bigg(\alpha-\frac{\tau}{1+\tau}\bigg)\cdot\log A=\log 2+\frac{1}{1+\tau}\log A.
		\end{split}
	\end{equation}
	Hence \eqref{eq-thm-construction-proof-3} holds for all $\mu\in\mathcal{M}_f(X)$ satisfying \eqref{eq-thm-construction-proof-6}.
	Combining with \eqref{eq-thm-construction-proof-1}, we conclude that \eqref{eq-thm-construction-proof-3} holds for all $\mu\in\overline{\mathcal{M}_0}$.
	This completes the proof of (5) for the case $0<\tau<\infty$.

	It remains to consider the case $\tau=\infty$.
	Given $\delta>0$, choose $\sigma>0$ sufficiently large such that
	\begin{equation}\label{eq-thm-construction-proof-12}
		\frac{1}{1+\sigma}<\delta.
	\end{equation}
	Consider the set of forbidden words $\mathcal{F}^{\sigma}_{\infty}\subset\widetilde{\mathcal{A}}^*$, consisting of:
	\begin{itemize}
		\item all words $w\in\mathcal{A}^*\setminus\mathcal{L}(Z)$;
		\item all words of the form $x_1\cdots x_s0^tx_{s+1}$, where $s,t\in\mathbb{N}$, $x_i\neq0$, and $t<\sigma s^2$.
	\end{itemize}
	Define
	\begin{equation*}
		X=X(\mathcal{F}_{\infty}^{\sigma})=\{x\in\widetilde{\mathcal{A}}^{\mathbb{N}_0}:x_ix_{i+1}\cdots x_l\notin\mathcal{F}_{\infty}^{\sigma},\;\forall\,0\leq i\leq l\}.
	\end{equation*}
	Clearly, $Z\subset X$, then (2) holds.
	Let
	\begin{equation*}
		M(n,\varepsilon)=
		\begin{cases}
			1 & \text{if }\varepsilon\geq2,\\
			\big\lceil\sigma(n+\lceil-\log_2\varepsilon\rceil)^2\big\rceil+\lceil-\log_2\varepsilon\rceil+1 & \text{if }0<\varepsilon<2.
		\end{cases}
	\end{equation*}
	Then $M$ is non-decreasing in $n$ and non-increasing in $\varepsilon$, and
	\begin{equation*}
		\lim_{n\to\infty}\frac{M(n,\varepsilon)}{n}=\infty,\quad\forall\,0<\varepsilon<2.
	\end{equation*}
	This shows (1).
	Using arguments analogous to the case $0<\tau<\infty$, one can show that $(X,f)$ satisfies the non-uniform periodic specification property with gap function $M$, and $\tau_*(X,f)=\infty$.
	In particular, $(X,f)$ satisfies the non-uniform periodic $\infty$-specification property.

	Note that a typical point in $X$ is of the form
	\begin{equation*}
		0^{t_0}w^{(1)}0^{t_1}w^{(2)}0^{t_2}\cdots,
	\end{equation*}
	where $w^{(i)}\in\mathcal{L}(Z)$ and
	\begin{equation*}
		t_i\geq\big\lceil\sigma|w^{(i)}|^2\big\rceil\geq\sigma|w^{(i)}|,\quad\forall\,i\geq1.
	\end{equation*}
	Let $\mathcal{M}_0=\mathcal{M}_f^{erg}(X)\setminus\mathcal{M}_f^{erg}(Z)$.
	From the proof of \eqref{eq-thm-construction-proof-1}, with $\sigma$ in place of $\tau$, it follows that
	\begin{equation}\label{eq-thm-construction-proof-infinity-4}
		\int_X\chi_{[0]}\mathrm{d}\mu\geq\frac{\sigma}{1+\sigma}>0,\quad\forall\,\mu\in\overline{\mathcal{M}_0}.
	\end{equation}
	On the other hand, \eqref{eq-thm-construction-proof-2} remains valid for every $\mu\in\mathcal{M}_f(Z)$.
	Thus the proof of (3) is the same as in the case $0<\tau<\infty$.
	The proof of (4) is also similar to the case $0<\tau<\infty$.

	The proof of the entropy estimate in the case $0<\tau<\infty$ uses only the lower bound for $\mu([0])$.
	Applying the same counting argument to \eqref{eq-thm-construction-proof-infinity-4}, with $\sigma$ in place of $\tau$, we obtain that
	\begin{equation*}
		h_{\mu}(f)\leq\log 2+\frac{1}{1+\sigma}\log A,\quad\forall\,\mu\in\overline{\mathcal{M}_0}.
	\end{equation*}
	By \eqref{eq-thm-construction-proof-12}, we conclude that
	\begin{equation*}
		h_{\mu}(f)\leq\log 2+\delta\log|\mathcal{A}|,\quad\forall\,\mu\in\overline{\mathcal{M}_0}.
	\end{equation*}
	This proves (5) and completes the proof.
\end{proof}

\begin{proof}[{\bf Proof of Corollary \ref{cor-entropy}}]
	Since $h>\log\Big(2\big\lceil 2^{\frac{1+\tau}{\tau}}+\frac1\tau\big\rceil^{\frac{1}{1+\tau}}\Big)$, we can choose a positive integer $A$ such that
	\begin{equation*}
		\log 2+\frac{1}{1+\tau}\log A\leq h\leq\log A.
	\end{equation*}
	Consider finite alphabet $\mathcal{A}=\{1,\cdots,A\}$ and a subshift $Z$ over $\mathcal{A}$ such that $h_{top}(Z)=h$.
	Let $(X,f)$ be the subshift constructed in Theorem \ref{thm-construction}, then $(X,f)$ satisfies (1) and the non-uniform $\tau$-specification property.
	By Theorem \ref{thm-construction} (5), for any $\mu\in\mathcal{M}_f^{erg}(X)\setminus\mathcal{M}_f^{erg}(Z)$, we obtain that
	\begin{equation*}
		h_{\mu}(f)\leq\log2+\frac{1}{1+\tau}\log A\leq h.
	\end{equation*}
	From traditional variational principle of topological entropy, it follows that
	\begin{equation*}
		h=h_{top}(Z)=\sup\{h_{\mu}(f):\mu\in\mathcal{M}_f^{erg}(Z)\}.
	\end{equation*}
	Therefore
	\begin{equation*}
		h_{top}(X)=\sup\{h_{\mu}(f):\mu\in\mathcal{M}_f^{erg}(X)\}=\sup\{h_{\mu}(f):\mu\in\mathcal{M}_f^{erg}(Z)\}=h.
	\end{equation*}
	This shows (2) and completes the proof of Corollary \ref{cor-entropy}.
\end{proof}

If $A>\big\lceil 2^{\frac{1+\tau}{\tau}}\big\rceil+1$, then we always have $$\log2+\frac{1}{1+\tau}\log (A-1)>\log\Big(2\big\lceil 2^{\frac{1+\tau}{\tau}}\big\rceil^{\frac{1}{1+\tau}}\Big).$$ Hence, from the above proof and Remark \ref{Remark 1.7}, we have the following corollary.
\begin{corollary}\label{Corollary 4.1}
	Given $0<\tau<\infty$, an alphabet $\mathcal{A}$ with $|\mathcal{A}|>\big\lceil 2^{\frac{1+\tau}{\tau}}\big\rceil+1$ and$$ \log 2+\frac{1}{1+\tau}\log (|\mathcal{A}|-1)\leq h\leq\log (|\mathcal{A}|-1),$$ there exists a subshift $(X,f)$ of $\mathcal{A}^{\mathbb{N}_0}$ (or $\mathcal{A}^\mathbb{Z}$) such that 
	\begin{enumerate}
		\item $(X,f)$ satisfies the non-uniform $\tau$-specification property with gap function $M(n,\varepsilon)$ such that for $\varepsilon>0$ sufficiently small, we have
		\begin{equation*}
			\lim_{n\rightarrow\infty}\frac{M(n,\varepsilon)}{n}=\tau;
		\end{equation*}
		\item $h_{top}(f)=h$.
	\end{enumerate}
\end{corollary}

\subsection{Construction of subsystems in differential dynamical systems}
Let $f:X\to X$ be a homeomorphism.
The system $(X,f)$ is said to be \emph{topologically Anosov} if $X$ is infinite, $f$ is expansive, and $f$ has the two-sided shadowing property (equivalently, the shadowing property for a homeomorphism on a compact metric space).
A standard example is the restriction of a diffeomorphism to an infinite locally maximal hyperbolic set.
The following lemma follows from \cite[Lemma 7.10]{Dong-Hou-Lin-Tian-2025} and its proof.
\begin{lemma}\label{Lemma 4.2}
	Suppose that $(X,f)$ is topologically transitive and topologically Anosov with $h_{top}(f)>0$.
	For any non-empty open set $U\subset X$, any $0<\alpha<h_{top}(f)$, and any $N_0\in\mathbb{N}$, there exist $m,k\in\mathbb{N}$, an alphabet $\mathcal{A}$ with $|\mathcal{A}|=m$, and a non-empty compact $f^k$-invariant set $\Lambda\subset U$ such that $m>N_0$, $\log(m-1)/k>\alpha$, and
	\begin{enumerate}
		\item $f^i(\Lambda)\cap f^j(\Lambda)=\varnothing$ for any $0\leq i<j\leq k-1$;
		\item there is a topological conjugacy $\pi:(\Lambda,f^k)\to(\mathcal{A}^\mathbb{Z},\sigma)$.
	\end{enumerate}
	In particular, $\log(m-1)/k<h_{top}(f)$.
\end{lemma}
\begin{proof}
	By the construction in the proof of \cite[Lemma 7.10]{Dong-Hou-Lin-Tian-2025}, there exist $r,n\in\mathbb{N}$ with $\log r/n>\alpha$ such that, for every sufficiently large $l$, if $m=r^{l-2}$ and $k=nl$, then there is a compact $f^k$-invariant set $\Lambda\subset U$ satisfying (1) and (2).
	Since $m\to\infty$ and $\log(m-1)/k\to\log r/n>\alpha$ as $l\to\infty$, we may choose $l$ sufficiently large so that $m>N_0$ and $\log(m-1)/k>\alpha$.
	Finally, the conjugacy in (2) gives $\log m=h_{top}(f^k|_\Lambda)\leq k h_{top}(f)$, which implies the last assertion.
\end{proof}
Under the assumptions of Corollary \ref{Corollary C New}, \cite[Theorem 1]{Gelfert-2016} and its correction \cite{Gelfert-2022}, together with Lemma \ref{Lemma 4.2}, give the following consequence.
\begin{lemma}\label{Lemma 4.3}
	Suppose that $f$ is a $C^1$ diffeomorphism on a compact $C^\infty$ Riemannian manifold $M$ and $\mu$ is a hyperbolic ergodic measure with $h_{\mu}(f)>0$. Suppose also that one of the following conditions is satisfied:
	\begin{enumerate}
		\item $f$ is a $C^{1+\alpha}$ diffeomorphism for some $\alpha>0$;
		\item $\mu$ admits a dominated splitting corresponding to the stable/unstable subspaces of its Oseledets splitting.
	\end{enumerate}
	Then for any $0<\beta<h_\mu(f)$ and any $\gamma>0$, there exists a topologically transitive topologically Anosov subsystem $(X,f|_X)$ of $(M,f)$ such that
	\begin{equation*}
		\beta<h_{top}(f|_X)<\beta+\gamma.
	\end{equation*}
\end{lemma}
\begin{proof}
	Fix
	\begin{equation*}
		e\in\big(\beta,\min\{\beta+\gamma,h_\mu(f)\}\big).
	\end{equation*}
	By \cite[Theorem 1 and its proof]{Gelfert-2016}, together with the correction \cite{Gelfert-2022}, there exist an integer $p\geq1$, an integer $q\geq2$, and a compact set $\Gamma\subset M$ such that $f^p|_\Gamma$ is conjugate to the full two-sided shift on $q$ symbols and
	\begin{equation*}
		\beta<\frac{\log q}{p}<\beta+\gamma.
	\end{equation*}

	Let $d$ be the smallest positive integer such that $f^d(\Gamma)=\Gamma$.
	Then $d\mid p$.
	Indeed, write $p=ad+b$ with $0\leq b<d$.
	Since $f^d(\Gamma)=\Gamma$ and $f^p(\Gamma)=\Gamma$, we have $f^b(\Gamma)=\Gamma$, and hence $b=0$.
	Then $(f^d|_\Gamma)^{p/d}=f^p|_\Gamma$ is conjugate to the full two-sided shift on $q$ symbols.
	It follows that $(\Gamma,f^d|_\Gamma)$ is topologically transitive and topologically Anosov, and
	\begin{equation*}
		h_{top}(f^d|_\Gamma)=\frac{d}{p}\log q.
	\end{equation*}
	In particular,
	\begin{equation*}
		d\beta<h_{top}(f^d|_\Gamma)<d(\beta+\gamma).
	\end{equation*}

	For $1\leq j<d$, set $A_j=\Gamma\cap f^j(\Gamma)$.
	Then $A_j$ is a proper closed $f^d$-invariant subset of $\Gamma$.
	Indeed, if $A_j=\Gamma$, then
	\begin{equation*}
		\Gamma\subset f^j(\Gamma)\subset\cdots\subset f^{dj}(\Gamma)=\Gamma,
	\end{equation*}
	so $f^j(\Gamma)=\Gamma$, contradicting the choice of $d$.
	Under the conjugacy between $f^p|_\Gamma$ and the full two-sided shift, $A_j$ corresponds to a proper subshift and hence has empty interior in $\Gamma$.
	Therefore
	\begin{equation*}
		U=\Gamma\setminus\bigcup_{j=1}^{d-1}A_j
	\end{equation*}
	is a non-empty $f^d$-invariant open subset of $\Gamma$.

	By Lemma \ref{Lemma 4.2}, there exist integers $s,k\geq1$ and a compact set $\Delta\subset U$ such that $f^{dk}(\Delta)=\Delta$, the sets $f^{ad}(\Delta)$, $0\leq a<k$, are pairwise disjoint, $(\Delta,f^{dk})$ is conjugate to the full two-sided shift on $s$ symbols, and
	\begin{equation*}
		\frac{\log(s-1)}{k}>d\beta,
		\qquad
		\frac{\log s}{k}\leq h_{top}(f^d|_\Gamma)<d(\beta+\gamma).
	\end{equation*}

	Let $N=dk$ and
	\begin{equation*}
		X=\bigcup_{r=0}^{N-1}f^r(\Delta).
	\end{equation*}
	We show that the sets $f^r(\Delta)$, $0\leq r<N$, are pairwise disjoint.
	For $0\leq r,r'<N$, write $r=ad+b$ and $r'=a'd+b'$, where $0\leq a,a'<k$ and $0\leq b,b'<d$.
	If $b=b'$, the pairwise disjointness of the sets $f^{ad}(\Delta)$ implies that $r=r'$ whenever $f^r(\Delta)\cap f^{r'}(\Delta)\neq\varnothing$.
	Suppose that $b<b'$ and $f^r(\Delta)\cap f^{r'}(\Delta)\neq\varnothing$.
	Applying $f^{-b}$ gives
	\begin{equation*}
		f^{ad}(\Delta)\cap f^{b'-b}f^{a'd}(\Delta)\neq\varnothing.
	\end{equation*}
	Since $\Delta\subset U$ and $U$ is $f^d$-invariant, this implies $U\cap A_{b'-b}\neq\varnothing$, a contradiction.

	Since $f^N(\Delta)=\Delta$, the map $f$ cyclically permutes the sets $f^r(\Delta)$, $0\leq r<N$.
	Moreover, $f^N$ restricted to each $f^r(\Delta)$ is conjugate to the full two-sided shift on $s$ symbols.
	Thus $(X,f|_X)$ is topologically transitive.
	The map $(f|_X)^N$ is expansive and has the two-sided shadowing property, and hence so does $f|_X$.
	Since $X$ is infinite, $(X,f|_X)$ is topologically Anosov.
	Finally,
	\begin{equation*}
		h_{top}(f|_X)=\frac{1}{N}\log s=\frac{\log s}{dk}.
	\end{equation*}
	Hence
	\begin{equation*}
		\beta<h_{top}(f|_X)<\beta+\gamma.
	\end{equation*}
\end{proof}

\begin{proof}[{\bf Proof of Corollary \ref{Corollary C New}}]
	Fix $0<\tau<\infty$ and $0<h<h_\mu(f)$, and choose $0<\gamma<\tau h/(\tau+2)$.
	By Lemma \ref{Lemma 4.3}, there exists a topologically transitive topologically Anosov subsystem $(X,f|_X)$ of $(M,f)$ with
	\begin{equation*}
		h<h_{top}(f|_X)<h+\gamma.
	\end{equation*}
	Applying Lemma \ref{Lemma 4.2} to $(X,f|_X)$ with the open set $U=X$, $\alpha=h$, and $N_0=\big\lceil2^{2(1+\tau)/\tau}\big\rceil+1$, we obtain $m,N\in\mathbb{N}$, an alphabet $\mathcal{A}$ with $|\mathcal{A}|=m$, and a non-empty compact $f^N$-invariant set $\Delta\subset X$ such that $m>\big\lceil2^{2(1+\tau)/\tau}\big\rceil+1$, $\log(m-1)/N>h$, and
	\begin{enumerate}
		\item $f^i(\Delta)\cap f^j(\Delta)=\varnothing$ for any $0\leq i<j\leq N-1$;
		\item there is a topological conjugacy $\pi:(\Delta,f^N)\to(\mathcal{A}^\mathbb{Z},\sigma)$.
	\end{enumerate}
	In particular, $\log(m-1)/N<h_{top}(f|_X)<h+\gamma$.
	Put $a=\log(m-1)/N$. Since $m-1>2^{2(1+\tau)/\tau}$, we have
	\begin{equation*}
		\frac{\log2}{\log(m-1)}<\frac{\tau}{2(1+\tau)}.
	\end{equation*}
	As $h<a<h+\gamma$, it follows that
	\begin{equation*}
		\frac{\log2}{N}+\frac{1}{1+\tau}\frac{\log(|\mathcal{A}|-1)}N
		<\frac{\tau+2}{2(1+\tau)}(h+\gamma)<h,
	\end{equation*}
	where the last inequality follows from $\gamma<\tau h/(\tau+2)$.
	Since $Nh<\log(m-1)$, the hypotheses of Corollary \ref{Corollary 4.1} are satisfied.
	By Corollary \ref{Corollary 4.1}, there exists a subshift $(\Gamma,\sigma)$ of $\mathcal{A}^\mathbb{Z}$ such that
	\begin{enumerate}
		\item $(\Gamma,\sigma)$ satisfies the non-uniform $\tau$-specification property with gap function $M(n,\varepsilon)$ such that, for all sufficiently small $\varepsilon>0$,
		\begin{equation*}
			\lim_{n\rightarrow\infty}\frac{M(n,\varepsilon)}{n}=\tau;
		\end{equation*}
		\item $h_{top}(\sigma|_\Gamma)=Nh$.
	\end{enumerate}

	Denote $Y=\pi^{-1}(\Gamma)$.
	By Proposition \ref{intro-prop:conjugacy-invariant}, $(Y,f^N)$ satisfies the non-uniform $\tau$-specification property with a gap function $M_Y(n,\varepsilon)$ such that, for all sufficiently small $\varepsilon>0$,
	\begin{equation*}
		\lim_{n\rightarrow\infty}\frac{M_Y(n,\varepsilon)}{n}=\tau.
	\end{equation*}
	Moreover, $Y$ is a non-empty compact $f^N$-invariant subset of $M$, and
	\begin{equation*}
		h_{top}(f^N|_Y)=h_{top}(\sigma|_\Gamma)=Nh.
	\end{equation*}
	Since $Y\subset\Delta$, we also have
	\begin{equation*}
		f^i(Y)\cap f^j(Y)=\varnothing
	\end{equation*}
	for any $0\leq i<j\leq N-1$.
	Thus $Y$ satisfies (1)--(3).

	To prove the last assertion, let
	\begin{equation*}
		\Lambda=\bigcup_{i=0}^{N-1}f^i(Y)
		\quad\text{and}\quad
		D_i=f^i(Y),\quad 0\leq i<N.
	\end{equation*}
	The sets $D_i$ are pairwise disjoint compact subsets of $\Lambda$ and hence are clopen in $\Lambda$.
	Since $\Gamma$ is a two-sided subshift, $\sigma(\Gamma)=\Gamma$, and the conjugacy gives $f^N(Y)=Y$.
	Therefore $f(D_i)=D_{i+1\ ({\rm mod}\ N)}$, and $f^i:Y\to D_i$ conjugates $f^N|_Y$ to $f^N|_{D_i}$.
	By Proposition \ref{intro-prop:conjugacy-invariant}, every $(D_i,f^N)$ satisfies the non-uniform $\tau$-specification property.
	Thus $\{D_0,\ldots,D_{N-1}\}$ is a regular periodic decomposition with disjoint steps.
	Finally,
	\begin{equation*}
	\begin{split}
		N h_{top}(f|_\Lambda)
		&=h_{top}((f|_\Lambda)^N)
		=\max_{0\leq i<N}h_{top}(f^N|_{D_i})\\
		&=h_{top}(f^N|_Y)=Nh.
	\end{split}
	\end{equation*}
	Therefore $(\Lambda,f|_\Lambda)$ satisfies the relative non-uniform $\tau$-specification property with disjoint steps and $h_{top}(f|_\Lambda)=h$.
\end{proof}

\subsection{Upper bounds for Bowen topological entropy of the irregular set}

The following well-known result is used to estimate the upper bound for Bowen topological entropy.

%prop1
\begin{proposition}[{\cite[Theorem 2]{Bowen-1973}}]\label{Bowen-upper-prop1}
	Let $(X,f)$ be a dynamical system and
	\begin{equation*}
		\mathrm{QR}(t)=\{x\in X: \text{there is a}\;\mu\in V_f(x)\;\text{such that}\;h_{\mu}(f)\leq t\}.
	\end{equation*}	
	Then $h_{top}^B(f,\mathrm{QR}(t))\leq t$.
\end{proposition}

\subsubsection{Proof of Theorem \ref{thm-Bowen-upper-1}}

We prove Theorem \ref{thm-Bowen-upper-1} by using Theorem \ref{thm-construction}.
The subshift $Z$ will be taken as a union of $m$ copies of a strictly ergodic subshift that has sufficiently large topological entropy.

\begin{proof}[{\bf Proof of Theorem \ref{thm-Bowen-upper-1}}]
	From Lemma \ref{pre-symbolic-lemma2}, there are alphabets $\mathcal{A}^{(j)}$ and sequences $z^{(j)}$ over $\mathcal{A}^{(j)}$ for all $j\in\mathbb{N}$ such that
	\begin{itemize}
		\item $|\mathcal{A}^{(j)}|\rightarrow\infty$ as $j\rightarrow\infty$;
		\item for any $j\in\mathbb{N}$, the subshift $(Z^{(j)},f)$ is strictly ergodic, where $Z^{(j)}=\omega(z^{(j)},f)$ is the $\omega$-limit set of $z^{(j)}$;
		\item for any $j\in\mathbb{N}$, we have
		\begin{equation}\label{eq-thm-Bowen-upper-1-proof-1}
			-\delta+\log|\mathcal{A}^{(j)}|<h_{top}(Z^{(j)})<\log|\mathcal{A}^{(j)}|.
		\end{equation}
	\end{itemize}
	Without loss of generality, we assume that $\mathcal{A}^{(j)}=\{1,2,\cdots,A_j\}$.
	For $1\leq k\leq m$, define
	\begin{equation*}
		\mathcal{A}^{(j)}_k=\{(k-1)\cdot A_j+1,\cdots,k\cdot A_j\}
	\end{equation*}
	and
	\begin{equation*}
		Z_k^{(j)}=\{x=(x_0x_1\cdots)\in(\mathcal{A}^{(j)}_k)^{\mathbb{N}_0}:\bar{x}=(x_0-(k-1)A_j,x_1-(k-1)A_j,\cdots)\in Z^{(j)}\}.
	\end{equation*}
	In other words, alphabets $\mathcal{A}^{(j)}_1,\cdots,\mathcal{A}^{(j)}_m$ and subshifts $Z^{(j)}_1,\cdots,Z^{(j)}_m$ are $m$ disjoint copies of alphabet $\mathcal{A}^{(j)}$ and subshift $Z^{(j)}$, respectively, obtained by renaming the symbols.

	First, we consider the case of $0<\tau<\infty$.
	Fix $j$ sufficiently large.
	For simplicity, we rewrite $\mathcal{A}^{(j)}_k$ and $Z^{(j)}_k$ as $\mathcal{A}_k$ and $Z_k$, respectively.
	Consider finite alphabet
	\begin{equation*}
		\mathcal{A}=\bigcup_{k=1}^m\mathcal{A}_k=\{1,\cdots,A_j,A_j+1,\cdots,mA_j\}
	\end{equation*}
	and subshift
	\begin{equation*}
		Z=\bigcup_{k=1}^mZ_k\subset\mathcal{A}^{\mathbb{N}_0}.
	\end{equation*}
	Let $\widetilde{\mathcal{A}}=\{0\}\cup\mathcal{A}=\{0,1,\cdots,mA_j\}$ and $X\subset\widetilde{\mathcal{A}}^{\mathbb{N}_0}$ be the subshift defined in Theorem \ref{thm-construction} associated to $\mathcal{A}$ and $Z$.
	Then $(X,f)$ has non-uniform periodic $\tau$-specification and satisfies (1).

	We proceed to show (2).
	Let $\nu_k$ be the unique ergodic measure of $Z_k$, then $\mathcal{M}_f^{erg}(Z)=\{\nu_1,\cdots,\nu_m\}$.
	By Theorem \ref{thm-construction} (5), for any $\mu\in\mathcal{M}_f^{erg}(X)\setminus\{\nu_1,\cdots,\nu_m\}$, we have
	\begin{equation}\label{eq-thm-Bowen-upper-1-proof-2}
		h_{\mu}(f)\leq\log2+\frac{1}{1+\tau}\log|\mathcal{A}|=\log2+\frac{1}{1+\tau}\log (mA_j).
	\end{equation}
	Choosing $j$ sufficiently large, we may assume that $A_j$ is sufficiently large such that
	\begin{equation}\label{eq-thm-Bowen-upper-1-proof-3}
		\log2+\frac{1}{1+\tau}\log (mA_j)\leq\bigg(\frac{1}{1+\tau}+\delta\bigg)\cdot(\log A_j-\delta).
	\end{equation}
	Combining \eqref{eq-thm-Bowen-upper-1-proof-1}, \eqref{eq-thm-Bowen-upper-1-proof-2} and \eqref{eq-thm-Bowen-upper-1-proof-3}, we obtain that
	\begin{equation*}
		h_{\mu}(f)\leq\bigg(\frac{1}{1+\tau}+\delta\bigg)\cdot h_{top}(Z^{(j)})<h_{top}(Z^{(j)}).
	\end{equation*}
	Note that $h_{\nu_k}(f)=h_{top}(Z^{(j)})$ for $k=1,\cdots,m$.
	Applying traditional variational principle of topological entropy and affine character of entropy, we obtain that
	\begin{equation*}
		h_{top}(X)=\sup\{h_{\mu}(f):\mu\in\mathcal{M}_f^{erg}(X)\}\leq h_{top}(Z^{(j)})\leq h_{top}(X).
	\end{equation*}
	Hence $h_{top}(X)=h_{top}(Z^{(j)})$ and $\nu_1,\nu_2,\cdots,\nu_m$ are $m$ ergodic measures of maximal entropy.
	This completes the proof of (2).

	Condition (3) follows directly from the construction of $X$ and Theorem \ref{thm-construction} (3).
	Indeed, it is clear that $\nu_1,\cdots,\nu_m$ are isolated in $\mathcal{M}_f^{erg}(Z)$.
	It follows from Theorem \ref{thm-construction} (3) that the following two sets $\overline{\mathcal{M}_0}:=\overline{\mathcal{M}_f^{erg}(X)\setminus\{\nu_1,\cdots,\nu_m\}}$ and $\overline{\mathcal{M}_f^{erg}(Z)}:=\{\nu_1,\cdots,\nu_m\}$ are disjoint.
	Hence $\nu_1,\cdots,\nu_m$ are isolated in $\mathcal{M}_f^{erg}(X)$.
	This shows (3).

	Now we show (4).
	Recall that $Z=Z_1\cup\cdots\cup Z_m$.
	Since $Z_k$ is strictly ergodic, we have $\bigcup_{i\geq 0}f^{-i}Z=G_{\nu_1}\cup\cdots\cup G_{\nu_m}$.
	Hence $\mathrm{IR}(f)\subset X\setminus\bigcup_{i\geq 0}f^{-i}Z$.
	Combining \eqref{eq-thm-Bowen-upper-1-proof-1}, \eqref{eq-thm-Bowen-upper-1-proof-2}, \eqref{eq-thm-Bowen-upper-1-proof-3}, Theorem \ref{thm-construction} (4) and (5), and Proposition \ref{Bowen-upper-prop1}, we conclude that
	\begin{equation*}
		h_{top}^B(f,\mathrm{IR}(f))\leq\log2+\frac{1}{1+\tau}\log(mA_j)\leq\bigg(\frac{1}{1+\tau}+\delta\bigg)(\log A_j-\delta)\leq\bigg(\frac{1}{1+\tau}+\delta\bigg)h_{top}(X).
	\end{equation*}
	This completes the proof of (4).
	
	It remains to consider the case $\tau=\infty$.
	Choose $\eta>0$ such that
	\begin{equation*}
		0<\eta<\min\{\delta,1/2\}.
	\end{equation*}
	We may retake $\mathcal{A}^{(j)}$ and $Z^{(j)}$ such that
	\begin{equation}\label{eq-thm-Bowen-upper-1-proof-infinity-1}
		-\eta+\log A_j<h_{top}(Z^{(j)})<\log A_j.
	\end{equation}
	Fix $j$ sufficiently large such that
	\begin{equation}\label{eq-thm-Bowen-upper-1-proof-infinity-2}
		\log2+\eta\log(mA_j)\leq\min\{\delta,1/2\}(\log A_j-\eta).
	\end{equation}
	Write $\mathcal{A}^{(j)}_k$ and $Z^{(j)}_k$ as $\mathcal{A}_k$ and $Z_k$, respectively, and let
	\begin{equation*}
		\mathcal{A}=\bigcup_{k=1}^m\mathcal{A}_k,\qquad Z=\bigcup_{k=1}^mZ_k.
	\end{equation*}
	Let $\widetilde{\mathcal{A}}=\{0\}\cup\mathcal{A}$ and $X\subset \widetilde{\mathcal{A}}^{\mathbb{N}_0}$ be the subshift given by Theorem \ref{thm-construction} for $\tau=\infty$, with $\eta$ in place of $\delta$ in Theorem \ref{thm-construction} (5).
	In other words, $X$ is induced by the set of forbidden words $\mathcal{F}^{\sigma}_{\infty}\subset\widetilde{\mathcal{A}}^*$, consisting of:
	\begin{itemize}
		\item all words $w\in\mathcal{A}^*\setminus\mathcal{L}(Z)$;
		\item all words of the form $x_1\cdots x_s0^tx_{s+1}$, where $s,t\in\mathbb{N}$, $x_i\neq0$, and $t<\sigma s^2$,
	\end{itemize}
	where $\sigma>0$ satisfies
	\begin{equation*}
		\frac{1}{1+\sigma}<\eta.
	\end{equation*}
	Then $(X,f)$ has non-uniform periodic $\infty$-specification and satisfies (1).
	The same argument as in the case $0<\tau<\infty$ shows (2) and (3).

	Note that $\mathrm{IR}(f)\subset X\setminus\bigcup_{i\geq0}f^{-i}Z$.
	Finally, Theorem \ref{thm-construction} (4) and (5), Proposition \ref{Bowen-upper-prop1}, \eqref{eq-thm-Bowen-upper-1-proof-infinity-1}, and \eqref{eq-thm-Bowen-upper-1-proof-infinity-2} imply that
	\begin{equation*}
		h_{top}^B(f,\mathrm{IR}(f))\leq\log2+\eta\log(mA_j)\leq\delta(\log A_j-\eta)\leq\delta h_{top}(X).
	\end{equation*}
	This proves (4) and completes the proof of Theorem \ref{thm-Bowen-upper-1}.
\end{proof}

\subsubsection{Proof of Theorem \ref{thm-Bowen-upper-2}}

We prove Theorem \ref{thm-Bowen-upper-2} by using Theorem \ref{thm-construction}.
The subshift $Z$ will be taken as a (one-sided) full shift over a sufficiently large finite alphabet.

\begin{proof}[{\bf Proof of Theorem \ref{thm-Bowen-upper-2}}]
	First, we consider the case $0<\tau<\infty$.
	Take a positive integer $A$ sufficiently large.
	Consider a finite alphabet $\mathcal{A}=\{1,\cdots,A\}$ and $Z=\mathcal{A}^{\mathbb{N}_0}$.
	Let $\tilde{\mathcal{A}}=\{0\}\cup\mathcal{A}=\{0,1,\cdots,A\}$ and $X\subset\tilde{\mathcal{A}}^{\mathbb{N}_0}$ be the subshift defined in Theorem \ref{thm-construction} associated to $\mathcal{A}$ and $Z$.
	Then $(X,f)$ has non-uniform $\tau$-specification and satisfies (1).

	We proceed to show (2).
	Since $Z=\mathcal{A}^{\mathbb{N}_0}$, we have $h_{top}(Z)=\log A$.
	When $A>2^{\frac{1+\tau}{\tau}}$, it follows from Theorem \ref{thm-construction} (5) that for any $\mu\in\mathcal{M}_f^{erg}(X)\setminus\mathcal{M}_f^{erg}(Z)$, we have
	\begin{equation*}
		h_{\mu}(f)\leq\log2+\frac{1}{1+\tau}\log A\leq\log A=h_{top}(Z).
	\end{equation*}
	Therefore,
	\begin{equation*}
		h_{top}(X)=\sup\{h_{\mu}(f):\mu\in\mathcal{M}_f^{erg}(X)\}\leq h_{top}(Z)\leq h_{top}(X).
	\end{equation*}
	Hence $h_{top}(X)=h_{top}(Z)$.
	Let $\psi=\chi_{[1]}\in C(X)$. Its restriction $\psi|_Z$ has a non-empty irregular set in the full shift $(Z,f|_Z)$. By \cite[Theorem 2.1]{Barreira-Schmeling-2000},
	\begin{equation*}
		h_{top}^B\big(f|_Z,I_{\psi|_Z}(f|_Z)\big)=\log A.
	\end{equation*}
	Since
	\begin{equation*}
		I_{\psi|_Z}(f|_Z)\subset I_{\psi}(f)\subset\mathrm{IR}(f),
	\end{equation*}
	we obtain
	\begin{equation*}
		\log A\leq h_{top}^B(f,I_{\psi}(f))
		\leq h_{top}^B(f,\mathrm{IR}(f))
		\leq h_{top}(f)=\log A.
	\end{equation*}
	Thus $\psi\in\mathfrak{C}(f)$, and this completes the proof of (2).

	Now we show (3).
	Consider continuous function $\varphi=\chi_{[0]}$.
	It is easy to see that the point
	\begin{equation*}
		x:=10^{1+\lceil\tau\rceil}1^{(1+\lceil\tau\rceil)^2}0^{(1+\lceil\tau\rceil)^3}\cdots
	\end{equation*}
	in $X$ satisfies
	\begin{equation*}
		\limsup_{n\to\infty}\frac1n\sum_{k=0}^{n-1}\varphi(f^kx)=\frac{1+\lceil\tau\rceil}{2+\lceil\tau\rceil}\quad\text{and}\quad\liminf_{n\to\infty}\frac{1}{n}\sum_{k=0}^{n-1}\varphi(f^kx)=\frac{1}{2+\lceil\tau\rceil}.
	\end{equation*}
	Hence $x\in I_{\varphi}(f)$.
	In particular, $I_{\varphi}(f)\neq\varnothing$.
	Clearly $I_{\varphi}(f)\subset X\setminus\bigcup_{i\geq 0}f^{-i}Z$.
	By Theorem \ref{thm-construction} (4) and (5), for every $y\in I_{\varphi}(f)$ there exists $\mu_y\in V_f(y)\cap\overline{\mathcal{M}_0}$ such that
	\begin{equation*}
		h_{\mu_y}(f)\leq\log2+\frac{1}{1+\tau}\log A.
	\end{equation*}
	When $A$ is taken sufficiently large, we have
	\begin{equation*}
		h_{\mu_y}(f)\leq\bigg(\frac{1}{1+\tau}+\delta\bigg)\cdot\log A=\bigg(\frac{1}{1+\tau}+\delta\bigg)\cdot h_{top}(f).
	\end{equation*}
	Consequently,
	\begin{equation*}
		I_{\varphi}(f)\subset \mathrm{QR}\bigg(\Big(\frac{1}{1+\tau}+\delta\Big)h_{top}(f)\bigg).
	\end{equation*}
	By Proposition \ref{Bowen-upper-prop1}, we have
	\begin{equation*}
		h_{top}^B(f,\mathrm{CI}(f))\leq h_{top}^B(f,I_{\varphi}(f))\leq\bigg(\frac{1}{1+\tau}+\delta\bigg)\cdot h_{top}(f),
	\end{equation*}
	which completes the proof of (3).

	Finally, we construct $(X,f)$ satisfying (1)--(3) for the case $\tau=\infty$.
	Choose $\eta>0$ such that
	\begin{equation*}
		0<\eta<\min\{\delta,1\}.
	\end{equation*}
	Take a positive integer $A$ sufficiently large such that
	\begin{equation}\label{eq-thm-Bowen-upper-2-proof-infinity-1}
		\log2+\eta\log A\leq\min\{\delta,1\}\log A.
	\end{equation}
	Consider the finite alphabet $\mathcal{A}=\{1,\cdots,A\}$ and the full shift $Z=\mathcal{A}^{\mathbb{N}_0}$.
	Let $\widetilde{\mathcal{A}}=\{0\}\cup\mathcal{A}$ and let
	\begin{equation*}
		X=X(\mathcal{F}^{\sigma}_{\infty})\subset\widetilde{\mathcal{A}}^{\mathbb{N}_0},
	\end{equation*}
	where $\mathcal{F}^{\sigma}_{\infty}\subset\widetilde{\mathcal{A}}^*$ consists of:
	\begin{itemize}
		\item all words $w\in\mathcal{A}^*\setminus\mathcal{L}(Z)$;
		\item all words of the form $x_1\cdots x_s0^tx_{s+1}$, where $s,t\in\mathbb{N}$, $x_i\neq0$, and $t<\sigma s^2$,
	\end{itemize}
	where $\sigma>0$ satisfies
	\begin{equation*}
		\frac{1}{1+\sigma}<\eta.
	\end{equation*}
	Then $(X,f)$ has non-uniform $\infty$-specification and satisfies (1).

	By Theorem \ref{thm-construction} (5) and \eqref{eq-thm-Bowen-upper-2-proof-infinity-1}, for any $\mu\in\mathcal{M}_f^{erg}(X)\setminus\mathcal{M}_f^{erg}(Z)$, we have
	\begin{equation*}
		h_{\mu}(f)\leq\log2+\eta\log A\leq\log A=h_{top}(Z).
	\end{equation*}
	As in the case $0<\tau<\infty$, the variational principle implies that $h_{top}(X)=h_{top}(Z)=\log A$, and for $\psi=\chi_{[1]}\in C(X)$, we have
	\begin{equation*}
		h_{top}^B(f,I_{\psi}(f))=h_{top}^B(f,\mathrm{IR}(f))=h_{top}(f).
	\end{equation*}
	Thus $\psi\in\mathfrak{C}(f)$, and this proves (2).

	Consider $\varphi=\chi_{[0]}$.
	It is easy to see that the point
	\begin{equation*}
		x:=10^{(1+\lceil\sigma\rceil)^3}1^{(1+\lceil\sigma\rceil)^9}0^{(1+\lceil\sigma\rceil)^{27}}1^{(1+\lceil\sigma\rceil)^{81}}0^{(1+\lceil\sigma\rceil)^{243}}\cdots
	\end{equation*}
	in $X$ satisfies
	\begin{equation*}
		\limsup_{n\to\infty}\frac1n\sum_{k=0}^{n-1}\varphi(f^kx)=1\quad\text{and}\quad\liminf_{n\to\infty}\frac{1}{n}\sum_{k=0}^{n-1}\varphi(f^kx)=0.
	\end{equation*}
	Hence $x\in I_{\varphi}(f)$.
	In particular, $I_{\varphi}(f)\neq\varnothing$ and $I_{\varphi}(f)\subset X\setminus\bigcup_{i\geq0}f^{-i}Z$.
	By Theorem \ref{thm-construction} (4) and (5), for every $x\in I_{\varphi}(f)$, there exists $\mu\in V_f(x)$ such that
	\begin{equation*}
		h_{\mu}(f)\leq\log2+\eta\log A\leq\delta\log A=\delta h_{top}(f).
	\end{equation*}
	Applying Proposition \ref{Bowen-upper-prop1}, we conclude that
	\begin{equation*}
		h_{top}^B(f,\mathrm{CI}(f))\leq h_{top}^B(f,I_{\varphi}(f))\leq\delta h_{top}(f),
	\end{equation*}
	which implies that $(X,f)$ satisfies (3).
	This completes the proof of Theorem \ref{thm-Bowen-upper-2}.
\end{proof}

\begin{remark}
	Similar to the proof of \cite[Theorem 1.1]{Pavlov-2016}, it can be shown that the example constructed in Theorem \ref{thm-Bowen-upper-2} has exactly one ergodic measure of maximal entropy.
	Using arguments analogous to Theorem \ref{thm-Bowen-upper-1}, for any $m\in\mathbb{N}$, we can construct subshifts satisfying Theorem \ref{thm-Bowen-upper-2} (1)-(3) and having exactly $m$ ergodic measures of maximal entropy.
	Indeed, let $Z$ be a union of $m$ copies of $\mathcal{A}^{\mathbb{N}_0}$, where $\mathcal{A}=\{1,\cdots,A\}$ is a finite alphabet.
	Then one can show that the subshift $(X,f)$ constructed in Theorem \ref{thm-construction} associated to alphabet $\{1,\cdots,mA\}$ and subshift $Z$ has $m$ ergodic measures of maximal entropy and satisfies conditions (1)-(3) in Theorem \ref{thm-Bowen-upper-2} whenever $A$ sufficiently large.
\end{remark}

\section{Packing topological entropy}\label{Sect-packing}

In this section, we will show Theorem \ref{thm-packing}.
The following lemma is used to estimate the lower bound for packing topological entropy.

\begin{lemma}\label{packing-lemma1}
	Let $\{a_i\}_{i=1}^{\infty}$ and $\{b_i\}_{i=1}^{\infty}$ be two increasing sequences of non-negative integers satisfying $a_i\leq b_i<a_{i+1}$ for all $i\in\mathbb{N}$.
	Suppose that
	\begin{equation}\label{eq-packing-lemma1-1}
		\lim_{j\rightarrow\infty}\frac{b_{2j}-a_{2j}+1}{b_{2j}+1}=1.
	\end{equation}
	Given $\beta>0$, let $\{\Omega_i\}_{i=1}^{\infty}$ be the index sets such that for any $i\in\mathbb{N}$,
	\begin{enumerate}[label=(\alph*)]
		\item $0<|\Omega_i|<\infty$;
		\item $|\Omega_{2i}|\geq\exp(\beta\cdot(b_{2i}-a_{2i}+1))$.
	\end{enumerate}
	Let $\varepsilon>0$.
	Suppose that
	\begin{equation*}
		Y=\bigg\{x_{\xi}\in X:\xi=(\xi_1,\xi_2,\cdots)\in\prod_{i=1}^{\infty}\Omega_i\bigg\}\subset X
	\end{equation*}
	satisfies the following condition: $f^{a_{2l}}x_{\xi}$ and $f^{a_{2l}}x_{\zeta}$ are $(b_{2l}-a_{2l}+1,\varepsilon)$-separated whenever $\xi_{2l}\neq\zeta_{2l}$ for some $l\in\mathbb{N}$.
	Then
	\begin{equation}\label{eq-packing-lemma1-2}
		h_{top}^P(f,\overline{Y})\geq\beta.
	\end{equation}
\end{lemma}

\begin{proof}
	Fix $\zeta=(\zeta_1,\zeta_2,\cdots)\in\prod_{i=1}^{\infty}\Omega_i$.
	For any $l\in\mathbb{N}$, let
	\begin{equation*}
		\Delta_l=\bigg\{\xi\in\prod_{i=1}^{\infty}\Omega_i:\xi_j=\zeta_j,\forall j\geq l+1\bigg\}
	\end{equation*}
	and
	\begin{equation*}
		\mu_l=\frac{1}{|\Delta_l|}\sum_{\xi\in\Delta_l}\delta_{x_{\xi}}=\frac{1}{\prod_{j=1}^l|\Omega_j|}\sum_{\xi\in\Delta_l}\delta_{x_{\xi}}\in\mathcal{M}(X).
	\end{equation*}
	Then $\mu_l(\overline{Y})=1$.
	Let $\mu\in\mathcal{M}(X)$ be a limit point of the sequence $\{\mu_l\}_{l=1}^{\infty}$ and $\{l_k\}_{k=1}^{\infty}$ be an increasing subsequence such that $\mu_{l_k}\rightarrow\mu$ weakly as $k\rightarrow\infty$.
	Then
	\begin{equation*}
		\mu(\overline{Y})\geq\limsup_{k\rightarrow\infty}\mu_{l_k}(\overline{Y})=1.
	\end{equation*}
	From Lemma \ref{pre-entropy-lemma2}, it follows that
	\begin{equation}\label{eq-packing-lemma1-proof-1}
		h_{top}^P(f,\overline{Y})\geq\overline{h}_{\mu}(f)
		=\int_{\overline{Y}}\overline{h}_{\mu}(f,x)\mathrm{d}\mu(x)
		\geq\int_{\overline{Y}}\limsup_{n\rightarrow\infty}-\frac1n\log\mu(B_n(x,\varepsilon/2))\mathrm{d}\mu(x).
	\end{equation}
	Fix $x\in\overline{Y}$ and $j\in\mathbb{N}$, then there exists $y\in Y$ such that $d_{b_{2j}+1}(x,y)<\varepsilon/2$, and hence $B_{b_{2j}+1}(x,\varepsilon/2)\subset B_{b_{2j}+1}(y,\varepsilon)$.
	Fix any $l_k\geq 2j+1$, then we have
	\begin{equation}\label{eq-packing-lemma1-proof-2}
		\mu_{l_k}(B_{b_{2j}+1}(x,\varepsilon/2))\leq\mu_{l_k}(B_{b_{2j}+1}(y,\varepsilon))=\mu_{l_k}(B_{b_{2j}+1}(y,\varepsilon)\cap Y).
	\end{equation}
	Write $y=x_{\theta}$, where $\theta\in\prod_{i=1}^{\infty}\Omega_i$.
	For $x_{\xi}\in Y$ and $j\in\mathbb{N}$, if $\xi_{2j}\neq\theta_{2j}$, then by the assumptions, $f^{a_{2j}}(x_{\xi})$ and $f^{a_{2j}}(y)$ are $(b_{2j}-a_{2j}+1,\varepsilon)$-separated, and hence $x_{\xi}\notin B_{b_{2j}+1}(y,\varepsilon)$.
	From condition (a) and (b), we obtain that
	\begin{equation}\label{eq-packing-lemma1-proof-3}
		\mu_{l_k}(B_{b_{2j}+1}(y,\varepsilon)\cap Y)\leq\frac{\bigg(\prod_{i=1}^{2j-1}|\Omega_i|\bigg)\bigg(\prod_{i=2j+1}^{l_k}|\Omega_i|\bigg)}{\prod_{i=1}^{l_k}|\Omega_i|}=\frac{1}{|\Omega_{2j}|}\leq\exp(-\beta(b_{2j}-a_{2j}+1)).
	\end{equation}
	Combining \eqref{eq-packing-lemma1-proof-2} and \eqref{eq-packing-lemma1-proof-3}, we obtain that
	\begin{equation*}
		-\frac{1}{b_{2j}+1}\log\mu(B_{b_{2j}+1}(x,\varepsilon/2))\geq-\frac{1}{b_{2j}+1}\log\bigg(\liminf_{k\rightarrow\infty}\mu_{l_k}(B_{b_{2j}+1}(x,\varepsilon/2))\bigg)\geq\frac{b_{2j}-a_{2j}+1}{b_{2j}+1}\beta.
	\end{equation*}
	Using \eqref{eq-packing-lemma1-1}, we obtain that
	\begin{equation}\label{eq-packing-lemma1-proof-4}
		\begin{split}
			\limsup_{n\rightarrow\infty}-\frac1n\log\mu(B_n(x,\varepsilon/2))
			&\geq\limsup_{j\rightarrow\infty}-\frac{1}{b_{2j}+1}\log\mu(B_{b_{2j}+1}(x,\varepsilon/2))\\
			&\geq\limsup_{j\rightarrow\infty}\frac{b_{2j}-a_{2j}+1}{b_{2j}+1}\beta=\beta.
		\end{split}
	\end{equation}
	Combining \eqref{eq-packing-lemma1-proof-1} and \eqref{eq-packing-lemma1-proof-4}, we conclude that
		\begin{equation}
		h_{top}^P(f,\overline{Y})\geq\beta.
	\end{equation}
\end{proof}

According to Proposition \ref{pre-saturated-prop5}, if $(X,f)$ satisfies the non-uniform specification property, then every non-empty over-saturated set contains $G_{\mathcal{K}}$, where $\mathcal{K}$ is the non-empty connected closed subset of $\mathcal{M}_f(X)$ defined by \eqref{eq-pre-saturated-prop5}.
From Proposition \ref{pre-entropy-prop1}, if $G^K\neq\varnothing$, then $h_{top}^P(f,G_{\mathcal{K}})\leq h_{top}^P(f,G^K)$.
Moreover, by Lemma \ref{pre-saturated-lemma6}, we have $G_{\mathcal{K}}\subset\mathrm{CI}(f)$.
Therefore, Theorem \ref{thm-packing} is a consequence of the following theorem.

\begin{theorem}\label{packing-thm2}
	Suppose that $(X,f)$ satisfies the non-uniform specification property, then
	\begin{equation*}
		h_{top}^P(f,G_{\mathcal{K}})=h_{top}(f).
	\end{equation*}
\end{theorem}

\begin{proof}
	Let $\{\nu_j\}_{j=1}^{\infty}$ be a sequence in $\mathcal{K}$ such that
	\begin{equation}\label{eq-packing-thm2-proof-1}
		\overline{\{\nu_j:j\geq n\}}=\mathcal{K},\quad\forall n\in\mathbb{N}.
	\end{equation}
	Fix $z\in G_{\mathcal{K}}$.
	For each $j\geq1$, take an increasing sequence $\{n_j(l)\}_{l=1}^{\infty}$ such that
	\begin{equation}\label{eq-packing-thm2-proof-2}
		\lim_{l\rightarrow\infty}\delta_z^{n_j(l)}=\nu_j.
	\end{equation}
	Arbitrarily given $0<\eta<h_{top}(f)$.
	By the definition of topological entropy, we can choose $\varepsilon>0$ such that
	\begin{equation}\label{eq-packing-thm2-proof-3}
		\limsup_{n\rightarrow\infty}\frac1n\log s_n(5\varepsilon)>\eta.
	\end{equation}

	Let $\varepsilon_i=\varepsilon/2^i$ for $i\in\mathbb{Z}$.
	Now we define:
	\begin{itemize}
		\item sequences of non-negative integers $\{a_j\}_{j=1}^{\infty}$, $\{b_j\}_{j=1}^{\infty}$, $\{m_j\}_{j=1}^{\infty}$;
		\item an array of non-negative integers $\{l_{jr}:j\in\mathbb{N},1\leq r\leq j\}$.
	\end{itemize}
	Let $a_1=0$.
	Choose a positive integer $l_{11}$ such that
	\begin{equation*}
		\rho(\delta_z^{n_1(l_{11})},\nu_1)\leq\varepsilon_1.
	\end{equation*}
	Let $b_1=n_{1}(l_{11})-1$ and $a_2=b_1+M(b_1+1,\varepsilon_1)$.
	Take $m_1>a_2^2$ such that
	\begin{equation*}
		s_{m_1}(5\varepsilon)>\exp(\eta\cdot m_1)
	\end{equation*}
	and let $b_2=a_2+m_1-1$.
	Assume that $\{a_j\}_{j=1}^{2k}$, $\{b_j\}_{j=1}^{2k}$, $\{m_j\}_{j=1}^k$ and $\{l_{jr}:1\leq j\leq k,1\leq r\leq j\}$ have been defined.
	Let $a_{2k+1}=b_{2k}+M(b_{2k}+1,\varepsilon_{k+1})$.
	Take positive integers $l_{k+1,1}<l_{k+1,2}<\cdots<l_{k+1,k+1}$ such that
	\begin{equation*}
		a_{2k+1}^2<n_1(l_{k+1,1})<n_2(l_{k+1,2})<\cdots<n_{k+1}(l_{k+1,k+1})
	\end{equation*}
	and
	\begin{equation*}
		\rho(\delta_z^{n_j(l_{k+1,j})},\nu_j)\leq\varepsilon_{k+1},\quad\text{for}\;j=1,2,\cdots,k+1.
	\end{equation*}
	Let
	\begin{equation*}
		\begin{split}
			b_{2k+1}&=a_{2k+1}+n_{k+1}(l_{k+1,k+1})-1,\\
			a_{2k+2}&=b_{2k+1}+M(b_{2k+1}-a_{2k+1}+1,\varepsilon_{k+1}).
		\end{split}
	\end{equation*}
	Take $m_{k+1}>a_{2k+2}^2$ such that
	\begin{equation*}
		s_{m_{k+1}}(5\varepsilon)>\exp(\eta\cdot m_{k+1})
	\end{equation*}
	and let $b_{2k+2}=a_{2k+2}+m_{k+1}-1$.
	Inductively, we can define $\{a_j\}_{j=1}^{\infty}$, $\{b_j\}_{j=1}^{\infty}$, $\{m_j\}_{j=1}^{\infty}$ and $\{l_{jr}:j\in\mathbb{N},1\leq r\leq j\}$.
	It is easy to see that for any $j\in\mathbb{N}$, we have $a_j<b_j<a_{j+1}$ and $m_j<m_{j+1}$.
	Moreover,
	\begin{equation*}
		\lim_{j\to\infty}\frac{b_{2j}-a_{2j}+1}{b_{2j}+1}=1.
	\end{equation*}

	For each $k\in\mathbb{N}$, let $\Gamma_k\subset X$ be an $(m_k,5\varepsilon)$-separated set with $|\Gamma_k|=s_{m_k}(5\varepsilon)$.
	Now we construct the index sets $\{\Omega_i\}_{i=1}^{\infty}$ satisfying the condition (a) and (b) in Lemma \ref{packing-lemma1}.
	For $i\in\mathbb{N}$, define
	\begin{equation*}
		\Omega_i=
		\begin{cases}
			\{z\} & \text{if}\quad i=2k-1,\\
			\Gamma_k & \text{if}\quad i=2k.
		\end{cases}
	\end{equation*}
	It follows from \eqref{eq-packing-thm2-proof-3} that
	\begin{itemize}
		\item $0<|\Omega_i|<\infty$;
		\item $|\Omega_{2i}|=|\Gamma_i|=s_{m_i}(5\varepsilon)>\exp(\eta\cdot m_i)=\exp(\eta(b_{2i}-a_{2i}+1))$.
	\end{itemize}
	Therefore $\{\Omega_i\}_{i=1}^{\infty}$ satisfies the condition (a) and (b) in Lemma \ref{packing-lemma1} by taking $\beta=\eta$.
	Let
	\begin{equation*}
		\Xi_l=\prod_{j=1}^l\Omega_j\quad\text{and}\quad\Xi=\prod_{j=1}^{\infty}\Omega_j.
	\end{equation*}
	
	For each $\xi=(x_1,x_2,\cdots)\in\Xi$, we define non-empty closed sets $E_k^{\xi}$ inductively for every $k\in\mathbb{N}$ starting with
	\begin{equation*}
		E_1^{\xi}=\Bigg(\bigcap_{i=a_1}^{b_1}f^{-i}\left(\overline{B(f^{i-a_1}x_1,\varepsilon_1)}\right)\Bigg)
		\cap\Bigg(\bigcap_{i=a_2}^{b_2}f^{-i}\left(\overline{B(f^{i-a_2}x_2,\varepsilon_1)}\right)\Bigg).
	\end{equation*}
	Clearly $E_1^{\xi}$ is closed.
	From the non-uniform specification property, it follows that $E_1^{\xi}$ is a non-empty and closed set.
	Assume that $E_k^{\xi}$ has been defined.
	Define $\widehat{E}_k^{\xi}=\bigcup_{y\in E_k^{\xi}}\overline{B_{b_{2k}+1}(y,\varepsilon_{k+1})}$.
	Note that the set $\widehat{E}_k^{\xi}$ is the projection onto the second coordinate of the closed subset
	\begin{equation*}
		\{(y,y')\in E_k^{\xi}\times X:d_{b_{2k}+1}(y,y')\leq\varepsilon_{k+1}\}
	\end{equation*}
	of the compact space $E_k^{\xi}\times X$.
	Hence $\widehat{E}_k^{\xi}$ is compact.
	Define
	\begin{equation*}
		E_{k+1}^{\xi}=\widehat{E}_k^{\xi}
		\cap\bigg(\bigcap_{i=a_{2k+1}}^{b_{2k+1}}f^{-i}\Big(\overline{B(f^{i-a_{2k+1}}x_{2k+1},\varepsilon_{k+1})}\Big)\bigg)
		\cap\bigg(\bigcap_{i=a_{2k+2}}^{b_{2k+2}}f^{-i}\Big(\overline{B(f^{i-a_{2k+2}}x_{2k+2},\varepsilon_{k+1})}\Big)\bigg),
	\end{equation*}
	From the non-uniform specification property, it follows that $E_{k+1}^{\xi}$ is a non-empty closed set.
	Similarly, for every $l\in\mathbb{N}$, $1\leq k\leq l$ and any $\zeta\in\Xi_{2l}$, we can define a non-empty closed set $E_k^{\zeta}$.
	Note that if $\xi=(x_1,x_2\cdots)\in\Xi$ satisfies the first $2l$ coefficients coincide with $\zeta$, that is, $\zeta=(x_1,\cdots,x_{2l})$, then for any $1\leq k\leq l$, we have $E_k^{\zeta}=E_k^{\xi}$.

	For $k\in\mathbb{N}$, define
	\begin{equation*}
		F_k=\bigcup_{\zeta=(z_1,\cdots,z_{2k})\in\Xi_{2k}}\bigcap_{j=1}^k
		\left(\bigcap_{i=a_{2j-1}}^{b_{2j-1}}f^{-i}\left(\overline{B(f^{i-a_{2j-1}}z_{2j-1},\varepsilon_{j-2})}\right)
		\cap\bigcap_{i=a_{2j}}^{b_{2j}}f^{-i}\left(\overline{B(f^{i-a_{2j}}z_{2j},\varepsilon_{j-2})}\right)\right)
	\end{equation*}
	and $F=\bigcap_{k\geq 1}F_k$.
	Clearly, every $F_k$ is a non-empty compact set, satisfying $F_k\supset F_{k+1}$.
	Thus $F$ is a non-empty compact set.
	We proceed to show that $F\subset G_{\mathcal{K}}$.
	By Lemma \ref{pre-saturated-lemma1}, Proposition \ref{pre-saturated-prop5}, and \eqref{eq-packing-thm2-proof-1}, it suffices to show that $x\in G^{\nu_j}$ for any $x\in F$ and any $j\in\mathbb{N}$.
	Fix $x\in F$ and $j\in\mathbb{N}$.
	For any $k>j$, since $x\in F_k$, there is a $\zeta=(z_1,\cdots,z_{2k})\in\Xi_{2k}$ such that for any $a_{2k-1}\leq i\leq b_{2k-1}$, we have
	\begin{equation*}
		d(f^ix,f^{i-a_{2k-1}}z)=d(f^ix,f^{i-a_{2k-1}}z_{2k-1})\leq\varepsilon_{k-2}.
	\end{equation*}
	Since $b_{2k-1}-a_{2k-1}+1=n_k(l_{kk})\geq n_j(l_{kj})$, we have
	\begin{equation}\label{eq-packing-thm2-proof-4}
		d(f^{a_{2k-1}+r}x,f^rz)\leq\varepsilon_{k-2},\quad \forall\,0\leq r\leq n_j(l_{kj})-1.
	\end{equation}
	Since $n_j(l_{kj})>n_1(l_{k1})>a_{2k-1}^2$, combining Lemma \ref{pre-metric-lemma2}, \eqref{eq-packing-thm2-proof-2} and \eqref{eq-packing-thm2-proof-4}, we obtain that
	\begin{align*}
		\rho(\delta_x^{a_{2k-1}+n_j(l_{kj})},\nu_j)
		&\leq\rho(\delta_x^{a_{2k-1}+n_j(l_{kj})},\delta_z^{n_j(l_{kj})})+\rho(\delta_z^{n_j(l_{kj})},\nu_j)\\
		&\leq\varepsilon_{k-2}+\frac{a_{2k-1}}{a_{2k-1}+n_j(l_{kj})}\mathrm{diam}(X)+\varepsilon_k
		\leq\varepsilon_{k-3}+\frac{1}{a_{2k-1}}\mathrm{diam}(X).
	\end{align*}
	Therefore $\delta_x^{a_{2k-1}+n_j(l_{kj})}\rightarrow\nu_j$ as $k\rightarrow\infty$, which implies that $\nu_j\in V_f(x)$.
	As a result, $F\subset G_{\mathcal{K}}$.
	
	Now we construct a subset $Y\subset G_{\mathcal{K}}$ to apply Lemma \ref{packing-lemma1}.
	For every $\xi=(x_1,x_2,\cdots)\in\Xi$ and $k\in\mathbb{N}$, define
	\begin{equation*}
		Y_k^{\xi}=\bigcap_{j=1}^k\Bigg(\bigcap_{i=a_{2j-1}}^{b_{2j-1}}f^{-i}\Big(\overline{B(f^{i-a_{2j-1}}x_{2j-1},\varepsilon_{j-2})}\Big)
		\cap\bigcap_{i=a_{2j}}^{b_{2j}}f^{-i}\Big(\overline{B(f^{i-a_{2j}}x_{2j},\varepsilon_{j-2})}\Big)\Bigg).
	\end{equation*}
	Clearly, every $Y_k^{\xi}$ is a non-empty closed set satisfying $E_k^{\xi}\subset Y_k^{\xi}$.
	Define
	\begin{equation*}
		Y^{\xi}=\bigcap_{k\geq 1}Y_k^{\xi}=\bigcap_{j\geq1}\Bigg(\bigcap_{i=a_{2j-1}}^{b_{2j-1}}f^{-i}\Big(\overline{B(f^{i-a_{2j-1}}x_{2j-1},\varepsilon_{j-2})}\Big)
		\cap\bigcap_{i=a_{2j}}^{b_{2j}}f^{-i}\Big(\overline{B(f^{i-a_{2j}}x_{2j},\varepsilon_{j-2})}\Big)\Bigg),
	\end{equation*}
	then $Y^{\xi}$ is a non-empty closed set.
	It is clear that for any $\xi\in\Xi$, we have $Y_k^{\xi}\subset F_k$ for all $k\in\mathbb{N}$ and $Y^{\xi}\subset F$.
	For each $\xi\in\Xi$, we fix an $x_{\xi}\in Y^{\xi}$.
	Define
	\begin{equation*}
		Y=\left\{x_{\xi}:\xi=(x_1,x_2,\cdots)\in\prod_{i=1}^{\infty}\Omega_i\right\}.
	\end{equation*}
	Then $Y\subset\bigcup_{\xi\in\Xi}Y^{\xi}\subset F$.
	To apply Lemma \ref{packing-lemma1}, it remains to verify that, for any $\xi,\zeta\in\Xi$, where $\xi=(x_1,x_2,\cdots)$ and $\zeta=(z_1,z_2,\cdots)$, if there exists $l\in\mathbb{N}$ with $x_{2l}\neq z_{2l}$, then $f^{a_{2l}}x_{\xi}$ and $f^{a_{2l}}x_{\zeta}$ are $(b_{2l}-a_{2l}+1,\varepsilon)$-separated.
	Since $x_{\xi}\in Y_l^{\xi}$, by the construction of $Y^{\xi}$, for any $0\leq i\leq b_{2l}-a_{2l}$, we have
	\begin{equation}\label{eq-packing-thm2-proof-5}
		\begin{split}
			d(f^{a_{2l}+i}x_{\xi},f^{a_{2l}+i}x_{\zeta})
			&\geq d(f^ix_{2l},f^iz_{2l})
			-\Big(d(f^{a_{2l}+i}x_{\xi},f^ix_{2l})+d(f^{a_{2l}+i}x_{\zeta},f^iz_{2l})\Big)\\
			&\geq d(f^ix_{2l},f^iz_{2l})-2\varepsilon_{l-2}\geq d(f^ix_{2l},f^iz_{2l})-4\varepsilon.
		\end{split}
	\end{equation}
	By the definition of $\{\Omega_i\}_{i=1}^{\infty}$, we have $x_{2l},z_{2l}\in\Gamma_l$.
	Recall that $\Gamma_l$ is an $(m_l,5\varepsilon)$-separated set with $|\Gamma_l|=s_{m_l}(5\varepsilon)$.
	Since $x_{2l}\neq z_{2l}$, we obtain that $x_{2l}$ and $z_{2l}$ are $(m_l,5\varepsilon)$-separated.
	Hence there exists an integer $0\leq j\leq m_l-1$ such that $d(f^jx_{2l},f^jz_{2l})>5\varepsilon$.
	Since $m_l=b_{2l}-a_{2l}+1$, we have $0\leq j\leq b_{2l}-a_{2l}$.
	Using \eqref{eq-packing-thm2-proof-5}, we conclude that
	\begin{equation*}
		d(f^{a_{2l}+j}x_{\xi},f^{a_{2l}+j}x_{\zeta})\geq d(f^jx_{2l},f^jz_{2l})-4\varepsilon>\varepsilon.
	\end{equation*}
	As a result, $f^{a_{2l}}x_{\xi}$ and $f^{a_{2l}}x_{\zeta}$ are $(b_{2l}-a_{2l}+1,\varepsilon)$-separated.
	Applying Lemma \ref{packing-lemma1}, we have $h_{top}^P(f,\overline{Y})\geq\eta$.
	Since $Y\subset F\subset G_{\mathcal{K}}$ and $F$ is compact, we obtain that
	\begin{equation}
		h_{top}^P(f,G_{\mathcal{K}})\geq h_{top}^P(f,F)\geq h_{top}^P(f,\overline{Y})\geq\eta.
	\end{equation}
	By the arbitrariness of $\eta\in(0,h_{top}(f))$, we conclude that
	\begin{equation*}
		h_{top}^P(f,G_{\mathcal{K}})=h_{top}(f).
	\end{equation*}
	This completes the proof of Theorem \ref{packing-thm2}.
\end{proof}

\begin{proof}[{\bf Proof of Theorem \ref{thm-packing}}]
	According to Proposition \ref{pre-saturated-prop5} and Lemma \ref{pre-saturated-lemma6}, we have $G_{\mathcal{K}}\subset G^K$ and $G_{\mathcal{K}}\subset\mathrm{CI}(f)$.
	From Theorem \ref{packing-thm2}, it follows that
	\begin{equation*}
		h_{top}(f)\geq h_{top}^P(f,G^K)\geq h_{top}^P(f,G_{\mathcal{K}})=h_{top}(f)
	\end{equation*}
	and
	\begin{equation*}
		h_{top}(f)\geq h_{top}^P(f,\mathrm{CI}(f))\geq h_{top}^P(f,G_{\mathcal{K}})=h_{top}(f).
	\end{equation*}
	This completes the proof of Theorem \ref{thm-packing}.
\end{proof}

\section{Conditional variational principle}\label{Sect-CVP}

\begin{proposition}[{\cite[Proposition 7.1]{Pfister-Sullivan-2007}}]\label{CVP-prop1}
	Let $(X,f)$ be a dynamical system.
	Suppose that
	\begin{equation*}
		h_{top}^B(f,G_{\mu})=h_{\mu}(f),\quad\forall\,\mu\in\mathcal{M}_f(X).
	\end{equation*}
	Then the conditional variational principle holds, i.e. for every $\varphi\in C(X)$ and every $a\in L_{\varphi}$,
	\begin{equation*}
		h_{top}^B(f,R_{\varphi}(a))
		=\sup\bigg\{h_{\mu}(f):\mu\in\mathcal{M}_f(X),\ \int_X\varphi\mathrm{d}\mu=a\bigg\}.
	\end{equation*}
\end{proposition}

\begin{theorem}\label{CVP-thm2}
	Suppose that $(X,f)$ satisfies the non-uniform $0$-specification property.
	Then
	\begin{equation*}
		h_{top}^B(f,G_{\mu})=h_{\mu}(f),\quad\forall\,\mu\in\mathcal{M}_f(X).
	\end{equation*}
\end{theorem}

\begin{proof}
	Let $M(n,\varepsilon)$ be an admissible gap function such that
	\begin{equation*}
		\tau(M)=\sup_{\varepsilon>0}\liminf_{n\rightarrow\infty}\frac{M(n,\varepsilon)}{n}=0.
	\end{equation*}
	Fix $\mu\in\mathcal{M}_f(X)$.
	If $h_{\mu}(f)<\infty$, then $G_{\mu}\subset\mathrm{QR}(h_\mu(f))$, where $\mathrm{QR}(t)$ is defined in Proposition \ref{Bowen-upper-prop1}.
	It follows from Proposition \ref{Bowen-upper-prop1} that $h_{top}^B(f,G_{\mu})\leq h_{\mu}(f)$.
	If $h_\mu(f)=\infty$, then the upper bound is trivial.

	It remains to show that
	\begin{equation*}
		h_{top}^B(f,G_{\mu})\geq h_{\mu}(f).
	\end{equation*}
	There is nothing to prove if $h_{\mu}(f)=0$.
	Assume that $h_{\mu}(f)>0$.
	Fix two positive numbers $h^*$ and $\gamma$ with
	\begin{equation*}
		0<h^*<\gamma<h_{\mu}(f).
	\end{equation*}
	For a neighborhood $F$ of $\nu\in\mathcal{M}_f(X)$ in $\mathcal{M}(X)$, recall that $N(F;\delta,n,\varepsilon)$ denotes the maximal cardinality of a $(\delta,n,\varepsilon)$-separated subset of $X_{n,F}$, where
	\begin{equation*}
		X_{n,F}=\bigg\{x\in X:\frac1n\sum_{i=0}^{n-1}\delta_{f^ix}\in F\bigg\}.
	\end{equation*}
	Recall that
	\begin{equation*}
		\underline{s}(\nu,\delta,\varepsilon)=\inf_{F\ni\nu}\liminf_{n\rightarrow\infty}\frac1n\log N(F;\delta,n,\varepsilon).
	\end{equation*}
	By Lemma \ref{pre-separation-lemma2}, there exist $\delta^*>0$ and $\varepsilon^*>0$ such that, for every $j\in\mathbb{N}$, one can choose a finite convex combination of ergodic measures with positive rational coefficients
	\begin{equation*}
		\mu_j=\sum_{u=1}^{m_j}p_{j,u}\mu_{j,u}
	\end{equation*}
	satisfying
	\begin{equation}\label{eq-CVP-thm2-proof-1}
		\rho(\mu_j,\mu)<\eta_j
		\quad\text{and}\quad
		\sum_{u=1}^{m_j}p_{j,u}\underline{s}(\mu_{j,u},\delta^*,\varepsilon^*)>\gamma,
	\end{equation}
	where $\eta_j\rightarrow0$.

	Choose a positive number $\varepsilon_0$ and a sequence $\{\theta_j\}_{j=1}^{\infty}$ with $0<\theta_j<1$ decreasing to $0$ such that
	\begin{equation}\label{eq-CVP-thm2-proof-2}
		\varepsilon_0<\min\bigg\{1,\frac{\varepsilon^*}{4}\bigg\}
		\quad\text{and}\quad
		\gamma-\theta_j>h^*(1+\theta_j)
	\end{equation}
	for every $j\in\mathbb{N}$.
	For any $i\in\mathbb{N}_0$, let $\varepsilon_i=2^{-i}\varepsilon_0$.
	Choose positive integer $b_j$ such that for any $1\leq u\leq m_j$, $b_jp_{j,u}$ is a positive integer.
	Let
	\begin{equation*}
		c_j=\frac12\min_{1\leq u\leq m_j}p_{j,u}.
	\end{equation*}
	Since $\tau(M)=0$, for any $j\geq1$, there exists an increasing sequence $\{\overline{N}_i\}_{i\in\mathbb{N}}$ of positive integers such that
	\begin{equation*}
		M(\overline{N}_i,\varepsilon_j)<\theta_jc_j\overline{N}_i,\forall\,i\in\mathbb{N}.
	\end{equation*}
	Let $N_j=\overline{N}_i$ for a sufficiently large $i$ and $n_j$ be the largest multiple of $b_j$ with $n_j\leq N_j$.
	We may assume that $n_j\geq N_j/2$.
	For each $1\leq u\leq m_j$, by the definition of $\underline{s}(\mu_{j,u},\delta^*,\varepsilon^*)$, we may assume that $N_j$ is sufficiently large by taking larger $i$ if necessary, such that there exists $(\delta^*,l_{j,u},\varepsilon^*)$-separated set $\Gamma_{j,u}\subset X_{l_{j,u},B(\mu_{j,u},\eta_j)}$ with
	\begin{equation}\label{eq-CVP-thm2-proof-3}
		|\Gamma_{j,u}|\geq \exp\Big(l_{j,u}\big(\underline{s}(\mu_{j,u},\delta^*,\varepsilon^*)-\theta_j\big)\Big),
	\end{equation}
	where $l_{j,u}=p_{j,u}n_j$.
	For $1\leq u\leq m_j$, let $q_{j,u}=M(l_{j,u},\varepsilon_j)$.
	Note that $c_jN_j\leq l_{j,u}\leq N_j$.
	Since $M(\,\cdot\,,\varepsilon_j)$ is non-decreasing, we obtain that
	\begin{equation}\label{eq-CVP-thm2-proof-4}
		q_{j,u}\leq M(N_j,\varepsilon_j)<\theta_jc_jN_j\leq\theta_jl_{j,u}.
	\end{equation}

	Define
	\begin{equation*}
		\Gamma_j=\prod_{u=1}^{m_j}\Gamma_{j,u},\quad L_j=n_j+\sum_{u=1}^{m_j}q_{j,u}.
	\end{equation*}
	Then $L_j<n_j(1+\theta_j)$.
	Combining \eqref{eq-CVP-thm2-proof-1}, \eqref{eq-CVP-thm2-proof-2}, \eqref{eq-CVP-thm2-proof-3} and \eqref{eq-CVP-thm2-proof-4}, we have
	\begin{equation}\label{eq-CVP-thm2-proof-5}
		\log|\Gamma_j|>n_j(\gamma-\theta_j)>h^* n_j(1+\theta_j)>h^*L_j.
	\end{equation}
	For $1\leq u\leq m_j$, let $a(j,1)=0$ and for $2\leq u\leq m_j$, let
	\begin{equation*}
		a(j,u)=\sum_{v=1}^{u-1}(l_{j,v}+q_{j,v}).
	\end{equation*}
	We next choose $r_j$ and $M_j$ inductively.
	Suppose that $r_1,\ldots,r_{j-1}$ and $M_1,\ldots,M_{j-1}$ have been chosen, and let
	\begin{equation*}
		R_{j-1}=\sum_{i=1}^{j-1}(r_iL_i+M_i)\;\text{for}\;j\geq 2,
	\end{equation*}
	where $R_0=0$.
	Since $\tau(M)=0$, there exists a strictly increasing sequence of positive integers $\{\widehat{N}_i\}_{i\in\mathbb{N}}$ such that
	\begin{equation*}
		M(\widehat{N}_i,\varepsilon_{j+1})<\frac{\theta_j}{2}\widehat{N}_i,\quad\forall\,i\in\mathbb{N}.
	\end{equation*}
	Let $\widetilde{N}_j=\widehat{N}_i$ for a sufficiently large $i$ and $r_j$ be the largest positive integer such that
	\begin{equation*}
		T_j:=R_{j-1}+r_jL_j\leq\widetilde{N}_j.
	\end{equation*}
	We may assume that $\widetilde{N}_j$ is sufficiently large by taking larger $i$ if necessary, such that $r_j\geq2$ and
	\begin{equation}\label{eq-CVP-thm2-proof-6}
		R_{j-1}+L_{j+1}<\theta_jr_jL_j\quad\text{and}\quad T_j\geq\widetilde{N}_j/2.
	\end{equation}
	Let $M_j=M(T_j,\varepsilon_{j+1})$.
	Since $T_j\leq\widetilde{N}_j$ and $M(\,\cdot\,,\varepsilon_{j+1})$ is non-decreasing, we have
	\begin{equation}\label{eq-CVP-thm2-proof-7}
		M_j\leq M(\widetilde{N}_j,\varepsilon_{j+1})<\frac{\theta_j}{2}\widetilde{N}_j\leq\theta_jT_j.
	\end{equation}

	For $j\in\mathbb{N}$ and $0\leq k\leq r_j-1$, let
	\begin{equation*}
		s(j,k)=\sum_{i=1}^{j-1}(r_iL_i+M_i)+kL_j
	\end{equation*}
	and, for $1\leq u\leq m_j$, let
	\begin{equation*}
		a(j,k,u)=s(j,k)+a(j,u).
	\end{equation*}
	Define
	\begin{equation*}
		\Xi_j=\prod_{i=1}^j\Gamma_i^{r_i}
		\quad\text{and}\quad
		\Xi=\prod_{i=1}^{\infty}\Gamma_i^{r_i}.
	\end{equation*}
	For every $\xi=(\xi_1,\xi_2,\ldots)\in\Xi$, write $\xi_j=(\xi_j^0,\ldots,\xi_j^{r_j-1})\in\Gamma_j^{r_j}$ and $\xi_j^k=(x_{j,k,1},\ldots,x_{j,k,m_j})\in\Gamma_j$ for every $j\in\mathbb{N}$ and $0\leq k\leq r_j-1$.
	Define $E_j^{\xi}$ inductively for each $j\in\mathbb{N}$, starting with
	\begin{equation}\label{eq-CVP-thm2-proof-8}
		E_1^{\xi}=\bigcap_{k=0}^{r_1-1}\bigcap_{u=1}^{m_1}\bigcap_{r=0}^{l_{1,u}-1}f^{-(a(1,k,u)+r)}\Big(\overline{B(f^rx_{1,k,u},\varepsilon_1)}\Big).
	\end{equation}
	From the non-uniform specification property, it follows that $E_1^{\xi}$ is a non-empty and closed set.
	Assume that $E_i^{\xi}$ has been defined for every $1\leq i\leq j-1$.
	Define
	\begin{equation*}
		\widehat{E}_{j-1}^{\xi}=\bigcup_{x\in E_{j-1}^{\xi}}\overline{B_{T_{j-1}}(x,\varepsilon_j)}.
	\end{equation*}
	Note that the set $\widehat{E}_{j-1}^{\xi}$ is the projection onto the second coordinate of the closed subset
	\begin{equation*}
		\{(x,x')\in E_{j-1}^{\xi}\times X:d_{T_{j-1}}(x,x')\leq\varepsilon_j\}
	\end{equation*}
	of the compact space $E_{j-1}^{\xi}\times X$.
	Hence $\widehat{E}_{j-1}^{\xi}$ is compact.
	Define
	\begin{equation}\label{eq-CVP-thm2-proof-9}
		E_j^{\xi}=\widehat{E}_{j-1}^{\xi}\cap\bigg(\bigcap_{k=0}^{r_j-1}\bigcap_{u=1}^{m_j}\bigcap_{r=0}^{l_{j,u}-1}f^{-(a(j,k,u)+r)}\Big(\overline{B(f^rx_{j,k,u},\varepsilon_j)}\Big)\bigg).
	\end{equation}
	By the definition of $M_{j-1}$ and the non-uniform specification property, $E_j^{\xi}$ is a non-empty and closed set.
	Inductively, we can define non-empty closed sets $E_j^{\xi}$ for every $j\in\mathbb{N}$.
	Similarly, for any $l\in\mathbb{N}$ and $\zeta\in\Xi_l$, we can define $E_j^{\zeta}$ for every $1\leq j\leq l$.

	Given $\zeta=(\zeta_1,\ldots,\zeta_l)\in\Xi_l$, write $\zeta_j=(\zeta_j^0,\ldots,\zeta_j^{r_j-1})\in\Gamma_j^{r_j}$ and $\zeta_j^k=(z_{j,k,1},\ldots,z_{j,k,m_j})\in\Gamma_j$ for every $1\leq j\leq l$ and $0\leq k\leq r_j-1$.
	Define
	\begin{equation}\label{eq-CVP-thm2-proof-10}
		F_l^{\zeta}=\bigcap_{j=1}^l\bigcap_{k=0}^{r_j-1}\bigcap_{u=1}^{m_j}\bigcap_{r=0}^{l_{j,u}-1}f^{-(a(j,k,u)+r)}\Big(\overline{B(f^rz_{j,k,u},\varepsilon_{j-1})}\Big),
	\end{equation}
	and let
	\begin{equation*}
		F_l=\bigcup_{\zeta\in\Xi_l}F_l^{\zeta}
		\quad\text{and}\quad
		F=\bigcap_{l=1}^{\infty}F_l.
	\end{equation*}
	It follows from \eqref{eq-CVP-thm2-proof-8}, \eqref{eq-CVP-thm2-proof-9} and \eqref{eq-CVP-thm2-proof-10} that $E_l^{\zeta}\subset F_l^{\zeta}$ for every $l\in\mathbb{N}$ and $\zeta\in\Xi_l$.
	Then the sets $F_l$ are non-empty compact sets satisfying $F_{l+1}\subset F_l$, and hence $F$ is non-empty and compact.
	Similarly, for every $\xi=(\xi_1,\xi_2,\ldots)\in\Xi$ and $l\in\mathbb{N}$, we can define the non-empty compact set $F_l^{\xi}$ as
	\begin{equation*}
		F_l^{\xi}=\bigcap_{j=1}^l\bigcap_{k=0}^{r_j-1}\bigcap_{u=1}^{m_j}\bigcap_{r=0}^{l_{j,u}-1}f^{-(a(j,k,u)+r)}\Big(\overline{B(f^rx_{j,k,u},\varepsilon_{j-1})}\Big),
	\end{equation*}
	satisfying $E_l^{\xi}\subset F_l^{\xi}$ and $F_l^{\xi}\supset F_{l+1}^{\xi}$.
	Consequently, the intersection
	\begin{equation*}
		F^{\xi}:=\bigcap_{l=1}^{\infty}F_l^{\xi}
	\end{equation*}
	is a non-empty compact subset of $F$.
	Moreover, by the construction of $F$ and $F^{\xi}$, it is easy to see that
	\begin{equation*}
		F=\bigcup_{\xi\in\Xi}F^{\xi}.
	\end{equation*}

	We claim that $F\subset G_{\mu}$.
	For any $\xi\in\Xi$ and $x\in F^{\xi}$, by the definition of $F_l^{\xi}$, for every $j\in\mathbb{N}$, $0\leq k\leq r_j-1$, $1\leq u\leq m_j$ and $0\leq r\leq l_{j,u}-1$, we have
	\begin{equation}\label{eq-CVP-thm2-proof-11}
		d(f^{a(j,k,u)+r}x,f^rx_{j,k,u})\leq\varepsilon_{j-1}.
	\end{equation}
	By \eqref{eq-CVP-thm2-proof-1}, \eqref{eq-CVP-thm2-proof-4}, \eqref{eq-CVP-thm2-proof-11} and Lemma \ref{pre-metric-lemma2}, for any $x\in F$, $j\in\mathbb{N}$ and $0\leq k\leq r_j-1$, we have
	\begin{equation}\label{eq-CVP-thm2-proof-12}
		\rho\big(\delta_{f^{s(j,k)}x}^{L_j},\mu\big)\leq\varepsilon_{j-1}+2\eta_j+\theta_j\mathrm{diam}(X).
	\end{equation}
	By Proposition \ref{pre-metric-prop1} and Lemma \ref{pre-metric-lemma2}, we have
	\begin{equation*}
		\rho(\delta_x^{T_j},\mu)\leq\frac{1}{T_j}\sum_{i=1}^jr_iL_i\bigg(\varepsilon_{i-1}+2\eta_i+\theta_i\mathrm{diam}(X)\bigg)+\frac{\sum_{i=1}^{j-1}M_i}{T_j}\mathrm{diam}(X).
	\end{equation*}
	By \eqref{eq-CVP-thm2-proof-6} and \eqref{eq-CVP-thm2-proof-7}, the right-hand side tends to $0$ as $j\rightarrow\infty$.
	Hence $\rho(\delta_x^{T_j},\mu)\rightarrow0$ as $j\to\infty$.
	Since $R_j=T_j+M_j$, by \eqref{eq-CVP-thm2-proof-7} and Proposition \ref{pre-metric-prop1}, we have
	\begin{equation*}
		\rho(\delta_x^{R_j},\mu)\leq\rho(\delta_x^{T_j},\mu)+\frac{M_j}{T_j}\mathrm{diam}(X)\rightarrow0.
	\end{equation*}
	For any $0\leq k\leq r_{j+1}-1$, using \eqref{eq-CVP-thm2-proof-12} and Proposition \ref{pre-metric-prop1}, we have
	\begin{equation*}
		\rho(\delta_x^{s(j+1,k)},\mu)\leq\rho(\delta_x^{R_j},\mu)+\varepsilon_j+2\eta_{j+1}+\theta_{j+1}\mathrm{diam}(X)\rightarrow0.
	\end{equation*}
	By \eqref{eq-CVP-thm2-proof-6}, we have $L_{j+1}/R_j\rightarrow0$ as $j\to\infty$.
	For every sufficiently large $n$, there is a $j\in\mathbb{N}$ such that either $T_j\leq n<R_j$, or there is a $0\leq k\leq r_{j+1}-1$ with
	\begin{equation*}
		s(j+1,k)\leq n<s(j+1,k)+L_{j+1}.
	\end{equation*}
	In the second alternative,
	\begin{equation*}
		\rho(\delta_x^n,\mu)\leq\rho(\delta_x^{s(j+1,k)},\mu)+\frac{L_{j+1}}{R_j}\mathrm{diam}(X).
	\end{equation*}
	In the first alternative,
	\begin{equation*}
		\rho(\delta_x^n,\mu)\leq\rho(\delta_x^{T_j},\mu)+\frac{M_j}{T_j}\mathrm{diam}(X).
	\end{equation*}
	Therefore $\rho(\delta_x^n,\mu)\rightarrow0$ as $n\to\infty$ for every $x\in F$.
	Thus $F\subset G_{\mu}$.

	It remains to estimate the Bowen topological entropy of $F$.
	Arbitrarily take $x_{\xi}\in F^{\xi}$ for each $\xi\in\Xi$.
	Let
	\begin{equation*}
		Y=\{x_{\xi}:\xi\in\Xi\}.
	\end{equation*}
	For $j\in\mathbb{N}$ and $0\leq k\leq r_j-1$, let
	\begin{equation*}
		\Omega_j(k)=\Gamma_j.
	\end{equation*}
	For any $\xi,\zeta\in\Xi$, similar to the proof of Theorem \ref{Bowen-thm2}, if $\xi_j^k\neq\zeta_j^k$, then $f^{s(j,k)}x_{\xi}$ and $f^{s(j,k)}x_{\zeta}$ are $(L_j,\varepsilon^*/2)$-separated.
	Furthermore, by \eqref{eq-CVP-thm2-proof-6} and \eqref{eq-CVP-thm2-proof-7}, we have $R_{j-1}/r_jL_j\to0$, $M_j/r_jL_j\to0$ and $L_{j+1}/r_jL_j\to 0$ as $j\to\infty$.
	Consequently,
	\begin{equation*}
		\liminf_{j\rightarrow\infty}\frac{\sum_{i=1}^jr_iL_i}{s(j+1,1)-1}=\liminf_{j\rightarrow\infty}\frac{\sum_{i=1}^jr_iL_i}{\sum_{i=1}^j(r_iL_i+M_i)+L_{j+1}-1}=1.
	\end{equation*}
	By \eqref{eq-CVP-thm2-proof-5}, for every $j\in\mathbb{N}$ and $0\leq k\leq r_j-1$, we have
	\begin{equation*}
		|\Omega_j(k)|=|\Gamma_j|\geq\exp(h^*L_j).
	\end{equation*}
	Applying Lemma \ref{Bowen-lemma1} with $n_j=r_j$, $t_j=L_j$, $\alpha=1$ and $\beta=h^*$, we conclude that
	\begin{equation*}
		h_{top}^B(f,F)\geq h_{top}^B(f,\overline{Y})\geq h^*.
	\end{equation*}
	Since $F\subset G_{\mu}$, it follows that $h_{top}^B(f,G_{\mu})\geq h^*$.
	By the arbitrariness of $0<h^*<\gamma<h_{\mu}(f)$, we obtain that
	\begin{equation*}
		h_{top}^B(f,G_{\mu})\geq h_{\mu}(f).
	\end{equation*}
	This completes the proof.
\end{proof}

\begin{corollary}\label{CVP-cor3}
	Suppose that $(X,f)$ satisfies the non-uniform $0$-specification property.
	Then the conditional variational principle holds.
\end{corollary}

\begin{proof}[{\bf Proof of Theorem \ref{thm-CVP}}]
	First, we consider the case $0<\tau<\infty$.
	Let $(X,f)$ be the subshift constructed in Theorem \ref{thm-Bowen-upper-1}.
	Recall that $Z=Z_1\cup\cdots\cup Z_m$ is a union of $m$ copies of a strictly ergodic subshift.
	These components are chosen with positive entropy, so $h_{\nu_k}(f)=h_{top}(f)>0$ for every $1\leq k\leq m$.
	Consider $\psi=1-\chi_{[0]}\in C(X)$.
	Clearly, $L_{\psi}=[0,1]$ and $\mathrm{Int}(L_{\psi})=(0,1)$.
	Note that for any $x=x_0x_1\cdots\in X\setminus(G_{\nu_1}\cup\cdots\cup G_{\nu_m})=X\setminus\bigcup_{i\geq0}f^{-i}Z$, the upper density of $0$ in $x_0,x_1,\cdots$ is at least $\tau/(1+\tau)$.
	It follows that
	\begin{equation}\label{eq-thm-CVP-proof-1}
		\begin{aligned}
			\liminf_{n\rightarrow\infty}\frac{1}{n}\sum_{i=0}^{n-1}\psi(f^ix)&\leq\frac{1}{1+\tau} && \text{if }x\in X\setminus(G_{\nu_1}\cup\cdots\cup G_{\nu_m}),\\
			\lim_{n\rightarrow\infty}\frac{1}{n}\sum_{i=0}^{n-1}\psi(f^ix)&=1 && \text{if }x\in G_{\nu_1}\cup\cdots\cup G_{\nu_m}.
		\end{aligned}
	\end{equation}
	Let $\gamma=1/(1+\tau)\in(0,1)$.
	Then $R_{\psi}(a)=\varnothing$ for every $\gamma<a<1$.
	Using \eqref{eq-thm-CVP-proof-1}, we obtain that
	\begin{equation*}
		\int_X\psi\mathrm{d}\nu_k=1\;\text{for}\;1\leq k\leq m
	\end{equation*}
	and for any $\mu\in\mathcal{M}_f^{erg}(X)\setminus\{\nu_1,\cdots,\nu_m\}$,
	\begin{equation}\label{eq-thm-CVP-proof-2}
		\int_X\psi\mathrm{d}\mu=1-\int_X\chi_{[0]}\mathrm{d}\mu\leq\frac{1}{1+\tau}.
	\end{equation}
	For any $0<a<1$, let $\mu_a=a\nu_1+(1-a)\delta_{0^\infty}$.
	Then $\mu_a\in\mathcal{M}_f(X)$ and $\int_X\psi\mathrm{d}\mu_a=a$.
	Thus the set $\{\mu\in\mathcal{M}_f(X):\int_X\psi\mathrm{d}\mu=a\}$ is non-empty.

	From the ergodic decomposition and the affine character of entropy, for any $\mu\in\mathcal{M}_f(X)$ satisfying $\int_X\psi\mathrm{d}\mu=a$, there exist $0\leq\alpha\leq1$ and a probability measure $\pi$ on $\mathcal{M}_f^{erg}(X)\setminus\{\nu_1,\cdots,\nu_m\}$ such that $\alpha$ is the total mass assigned by the ergodic decomposition of $\mu$ to $\mathcal{M}_f^{erg}(X)\setminus\{\nu_1,\cdots,\nu_m\}$ and
	\begin{equation*}
		h_{\mu}(f)=(1-\alpha)h_{top}(f)+\alpha\int_{\mathcal{M}_f^{erg}(X)\setminus\{\nu_1,\cdots,\nu_m\}}h_{\omega}(f)\mathrm{d}\pi(\omega)\geq(1-\alpha)h_{top}(f).
	\end{equation*}
	Since $a>\gamma=1/(1+\tau)$, it follows from \eqref{eq-thm-CVP-proof-2} that $\alpha<1$.
	Thus
	\begin{equation*}
		h_{\mu}(f)\geq(1-\alpha)h_{top}(f)>0=h_{top}^B(f,R_{\psi}(a)).
	\end{equation*}
	Therefore,
	\begin{equation*}
		h_{top}^B(f,R_{\psi}(a))=0<\sup\bigg\{h_{\mu}(f):\mu\in\mathcal{M}_f(X),\int_X\psi\mathrm{d}\mu=a\bigg\},
	\end{equation*}
	and hence the conditional variational principle does not hold.
	This completes the proof for the case $0<\tau<\infty$.

	It remains to consider the case $\tau=\infty$.
	Let $(X,f)$ be the subshift constructed in the proof of Theorem \ref{thm-Bowen-upper-1} for $\tau=\infty$.
	Recall that $Z=Z_1\cup\cdots\cup Z_m$, and let $\nu_k$ be the unique ergodic measure of $Z_k$ for $1\leq k\leq m$.
	Consider $\psi=1-\chi_{[0]}\in C(X)$.
	Then $L_{\psi}=[0,1]$.
	Choose the constant $\sigma>0$ used in the construction of $X$ sufficiently large, and let
	\begin{equation*}
		\gamma=\frac{1}{1+\sigma}.
	\end{equation*}
	By the proof of Theorem \ref{thm-construction}, for any $x\in X\setminus\bigcup_{i\geq0}f^{-i}Z$, the upper density of $0$ in $x$ is at least $\sigma/(1+\sigma)$.
	Thus \eqref{eq-thm-CVP-proof-1} and \eqref{eq-thm-CVP-proof-2} remain valid with $\sigma$ in place of $\tau$.
	The preceding arguments therefore show that $R_{\psi}(a)=\varnothing$ and
	\begin{equation*}
		\sup\bigg\{h_{\mu}(f):\mu\in\mathcal{M}_f(X),\int_X\psi\mathrm{d}\mu=a\bigg\}>0=h_{top}^B(f,R_{\psi}(a))
	\end{equation*}
	for every $\gamma<a<1$.
	Since $\sigma$ can be chosen arbitrarily large, $\gamma$ can be chosen arbitrarily small.
	This completes the proof of Theorem \ref{thm-CVP}.
\end{proof}

\section{Transitive points}\label{Sect-trans}

\begin{proof}[{\bf Proof of Theorem \ref{thm-trans-1}}]
	We only consider the case $0<\tau<\infty$.
	The case $\tau=\infty$ is similar and hence omitted.
	Let $(X,f)$ be the dynamical system constructed in the proof of Theorem \ref{thm-Bowen-upper-2}.
	Then
	\begin{equation*}
		\mathrm{Trans}(f)\subset X\setminus\bigcup_{i\geq0}f^{-i}Z,
	\end{equation*}
	where $Z$ is the one-sided full shift used in the proof of Theorem \ref{thm-Bowen-upper-2}.
	By Theorem \ref{thm-construction} (4) and (5), for any $x\in\mathrm{Trans}(f)$, there exists $\mu\in V_f(x)$ such that
	\begin{equation*}
		h_{\mu}(f)\leq\bigg(\frac{1}{1+\tau}+\delta\bigg)h_{top}(f).
	\end{equation*}
	Applying Proposition \ref{Bowen-upper-prop1}, we conclude that
	\begin{equation}\label{eq-thm-trans-1-proof-1}
		h_{top}^B(f,\mathrm{Trans}(f))\leq\bigg(\frac{1}{1+\tau}+\delta\bigg)h_{top}(f)<h_{top}(f).
	\end{equation}

%	It remains to consider the case $\tau=\infty$.
%	Let $(X,f)$ be the dynamical system constructed in the proof of Theorem \ref{thm-Bowen-upper-2} for $\tau=\infty$.
%	As in the case $0<\tau<\infty$, we have
%	\begin{equation*}
%		\mathrm{Trans}(f)\subset X\setminus\bigcup_{i\geq0}f^{-i}Z.
%	\end{equation*}
%	By Theorem \ref{thm-construction} (4) and (5), for any $x\in\mathrm{Trans}(f)$, there exists $\mu\in V_f(x)$ such that
%	\begin{equation*}
%		h_{\mu}(f)\leq\delta h_{top}(f).
%	\end{equation*}
%	Applying Proposition \ref{Bowen-upper-prop1}, we conclude that
%	\begin{equation}\label{eq-thm-trans-1-proof-2}
%		h_{top}^B(f,\mathrm{Trans}(f))\leq\delta h_{top}(f)<h_{top}(f).
%	\end{equation}

	According to \cite[Theorem~3.1]{Tian-2016}, in both cases we have
	\begin{equation*}
		h_{top}^B(f,\mathrm{Rec}(f))=h_{top}(f).
	\end{equation*}
	Combining this equality with \eqref{eq-thm-trans-1-proof-1}, we obtain that
	\begin{equation*}
		h_{top}^B(f,\mathrm{Trans}(f))<h_{top}^B(f,\mathrm{Rec}(f)).
	\end{equation*}
	This completes the proof of Theorem \ref{thm-trans-1}.

\end{proof}

\begin{proof}[{\bf Proof of Theorem \ref{thm-trans-2}}]
	It suffices to consider the case $\tau<\infty$.
	By the traditional variational principle, for any $0<\gamma<h_{top}(f)$, there is an ergodic measure $\lambda\in\mathcal{M}_f^{erg}(X)$ such that $h_{\lambda}(f)>\gamma$.
	Take $q\in\mathbb{N}$ sufficiently large such that
	\begin{equation}\label{eq-thm-trans-2-proof-1}
		\gamma<\frac{q-1}{q}h_{\lambda}(f).
	\end{equation}
	Take $0<\sigma<1$ such that
	\begin{equation}\label{eq-thm-trans-2-proof-2}
		(1+3q(1+\tau))(1-\sigma)<1\Longleftrightarrow\sigma>3q(1+\tau)(1-\sigma).
	\end{equation}
	
	Consider a function $\Theta:[0,1/5]\rightarrow\mathbb{R}$ defined by
	\begin{equation*}
		\Theta(\theta)=\min\Bigg\{\sigma-\bigg(\frac{1-\sigma}{3q(1+\tau-2\theta)}+\frac{3q(1+\tau+2\theta)-1}{3q(1+\tau+2\theta)}\bigg),\frac{\sigma}{3q(1+\tau+2\theta)}-(1-\sigma)\Bigg\}.
	\end{equation*}
	Then $\Theta$ is continuous.
	Using \eqref{eq-thm-trans-2-proof-2}, we obtain that
	\begin{equation*}
		\Theta(0)=\min\Bigg\{\sigma-\frac{3q(1+\tau)-\sigma}{3q(1+\tau)},\frac{\sigma}{3q(1+\tau)}-(1-\sigma)\Bigg\}=\frac{(1+3q(1+\tau))(\sigma-1)+1}{3q(1+\tau)}>0.
	\end{equation*}
	Take $\eta>0$ sufficiently small such that $\Theta(0)>6\eta$.
	Then there exists $0<\theta_0<\min\{\eta,1/5\}$ such that for any $0<\theta<\theta_0$, we have $\Theta(\theta)>5\eta$.
	In other words, for any $0<\theta<\theta_0$, we have
	\begin{equation}\label{eq-thm-trans-2-proof-3}
		\frac{1-\sigma}{3q(1+\tau-2\theta)}+\frac{3q(1+\tau+2\theta)-1}{3q(1+\tau+2\theta)}<\sigma-5\eta
	\end{equation}
	and
	\begin{equation}\label{eq-thm-trans-2-proof-4}
		\frac{\sigma}{3q(1+\tau+2\theta)}>1-\sigma+5\eta.
	\end{equation}
	
	Fix two distinct ergodic measures $\mu,\nu\in\mathcal{M}_f^{erg}(X)$.
	Using Corollary \ref{Bowen-cor4}, we can find a $\varphi\in C_{\sigma}(X)$ such that $\inf_{x\in X}\varphi(x)=0$, $\|\varphi\|=\sup_{x\in X}\varphi(x)=1$ and
	\begin{equation*}
		\int_X\varphi\mathrm{d}\nu-\int_X\varphi\mathrm{d}\mu>\sigma\Big(\sup_{x\in X}\varphi(x)-\inf_{x\in X}\varphi(x)\Big)=\sigma.
	\end{equation*}
	Take an open neighborhood $\mathcal{W}$ of $\nu$ such that
	\begin{equation*}
		\int_X\varphi\mathrm{d}\omega-\int_X\varphi\mathrm{d}\mu>\sigma,\quad\forall\omega\in\mathcal{W}.
	\end{equation*}
	Hence
	\begin{equation}\label{eq-thm-trans-2-proof-5}
		\sigma\leq\int_X\varphi\mathrm{d}\omega\leq1\quad\text{for all}\;\omega\in\mathcal{W}\;\text{and}\quad 0\leq\int_X\varphi\mathrm{d}\mu\leq 1-\sigma.
	\end{equation}
	Similarly, for any $\pi\in\mathcal{M}_f^{erg}(X)\setminus\mathcal{W}$, there exist $\psi_{\pi}\in C_{\sigma}(X)$ and an open neighborhood $\mathcal{U}_{\pi}$ of $\pi$ such that
	\begin{itemize}
		\item $\inf_{x\in X}\psi_{\pi}(x)=0$ and $\sup_{x\in X}\psi_{\pi}(x)=1$;
		\item for any $\omega\in\mathcal{U}_{\pi}$, we have
		\begin{equation*}
			\sigma\leq\int_X\psi_{\pi}\mathrm{d}\nu\leq1\quad\text{and}\quad 0\leq\int_X\psi_{\pi}\mathrm{d}\omega\leq 1-\sigma.
		\end{equation*}	
	\end{itemize}

	The family $\{\mathcal{U}_{\pi}:\pi\in\mathcal{M}_f^{erg}(X)\setminus\mathcal{W}\}$ is an open cover of $\mathcal{M}_f^{erg}(X)\setminus\mathcal{W}$.
	Choose a countable subcover $\{\mathcal{U}_u\}_{u=1}^{\infty}$ of $\mathcal{M}_f^{\mathrm{erg}}(X)\setminus\mathcal{W}$.
	Let $\psi_u$ be the function corresponding to $\mathcal{U}_u$ for each $u\in\mathbb{N}$.
	Then
	\begin{itemize}
		\item $\mathcal{M}_f^{erg}(X)\subset\mathcal{W}\cup\Big(\bigcup_{u\geq1}\mathcal{U}_u\Big)$;
		\item for any $u\in\mathbb{N}$ and any $\omega\in\mathcal{U}_u$, we have
		\begin{equation}\label{eq-thm-trans-2-proof-6}
			\sigma\leq\int_X\psi_u\mathrm{d}\nu\leq 1 \quad\text{and}\quad 0\leq\int_X\psi_u\mathrm{d}\omega\leq1-\sigma.
		\end{equation}
	\end{itemize}
	
	Let
	\begin{equation*}
		\mathfrak{L}=\Bigg\{\omega\in\mathcal{M}(X):\bigg|\int_X\varphi\mathrm{d}\omega-\int_X\varphi\mathrm{d}\lambda\bigg|<\eta\Bigg\}.
	\end{equation*}
	Then $\mathfrak{L}$ is an open neighborhood in $\mathcal{M}(X)$ of $\lambda$.
	By \eqref{eq-thm-trans-2-proof-1} and Proposition \ref{pre-separation-prop1}, there exist $\delta^*>0$, $\varepsilon^*>0$ and $N^*\in\mathbb{N}$ such that for any $n\geq N^*$, there is a $(\delta^*,n,\varepsilon^*)$-separated set $\Gamma_n\subset X_{n,\mathfrak{L}}$ satisfying
	\begin{equation}\label{eq-thm-trans-2-proof-7}
		|\Gamma_n|\geq\exp\frac{qn\gamma}{q-1}.
	\end{equation}
	Let
	\begin{equation*}
		\begin{split}
			\mathfrak{M}&=\Bigg\{\omega\in\mathcal{M}(X):\bigg|\int_X\varphi\mathrm{d}\omega-\int_X\varphi\mathrm{d}\mu\bigg|<\eta\Bigg\},\\
			\mathfrak{N}&=\Bigg\{\omega\in\mathcal{M}(X):\bigg|\int_X\varphi\mathrm{d}\omega-\int_X\varphi\mathrm{d}\nu\bigg|<\eta\Bigg\}.
		\end{split}
	\end{equation*}
	Then $\mathfrak{M}$ and $\mathfrak{N}$ are open neighborhoods in $\mathcal{M}(X)$ of $\mu$ and $\nu$, respectively.
	Since $\mu$ and $\nu$ are ergodic, it follows from Birkhoff's ergodic theorem that we can take $n^*\geq N^*$ sufficiently large such that both $X_{n,\mathfrak{M}}$ and $X_{n,\mathfrak{N}}$ are non-empty for all $n\geq n^*$.
	Fix a generic point $y^{\nu}\in G_{\nu}$.
	
	Arbitrarily given $0<\theta<\theta_0$.
	Since $M(n,\varepsilon)$ is non-decreasing as $\varepsilon$ decreases, we have
	\begin{equation*}
		\tau=\sup_{\varepsilon>0}\liminf_{n\rightarrow\infty}\frac{M(n,\varepsilon)}{n}=\lim_{\varepsilon\rightarrow0}\liminf_{n\rightarrow\infty}\frac{M(n,\varepsilon)}{n}.
	\end{equation*}
	Therefore we can take $\varepsilon_0>0$ sufficiently small such that
	\begin{itemize}
		\item $\varepsilon_0<\varepsilon^*/5$;
		\item if two points $x,y\in X$ satisfy $d(x,y)\leq\varepsilon_0$, then $|\varphi(x)-\varphi(y)|\leq\eta$;
		\item for any $0<\varepsilon\leq2\varepsilon_0$, one has
		\begin{equation*}
			\bigg|\tau-\liminf_{n\rightarrow\infty}\frac{M(n,\varepsilon)}{n}\bigg|<\theta.
		\end{equation*}
	\end{itemize}
	We can take a decreasing sequences of positive real numbers $\{\varepsilon_i\}_{i=1}^{\infty}$ such that for any $i\in\mathbb{N}$,
	\begin{itemize}
		\item $\varepsilon_i<\varepsilon_{i-1}/2$;
		\item if two points $x,y\in X$ satisfy $d(x,y)\leq\varepsilon_i$, then $|\psi_u(x)-\psi_u(y)|\leq\eta$ for each $1\leq u\leq i$.
	\end{itemize}
	Then
	\begin{equation*}
		\bigg|\tau-\liminf_{n\rightarrow\infty}\frac{M(n,\varepsilon_i)}{n}\bigg|<\theta.
	\end{equation*}
	For $\varepsilon>0$, an \emph{$\varepsilon$-cover} is a cover of $X$ consisting of (open) balls with radii $\varepsilon$.
	For any $i\in\mathbb{N}_0$, let $N_i$ denote the minimal cardinality of $\varepsilon_i$-cover.
	Write $\varepsilon_{-1}=2\varepsilon_0$.
	
	Now we define sequences of positive integers $\{l_i\}_{i=0}^{\infty}$, $\{m_i\}_{i=0}^{\infty}$, $\{t_i\}_{i=0}^{\infty}$, $\{n_i\}_{i=0}^{\infty}$, $\{T_i\}_{i=0}^{\infty}$ and $\{M_i\}_{i=0}^{\infty}$.
	Take $l_0\geq\max\{n^*,1/\theta\}$ such that
	\begin{equation*}
		\bigg|\tau-\frac{M(l_0,\varepsilon_0)}{l_0}\bigg|<2\theta.
	\end{equation*}
	Let $m_0=M(l_0,\varepsilon_0)$ and $t_0=3q(l_0+m_0)$.
	Take $l_1>(l_0+m_0)^2$ such that
	\begin{equation*}
		\bigg|\tau-\frac{M(l_1,\varepsilon_1)}{l_1}\bigg|<2\theta.
	\end{equation*}
	Since $y^{\nu}\in G_{\nu}$, we can replace $l_1$ with a larger integer such that
	\begin{equation*}
		\left|\frac{1}{l_1}\sum_{i=0}^{l_1-1}\varphi(f^iy^{\nu})-\int_X\varphi\mathrm{d}\nu\right|<\eta
		\quad\text{and}\quad
		\left|\frac{1}{l_1}\sum_{i=0}^{l_1-1}\psi_1(f^iy^{\nu})-\int_X\psi_1\mathrm{d}\nu\right|<\eta.
	\end{equation*}
	Let $m_1=M(l_1,\varepsilon_1)$ and $t_1=3q(l_1+m_1)$.
	Take $n_0\geq N_0$ sufficiently large such that
	\begin{equation*}
		\frac{t_1}{n_0t_0}=\frac{l_1+m_1}{n_0(l_0+m_0)}<\theta
	\end{equation*}
	and
	\begin{equation*}
		\bigg|\tau-\frac{M(T_0,\varepsilon_1)}{T_0}\bigg|<2\theta,
	\end{equation*}
	where $T_0=3qn_0(l_0+m_0)=n_0t_0$.
	Let $M_0=M(T_0,\varepsilon_1)$.
	Assume that $\{l_i\}_{i=0}^k$, $\{m_i\}_{i=0}^k$, $\{t_i\}_{i=0}^k$, $\{n_i\}_{i=0}^{k-1}$, $\{T_i\}_{i=0}^{k-1}$ and $\{M_i\}_{i=0}^{k-1}$ have been defined.
	Take $l_{k+1}>(l_k+m_k)^2$ sufficiently large such that
	\begin{equation*}
		\bigg|\tau-\frac{M(l_{k+1},\varepsilon_{k+1})}{l_{k+1}}\bigg|<2\theta.
	\end{equation*}
	Since $y^{\nu}\in G_{\nu}$, we can replace $l_{k+1}$ with a larger integer such that
	\begin{equation*}
		\left|\frac{1}{l_{k+1}}\sum_{i=0}^{l_{k+1}-1}\varphi(f^iy^{\nu})-\int_X\varphi\mathrm{d}\nu\right|<\eta
		\quad\text{and}\quad
		\left|\frac{1}{l_{k+1}}\sum_{i=0}^{l_{k+1}-1}\psi_u(f^iy^{\nu})-\int_X\psi_u\mathrm{d}\nu\right|<\eta\;\text{for}\;1\leq u\leq k+1.
	\end{equation*}
	Let $m_{k+1}=M(l_{k+1},\varepsilon_{k+1})$ and $t_{k+1}=3q(l_{k+1}+m_{k+1})$.
	Take $n_k\geq N_k$ sufficiently large such that
	\begin{equation*}
		\frac{t_{k+1}}{n_kt_k}=\frac{l_{k+1}+m_{k+1}}{n_k(l_k+m_k)}<\theta
		\quad\text{and}\quad
		\bigg|\tau-\frac{M(T_k,\varepsilon_{k+1})}{T_k}\bigg|<2\theta,
	\end{equation*}
	where
	\begin{equation*}
		T_k=\sum_{i=0}^kn_it_i+\sum_{i=0}^{k-1}M_i.
	\end{equation*}
	Note that $T_k$ can be written as $T_k=n_kt_k+C_k$, where $C_k=\sum_{i=0}^{k-1}(n_it_i+M_i)$ is a constant dependent on sequences that have already been defined.
	Thus $n_k$ and $T_k$ are well-defined.
	Replacing with larger $n_k$ and $T_k$, we can also assume that
	\begin{equation*}
		1-\theta<\frac{n_kt_k}{T_k}<1\quad\text{and}\quad \frac{t_{k+1}}{T_k}<\theta.
	\end{equation*}
	Let $M_k=M(T_k,\varepsilon_{k+1})$.
	Inductively, we can define sequences of positive integers $\{l_i\}_{i=0}^{\infty}$, $\{m_i\}_{i=0}^{\infty}$, $\{t_i\}_{i=0}^{\infty}$, $\{n_i\}_{i=0}^{\infty}$, $\{T_i\}_{i=0}^{\infty}$ and $\{M_i\}_{i=0}^{\infty}$ satisfying that, for any $k\in\mathbb{N}_0$, we have
	\begin{equation}\label{eq-thm-trans-2-proof-8}
		\begin{gathered}
			n_k\geq N_k;\\
			(\tau-2\theta)l_k<m_k<(\tau+2\theta)l_k;\\
			\max\bigg\{\frac{t_{k+1}}{n_kt_k},\frac{t_{k+1}}{T_k}\bigg\}<\theta;\\
			(\tau-2\theta)T_k<M_k<(\tau+2\theta)T_k;\\
			(1-\theta)T_{k+1}<n_{k+1}t_{k+1}<T_{k+1}.
		\end{gathered}
	\end{equation}

	For $j\in\mathbb{N}_0$, $0\leq k\leq n_j-1$, define
	\begin{equation*}
		s(j,k)=
		\begin{cases}
			kt_0 & \text{if}\;j=0,\\
			\sum_{i=0}^{j-1}(n_it_i+M_i)+kt_j & \text{if}\;j\geq1,
		\end{cases}
	\end{equation*}
	and for $j\in\mathbb{N}_0$, $0\leq k\leq n_j-1$ and $0\leq r\leq 3q-1$, define
	\begin{equation*}
		\begin{cases}
			a(j,k,r)=s(j,k)+r(l_j+m_j),\\
			b(j,k,r)=s(j,k)+r(l_j+m_j)+l_j-1=a(j,k,r)+l_j-1.
		\end{cases}
	\end{equation*}
	Note that for $j\geq 1$, and $0\leq k\leq n_j-1$, we have
	\begin{equation*}
		s(j,k)=T_{j-1}+M_{j-1}+kt_j.
	\end{equation*}
	For any $j\in\mathbb{N}_0$, let $\Gamma_j=\Gamma_{l_j}\subset X_{l_j,\mathfrak{L}}$ be the $(\delta^*,l_j,\varepsilon^*)$-separated set satisfying \eqref{eq-thm-trans-2-proof-7}.
	For each $j\geq0$, fix $y^{\mu}_j\in X_{l_j,\mathfrak{M}}$.
	Clearly, $y^{\nu}\in X_{l_j,\mathfrak{N}}$.
	Since $n_j\geq N_j$, we can choose an $\varepsilon_j$-cover $\{B(z_j(r),\varepsilon_j)\}_{r=0}^{n_j-1}$ of $X$, allowing repetitions if $n_j>N_j$.
	Define
	\begin{equation*}
		\Omega_j(r)=\{y^{\mu}_j\}\times\Gamma_j^{q-1}\times\{y^{\nu}\}\times\Gamma_j^{q-1}\times\{z_j(r)\}\times\Gamma_j^{q-1}.
	\end{equation*}
	
	Let
	\begin{equation*}
		\Xi_l=\prod_{j=0}^{l-1}\prod_{r=0}^{n_j-1}\Omega_j(r)\quad\text{and}\quad\Xi=\prod_{j=0}^{\infty}\prod_{r=0}^{n_j-1}\Omega_j(r).
	\end{equation*}
	For every $\xi=(\xi_0,\xi_1,\cdots)\in\Xi$, we write
	\begin{equation*}
		\xi_j=(\xi_j^0,\cdots,\xi_j^{n_j-1})\in\prod_{r=0}^{n_j-1}\Omega_j(r)\quad\text{and}\quad\xi_j^r=(x_j^0(r),\cdots,x_j^{3q-1}(r))\in\Omega_j(r).
	\end{equation*}
	For every $k\in\mathbb{N}$, we define $E^{\xi}_k$ inductively starting with
	\begin{equation*}
		E^{\xi}_1=\bigcap_{r=0}^{n_0-1}\bigcap_{s=0}^{3q-1}\bigcap_{i=a(0,r,s)}^{b(0,r,s)}f^{-i}\Big(\overline{B(f^{i-a(0,r,s)}x_0^s(r),\varepsilon_0)}\Big).
	\end{equation*}
	From the non-uniform specification property, it follows that $E^{\xi}_1$ is a non-empty and closed set.
	Assume that $E^{\xi}_j$ have been defined for $1\leq j\leq k$.
	Define
	\begin{equation*}
		E^{\xi}_{k+1}=\bigg(\bigcup_{x\in E^{\xi}_k}\overline{B_{T_{k-1}}(x,\varepsilon_k)}\bigg)\cap\bigg(\bigcap_{r=0}^{n_k-1}\bigcap_{s=0}^{3q-1}\bigcap_{i=a(k,r,s)}^{b(k,r,s)}f^{-i}\Big(\overline{B(f^{i-a(k,r,s)}x_k^s(r),\varepsilon_k)}\Big)\bigg).
	\end{equation*}
	Clearly, the first set in this intersection is compact.
	From the non-uniform specification property, it follows that each $E^{\xi}_{k+1}$ is a non-empty and closed set.
	Inductively, we can define non-empty closed sets $E^{\xi}_k$ for each $k\in\mathbb{N}$.
	Similarly, for any $l\in\mathbb{N}$, we can define $E^{\zeta}_k$ for $\zeta\in\Xi_l$ and $1\leq k\leq l$.
	
	For $l\in\mathbb{N}$, define
	\begin{equation*}
		F_l=\bigcup_{\zeta\in\Xi_l}\bigg(\bigcap_{k=0}^{l-1}\bigcap_{r=0}^{n_k-1}\bigcap_{s=0}^{3q-1}\bigcap_{i=a(k,r,s)}^{b(k,r,s)}f^{-i}\Big(\overline{B(f^{i-a(k,r,s)}z_k^s(r),\varepsilon_{k-1})}\Big)\bigg)
	\end{equation*}
	and $F=\bigcap_{l\geq1}F_l$, where $\zeta=(\zeta_0,\zeta_1,\cdots,\zeta_{l-1})$, $\zeta_k=(\zeta_k^0,\cdots,\zeta_k^{n_k-1})$ and $\zeta_k^r=(z_k^0(r),\cdots,z_k^{3q-1}(r))$.
	Recall that $\varepsilon_k<\varepsilon_{k-1}/2$, then we have $\varepsilon_{k-1}\geq\sum_{j\geq k}\varepsilon_j$.
	Therefore, we have $F_l\supset\bigcup_{\zeta\in\Xi_l}E^{\zeta}_l$.
	The sets $F_l$ are non-empty compact sets, satisfying $F_{l+1}\subset F_l$, and hence $F$ is non-empty and compact.
	One can easily see that
	\begin{equation*}
		F\supset\bigcup_{\xi\in\Xi}\bigcap_{k\geq1}E^{\xi}_k.
	\end{equation*}
	We proceed to verify that
	\begin{equation}\label{eq-thm-trans-2-proof-9}
		F\subset X\setminus\bigg(\bigcup_{\omega\in\mathcal{M}_f^{erg}(X)}G_{\omega}\bigg).
	\end{equation}
	Fix $x\in F$ and $\omega\in\mathcal{M}_f^{erg}(X)$.
	
	{\bf Case 1.}
	If $\omega\in\mathcal{W}$, then for any $j\geq2$ and any $0\leq k\leq n_j-1$, we have $f^{a(j,k,0)}x\in B_{l_j}(y^{\mu}_j,\varepsilon_0)$.
	Recall that $|\varphi(x)-\varphi(y)|\leq\eta$ whenever $d(x,y)\leq\varepsilon_0$.
	Since $y^{\mu}_j\in X_{l_j,\mathfrak{M}}$, for any $j\geq2$, we have
	\begin{equation}\label{eq-thm-trans-2-proof-10}
		\frac{1}{l_j}\sum_{i=a(j,k,0)}^{b(j,k,0)}\varphi(f^ix)\leq\frac{1}{l_j}\sum_{i=0}^{l_j-1}\varphi(f^iy^{\mu}_j)+\eta<\int_X\varphi\mathrm{d}\mu+2\eta.
	\end{equation}
	Recall that $\sup_{x\in X}\varphi(x)=1$ and $\theta<\theta_0<\eta$.
	Combining \eqref{eq-thm-trans-2-proof-3}, \eqref{eq-thm-trans-2-proof-5}, \eqref{eq-thm-trans-2-proof-8} and \eqref{eq-thm-trans-2-proof-10}, for any $j\geq2$, using the arguments similar to the proofs of Theorem \ref{Bowen-thm2}, we obtain that
	\begin{equation*}
		0\leq \frac{l_j}{3q(l_j+m_j)}-\frac{n_jl_j}{T_j}\leq\theta
	\end{equation*}
	and
	\begin{equation*}
		\begin{split}
			\frac{1}{T_j}\sum_{i=0}^{T_j-1}\varphi(f^ix)&<\frac{n_jl_j}{T_j}\int_X\varphi\mathrm{d}\mu+\frac{T_j-n_jl_j}{T_j}+2\eta\\
			&\leq\frac{1}{3q(1+\tau-2\theta)}\int_X\varphi\mathrm{d}\mu+\bigg(\theta+\frac{3q(1+\tau+2\theta)-1}{3q(1+\tau+2\theta)}\bigg)+2\eta\\
			&\leq\frac{1-\sigma}{3q(1+\tau-2\theta)}+\frac{3q(1+\tau+2\theta)-1}{3q(1+\tau+2\theta)}+3\eta<\sigma-2\eta\leq\int_X\varphi\mathrm{d}\omega-2\eta.
		\end{split}
	\end{equation*}
	Hence $x\notin G_{\omega}$.

	{\bf Case 2.}
	If $\omega\in\mathcal{U}_u$ for some $u\in\mathbb{N}$, then for any $j\geq u+2$ and any $0\leq k\leq n_j-1$, we have $f^{a(j,k,q)}x\in B_{l_j}(y^{\nu},\varepsilon_u)$.
	Recall that $|\psi_u(x)-\psi_u(y)|\leq\eta$ whenever $d(x,y)\leq\varepsilon_u$.
	From the choice of $l_j$, it follows that
	\begin{equation}\label{eq-thm-trans-2-proof-11}
		\frac{1}{l_j}\sum_{i=a(j,k,q)}^{b(j,k,q)}\psi_u(f^ix)\geq\frac{1}{l_j}\sum_{i=0}^{l_j-1}\psi_u(f^iy^{\nu})-\eta>\int_X\psi_u\mathrm{d}\nu-2\eta.
	\end{equation}
	Recall that $\inf_{x\in X}\psi_u(x)=0$ and $\theta<\theta_0<\eta$.
	Combining \eqref{eq-thm-trans-2-proof-4}, \eqref{eq-thm-trans-2-proof-6}, \eqref{eq-thm-trans-2-proof-8} and \eqref{eq-thm-trans-2-proof-11}, for any $j\geq u+2$, using an argument similar to the proofs of Theorem \ref{Bowen-thm2}, we obtain that
	\begin{equation*}
		0\leq \frac{l_j}{3q(l_j+m_j)}-\frac{n_jl_j}{T_j}\leq\theta
	\end{equation*}
	and
	\begin{equation*}
		\begin{split}
			\frac{1}{T_j}\sum_{i=0}^{T_j-1}\psi_u(f^ix)&>\frac{n_jl_j}{T_j}\int_X\psi_u\mathrm{d}\nu-2\eta\\
			&\geq\bigg(\frac{1}{3q(1+\tau+2\theta)}-\theta\bigg)\cdot\int_X\psi_u\mathrm{d}\nu-2\eta\\
			&\geq\frac{\sigma}{3q(1+\tau+2\theta)}-3\eta>1-\sigma+2\eta\geq\int_X\psi_u\mathrm{d}\omega+2\eta.
		\end{split}
	\end{equation*}
	Hence $x\notin G_{\omega}$.
	Therefore \eqref{eq-thm-trans-2-proof-9} holds.
	
	Now we show that $F\subset\mathrm{Trans}(f)$.
	Fix $x\in F$ and $y\in X$.
	For any $j\in\mathbb{N}$, since $\{B(z_j(r),\varepsilon_j)\}_{r=0}^{n_j-1}$ covers $X$, there is an integer $0\leq r_j\leq n_j-1$ such that $y\in B(z_j(r_j),\varepsilon_j)$.
	Then from the construction of $F$, we have
	\begin{equation*}
		d(f^{a(j,r_j,2q)}x,y)\leq d(f^{a(j,r_j,2q)}x,z_j(r_j))+d(z_j(r_j),y)<2\varepsilon_{j-1}.
	\end{equation*}
	Since $\varepsilon_j\rightarrow0$ and $a(j,r_j,2q)\rightarrow\infty$ as $j\rightarrow\infty$, we have $y\in\omega(x,f)$.
	This implies that $F\subset\mathrm{Trans}(f)$.
	Combining with \eqref{eq-thm-trans-2-proof-9}, we conclude that
	\begin{equation*}
		F\subset \mathrm{Trans}(f)\setminus\bigg(\bigcup_{\omega\in\mathcal{M}_f^{erg}(X)}G_{\omega}\bigg).
	\end{equation*}

	Finally we use Lemma \ref{Bowen-lemma1} to estimate the Bowen topological entropy of $F$.
	For $\xi\in\Xi$ and $k\in\mathbb{N}$, define
	\begin{equation*}
		Y_k^{\xi}=\bigcap_{j=0}^{k-1}\bigcap_{r=0}^{n_j-1}\bigcap_{s=0}^{3q-1}\bigcap_{i=a(j,r,s)}^{b(j,r,s)}f^{-i}\Big(\overline{B(f^{i-a(j,r,s)}x_j^s(r),\varepsilon_{j-1})}\Big).
	\end{equation*}
	Then $E_k^{\xi}\subset Y_k^{\xi}$.
	Thus the sets $Y_k^{\xi}$ are non-empty compact sets, satisfying $Y_{k+1}^{\xi}\subset Y_k^{\xi}$.
	Consequently,
	\begin{equation*}
		Y^{\xi}:=\bigcap_{k\geq1}Y_k^{\xi}
	\end{equation*}
	is a non-empty compact subset of $F$.
	Arbitrarily take $x_{\xi}\in Y^{\xi}$ for each $\xi\in\Xi$.
	Let $Y=\{x_{\xi}:\xi\in\Xi\}$.
	Recall that $t_j=3q(l_j+m_j)$.
	Similar to the proof of Theorem \ref{Bowen-thm2}, for any $\xi,\zeta\in\Xi$, we have $f^{s(j,k)}x_{\xi}$ and $f^{s(j,k)}x_{\zeta}$ are $(t_j,\varepsilon_0)$-separated whenever $\xi_j^k\neq\zeta_j^k$ for some $j\in\mathbb{N}_0$ and $0\leq k\leq n_j-1$.
	Using \eqref{eq-thm-trans-2-proof-8}, we have
	\begin{equation*}
		\liminf_{j\rightarrow\infty}\frac{\sum_{i=0}^jn_it_i}{s(j+1,1)-1}\geq\liminf_{j\rightarrow\infty}\frac{T_j}{T_j+M_j-1+t_{j+1}}\cdot\frac{n_jt_j}{T_j}\geq\frac{1-\theta}{1+\tau+3\theta}
	\end{equation*}
	and 
	\begin{equation*}
		|\Omega_j(r)|=|\Gamma_j|^{3(q-1)}\geq\exp(3ql_j\gamma)=\exp\bigg(\frac{l_j}{l_j+m_j}\cdot(t_j\cdot\gamma)\bigg)\geq\exp\frac{t_j\cdot\gamma}{1+\tau+2\theta}.
	\end{equation*}
	Applying Lemma \ref{Bowen-lemma1}, we conclude that
	\begin{equation*}
		h_{top}^B(f,F)\geq h_{top}^B(f,\overline{Y})\geq\frac{(1-\theta)\cdot\gamma}{(1+\tau+3\theta)(1+\tau+2\theta)}.
	\end{equation*}

	Note that
	\begin{equation*}
		\mathrm{IR}(f)\cup\bigg(\bigcup_{\omega\in\mathcal{M}_f(X)\setminus\mathcal{M}_f^{erg}(X)}G_{\omega}\bigg)=X\setminus\bigg(\bigcup_{\omega\in\mathcal{M}_f^{erg}(X)}G_{\omega}\bigg).
	\end{equation*}
	From Proposition \ref{pre-entropy-prop1}, it follows that
	\begin{equation*}
		\begin{split}
			\max\,&\Bigg\{h_{top}^B\Big(f,\mathrm{IR}(f)\cap\mathrm{Trans}(f)\Big),h_{top}^B\bigg(f,\bigcup_{\omega\in\mathcal{M}_f(X)\setminus\mathcal{M}_f^{erg}(X)}G_{\omega}\cap\mathrm{Trans}(f)\bigg)\Bigg\}\\
			&=h_{top}^B\Bigg(f,\mathrm{Trans}(f)\setminus\bigg(\bigcup_{\omega\in\mathcal{M}_f^{erg}(X)}G_{\omega}\bigg)\Bigg)\geq h_{top}^B(f,F)\geq\frac{(1-\theta)\cdot\gamma}{(1+\tau+3\theta)(1+\tau+2\theta)}.
		\end{split}
	\end{equation*}
	By the arbitrariness of $0<\gamma<h_{top}(f)$ and $0<\theta\leq\theta_0$, we conclude that
	\begin{equation*}
		\max\Bigg\{h_{top}^B\Big(f,\mathrm{IR}(f)\cap\mathrm{Trans}(f)\Big),h_{top}^B\bigg(f,\bigcup_{\omega\in\mathcal{M}_f(X)\setminus\mathcal{M}_f^{erg}(X)}G_{\omega}\cap\mathrm{Trans}(f)\bigg)\Bigg\}\geq \frac{1}{(1+\tau)^2}h_{top}(f).
	\end{equation*}
	This completes the proof of Theorem \ref{thm-trans-2}.
\end{proof}

\section{Growth of periodic orbits}\label{Sect-growth}

\begin{proof}[{\bf Proof of Theorem \ref{thm-growth-1}}]
	It suffices to consider the case $0\leq\tau<\infty$.
	Arbitrarily given $\kappa>0$ and $0<\gamma<h_{top}(f)$, we aim to show that
	\begin{equation*}
		\liminf_{n\rightarrow\infty}\frac{1}{n}\log|\mathrm{Per}_n(f)|\geq\frac{\gamma}{1+\tau+\kappa}.
	\end{equation*}
	Recall that $s_n(\varepsilon)$ denotes the maximal cardinality of $(n,\varepsilon)$-separated sets.
	We choose $\varepsilon_0>0$ sufficiently small and $N\in\mathbb{N}$ sufficiently large such that
	\begin{equation*}
		s_n(3\varepsilon_0)\geq\exp(n\gamma),\quad\forall n\geq N.
	\end{equation*}
	Note that $s_n(3\varepsilon)\geq s_n(3\varepsilon_0)$ for all $0<\varepsilon<\varepsilon_0$.
	We may make $\varepsilon_0$ smaller such that there is a positive integer $l_0\geq N$ satisfying
	\begin{equation*}
		\bigg|\frac{M(l_0,\varepsilon_0)}{l_0}-\tau\bigg|\leq\kappa.
	\end{equation*}
	Write $m_0=M(l_0,\varepsilon_0)$ and $t_0=l_0+m_0$, then
	\begin{equation*}
		(1+\tau-\kappa)l_0\leq t_0\leq(1+\tau+\kappa)l_0.
	\end{equation*}
	Let $\Gamma$ be an $(l_0,3\varepsilon_0)$-separated set with $|\Gamma|=s_{l_0}(3\varepsilon_0)$.
	For any $n\geq 2t_0$, we can write $n=qt_0+r$, where $q,r\in\mathbb{Z}$ with $q\geq 2$ and $0\leq r\leq t_0-1$.
	By the non-uniform periodic specification property, for any $x=(x_0,\cdots,x_{q-1})\in\Gamma^q$, there is a $z=z(x)\in\mathrm{Per}_n(f)$ such that for any $0\leq j\leq q-1$, we have
	\begin{equation*}
		d(f^iz,f^{i-jt_0}x_j)\leq\varepsilon_0\quad\text{for}\quad jt_0\leq i\leq jt_0+l_0-1.
	\end{equation*}
	Since $\Gamma$ is $(l_0,3\varepsilon_0)$-separated, if $x\neq y$, then we have $z(x)\neq z(y)$.
	This implies that
	\begin{equation*}
		|\mathrm{Per}_n(f)|\geq|\Gamma^q|=(s_{l_0}(3\varepsilon_0))^q\geq\exp(ql_0\gamma).
	\end{equation*}
	Hence
	\begin{equation*}
		\frac1n\log|\mathrm{Per}_n(f)|\geq\frac{ql_0\gamma}{qt_0+r}\geq\frac{q\gamma}{(q+1)(1+\tau+\kappa)}.
	\end{equation*}
	Therefore
	\begin{equation*}
		\liminf_{n\rightarrow\infty}\frac{1}{n}\log|\mathrm{Per}_n(f)|\geq\liminf_{q\rightarrow\infty}\frac{q\gamma}{(q+1)(1+\tau+\kappa)}=\frac{\gamma}{1+\tau+\kappa}.
	\end{equation*}
	By the arbitrariness of $\kappa$ and $\gamma$, we conclude that
	\begin{equation*}
		\liminf_{n\rightarrow\infty}\frac{1}{n}\log|\mathrm{Per}_n(f)|\geq\frac{1}{1+\tau}h_{top}(f).
	\end{equation*}
	The inequality
	\begin{equation*}
		\limsup_{n\rightarrow\infty}\frac{1}{n}\log|\mathrm{Per}_n(f)|\leq h_{top}(f)
	\end{equation*}
	follows directly from \cite[Theorem 3.2 (3)]{Lin-Tian-Yu-2024}, as positive expansiveness implies pseudo-expansiveness.
	This completes the proof of Theorem \ref{thm-growth-1}.
\end{proof}

Before the proof of Theorem \ref{thm-growth-2}, let us introduce a lemma, which is used to estimate the upper bound for the growth of orbits.

%lemma1
\begin{lemma}[{\cite[Lemma A.2]{Crovisier-Yang-Zhang}}]\label{growth-lemma1}
	For any $\varepsilon,\eta>0$ and $\mu\in\mathcal{M}_f(X)$, there exist $\kappa>0$ and $n_0\in\mathbb{N}$ satisfying the following condition:
	if $\Gamma$ is an $(n,\varepsilon)$-separated set with $n\geq n_0$ and
	\begin{equation*}
		\rho\bigg(\mu,\frac{1}{|\Gamma|}\sum_{x\in\Gamma}\delta^n_x\bigg)<\kappa,
	\end{equation*}
	then $|\Gamma|\leq\exp(n(h_{\mu}(f)+\eta))$.
\end{lemma}

\begin{proof}[{\bf Proof of Theorem \ref{thm-growth-2}}]
	We only consider the case $0<\tau<\infty$.
	The case $\tau=\infty$ is similar and hence omitted.
	Let $(X,f)$ be the subshift defined in Theorem \ref{thm-Bowen-upper-1} for $m=1$, with $\delta/2$ in place of $\delta$.
	Then $(X,f)$ satisfies the non-uniform periodic $\tau$-specification property and (1) holds.
	For any $n\in\mathbb{N}$, define
	\begin{equation}\label{eq-thm-growth-2-proof-1}
		\omega_n=\frac{1}{|\mathrm{Per}_n(f)|}\sum_{x\in\mathrm{Per}_n(f)}\delta_x^n.
	\end{equation}
	Then $\omega_n\in\mathcal{M}_f(X)$.
	Let $\mu$ be an accumulation point of $\{\omega_n\}_{n=1}^{\infty}$.
	Then $\mu\in\mathcal{M}_f(X)$.
	Recall that $X$ is the subshift constructed in Theorem \ref{thm-construction} containing $Z$, where $Z$ is a strictly ergodic subshift that has large entropy given in Theorem \ref{thm-Bowen-upper-1}.
	Thus $\mathcal{M}_f^{co}(X)\subset\mathcal{M}_0$.
	This implies that \eqref{eq-thm-construction-proof-6} holds for all $\mu\in\mathcal{M}_f^{co}(X)$.
	Since $\omega_n$ is a convex combination of periodic measures, for any $n\in\mathbb{N}$, we obtain that $\omega_n$ satisfies \eqref{eq-thm-construction-proof-6}.
	According to the proof of Theorem \ref{thm-construction} (5), we obtain that
	\begin{equation*}
		h_{\mu}(f)\leq\log2+\frac{1}{1+\tau}\log A_j.
	\end{equation*}
	Similar to Theorem \ref{thm-Bowen-upper-1}, when $A_j$ is sufficiently large, one can show that
	\begin{equation*}
		h_{\mu}(f)\leq\bigg(\frac{1}{1+\tau}+\frac{\delta}{2}\bigg)\cdot h_{top}(f).
	\end{equation*}
	Let $\eta=(\delta\cdot h_{top}(f))/2>0$ and $\varepsilon=1/3$, then $\mathrm{Per}_n(f)$ is an $(n,\varepsilon)$-separated set.
	Let $K$ be the set of accumulation points of $\{\omega_n\}_{n=1}^{\infty}$.
	Then $K$ is a non-empty compact subset of $\mathcal{M}_f(X)$.
	For every $\mu\in K$, let $\kappa_{\mu}>0$ and $n_{\mu}\in\mathbb{N}$ be the constants given by Lemma \ref{growth-lemma1}.
	By the compactness of $K$, there exist $\mu_1,\ldots,\mu_k\in K$ such that $\{B(\mu_i,\kappa_{\mu_i})\}_{i=1}^k$ covers $K$.
	It follows from the definition of $K$ that $\omega_n$ belongs to this union for all sufficiently large $n$.
	By Lemma \ref{growth-lemma1}, for any $n\geq\max\{n_{\mu_i}:1\leq i\leq k\}$, we conclude that
	\begin{equation*}
		\limsup_{n\rightarrow\infty}\frac{1}{n}\log|\mathrm{Per}_n(f)|\leq\bigg(\frac{1}{1+\tau}+\delta\bigg)\cdot h_{top}(f).
	\end{equation*}
%
%	It remains to consider the case $\tau=\infty$.
%	Let $(X,f)$ be the subshift defined in Theorem \ref{thm-Bowen-upper-1} for $m=1$ and $\tau=\infty$, with $\delta/2$ in place of $\eta$.
%	Then $(X,f)$ satisfies the non-uniform periodic $\infty$-specification property and (1) holds.
%	The proof is similar to the case $0<\tau<\infty$.
%	Indeed, replacing \eqref{eq-thm-construction-proof-6} by \eqref{eq-thm-construction-proof-infinity-4}, it follows from the proof of Theorem \ref{thm-construction} (5) that
%	\begin{equation*}
%		h_{\mu}(f)\leq\frac{\delta}{2}h_{top}(f)
%	\end{equation*}
%	for every accumulation point $\mu$ of the sequence $\{\omega_n\}_{n\geq1}$, where $\omega_n$ is the measure defined by \eqref{eq-thm-trans-2-proof-7}.
%	Applying Lemma \ref{growth-lemma1} as above, we obtain that
%	\begin{equation*}
%		\limsup_{n\rightarrow\infty}\frac{1}{n}\log|\mathrm{Per}_n(f)|\leq\delta\cdot h_{top}(f)<h_{top}(f).
%	\end{equation*}
	This completes the proof of Theorem \ref{thm-growth-2}.
\end{proof}

\section{Intermediate entropy property}\label{Sect-intermediate}

\begin{proof}[{\bf Proof of Theorem \ref{thm-intermediate-1}}]
	We only consider the case that $0<\tau<\infty$.
	Let $(X,f)$ be the subshift defined in Theorem \ref{thm-Bowen-upper-1} for the case $m=1$.
	Recall that $X$ is the subshift given in Theorem \ref{thm-construction} associated to a strictly ergodic subshift $Z$ that has large entropy.
	Then $\mathcal{M}_0=\mathcal{M}_f^{erg}(X)\setminus\{\nu_1\}$, where $\nu_1$ is the unique measure of maximal entropy.
	By Theorem \ref{thm-construction} (5), if $\mu\in\mathcal{M}_0$, then
	\begin{equation*}
		h_{\mu}(f)\leq\bigg(\frac{1}{1+\tau}+\delta\bigg)\cdot h_{top}(f)<h_{top}(f).
	\end{equation*}
	Hence for any $(1/(1+\tau)+\delta)\cdot h_{top}(f)<h<h_{top}(f)$, there is no ergodic measure with measure-theoretic entropy $h$, which implies that the intermediate entropy property does not hold.
	This completes the proof of Theorem \ref{thm-intermediate-1}.
\end{proof}

\begin{proof}[{\bf Proof of Theorem \ref{thm-intermediate-2}}]
	First, we consider the case $0<\tau<\infty$.
	Let $(X,f)$ be the subshift defined in Theorem \ref{thm-Bowen-upper-2}.
	Recall that $X$ is the subshift constructed in Theorem \ref{thm-construction} associated to $\mathcal{A}=\{1,\cdots,A\}$ and $Z=\mathcal{A}^{\mathbb{N}_0}$, where $A>2^{\frac{1+\tau}{\tau}}$.
	In other words, $X=X(\mathcal{F})$, where $\mathcal{F}$ is the set of forbidden words over $\widetilde{\mathcal{A}}$ consisting of all words of the form $x_1\cdots x_s0^tx_{s+1}$, where $x_i>0$ and $\tau \cdot s>t$, and $\widetilde{\mathcal{A}}=\{0\}\cup\mathcal{A}=\{0,1,\cdots,A\}$.
	See Theorem \ref{thm-construction} and Theorem \ref{thm-Bowen-upper-2}.
	We may retake a larger $A$ in this proof.

	It has been shown in Theorem \ref{thm-Bowen-upper-2} that $h_{top}(f)=h_{top}(X)=\log A$.
	Since $Z=\mathcal{A}^{\mathbb{N}_0}$ is a one-sided full shift, it is known that $h_{top}(Z)=\log A$ and $(Z,f)$ satisfies the intermediate entropy property.
	Thus for any $0\leq h<h_{top}(f)=\log A$, there is a $\nu\in\mathcal{M}_f^{erg}(Z)\subset\mathcal{M}_f^{erg}(X)$ with $h_{\nu}(f)=h$, which implies that $(X,f)$ satisfies the intermediate entropy property.

	Now we consider the continuous function $\phi=(2\tau^{-1}\log A)\cdot(1-\chi_{[0]})$.
	Clearly $\phi|_{Z}\equiv (2\tau^{-1}\log A)$ and $\phi(x)\geq 0$ for all $x\in X$.
	By Theorem \ref{thm-construction} (5), for any $\mu\in\mathcal{M}_f^{erg}(X)\setminus\mathcal{M}_f^{erg}(Z)$, we have
	\begin{equation*}
		h_{\mu}(f)\leq\log2+\frac{1}{1+\tau}\log A
	\end{equation*}
	and thus
	\begin{equation*}
		h_{\mu}(f)+\int_X\phi\mathrm{d}\mu\leq\log2+\frac{\log A}{1+\tau}+\frac{2\tau^{-1}\log A}{1+\tau}=\log2+\frac{2+\tau}{(1+\tau)\tau}\log A.
	\end{equation*}
	On the other hand, for any $\nu\in\mathcal{M}_f^{erg}(Z)$, we have
	\begin{equation*}
		h_{\nu}(f)+\int_X\phi\mathrm{d}\nu\geq\int_X\phi\mathrm{d}\nu=2\tau^{-1}\log A.
	\end{equation*}
	When $A$ is sufficiently large, we have
	\begin{equation*}
		\log2+\frac{2+\tau}{(1+\tau)\tau}\log A<2\tau^{-1}\log A.
	\end{equation*}
	As a result,
	\begin{equation*}
		\sup_{\mu\in\mathcal{M}_f^{erg}(X)\setminus\mathcal{M}_f^{erg}(Z)}\bigg\{h_{\mu}(f)+\int_X\phi\mathrm{d}\mu\bigg\}<\inf_{\nu\in\mathcal{M}_f^{erg}(Z)}\bigg\{h_{\nu}(f)+\int_X\phi\mathrm{d}\nu\bigg\}.
	\end{equation*}
	Hence the intermediate pressure property does not hold.

	Now assume that $\tau=\infty$.
	Let $(X,f)$ be the subshift constructed in Theorem \ref{thm-Bowen-upper-2} for $\tau=\infty$, with $Z=\mathcal{A}^{\mathbb{N}_0}$ and $|\mathcal{A}|=A$.
	Since $h_{top}(f)=h_{top}(Z)=\log A$, the same argument as above shows that $(X,f)$ satisfies the intermediate entropy property.
	Let $\sigma>0$ be the constant used in the construction of $X$.
	By the proof of Theorem \ref{thm-construction}, for any $\mu\in\mathcal{M}_f^{erg}(X)\setminus\mathcal{M}_f^{erg}(Z)$, we have
	\begin{equation*}
		h_\mu(f)\leq \log 2+\frac{1}{1+\sigma}\log A\quad\text{and}\quad\int_X(1-\chi_{[0]})\mathrm{d}\mu\leq \frac{1}{1+\sigma}.
	\end{equation*}
	Choose $a>((1+\sigma)\log 2+\log A)/\sigma$ and let $\phi=a(1-\chi_{[0]})$.
	Then
	\begin{equation*}
		\sup_{\mu\in\mathcal{M}_f^{erg}(X)\setminus\mathcal{M}_f^{erg}(Z)}P_\mu(f,\phi)\leq \log2+\frac{1}{1+\sigma}\log A+\frac{a}{1+\sigma}<a\leq\inf_{\nu\in\mathcal{M}_f^{erg}(Z)}P_\nu(f,\phi).
	\end{equation*}
	Since $P_{\delta_{0^\infty}}(f,\phi)=0$ and $P(f,\phi)\geq a$, the intermediate pressure property does not hold.
	This completes the proof of Theorem \ref{thm-intermediate-2}.
\end{proof}

\begin{proof}[{\bf Proof of Theorem \ref{thm-intermediate-3}}]
	The case of $\tau=0$ is a consequence of \cite[Lemma 3.5]{Lin-Tian-Yu-2024} and \cite[Theorem 1.4]{Sun-2025}.
	We assume that $0<\tau<\infty$.
	If $\gamma=0$, choose $\widetilde\gamma>0$ sufficiently small and set $\widetilde\beta=\beta-2\widetilde\gamma$ so that
	\begin{equation*}
		\frac{\tau}{1+\tau}h_{top}(f)<\widetilde\beta<h_{top}(f)-\widetilde\gamma.
	\end{equation*}
	The construction below with $(\widetilde\gamma,\widetilde\beta)$ in place of $(\gamma,\beta)$ gives
	\begin{equation*}
		0<\widetilde\gamma<h_{top}(f|_Y)<\widetilde\gamma+\widetilde\beta=\beta-\widetilde\gamma<\beta.
	\end{equation*}
	Thus it suffices to consider $0<\gamma<h_{top}(f)$.

	The construction of $Y$ and the proof of the first inequality are similar to Proposition \ref{Bowen-prop5}.
	Recall that $s_n(\varepsilon)$ denotes the maximal cardinality of $(n,\varepsilon)$-separated sets.
	Take $0<\theta^*<\tau/2$ sufficiently small such that
	\begin{equation}\label{eq-thm-intermediate-3-proof-1}
		\gamma<\frac{1}{1+\tau+4\theta^*}h_{top}(f)
	\end{equation}
	and
	\begin{equation}\label{eq-thm-intermediate-3-proof-2}
		\frac{(1+\tau+4\theta^*)\gamma}{1+\tau-2\theta^*}+\frac{(\tau+2\theta^*)(h_{top}(f)+\theta^*)}{1+\tau+2\theta^*}<\gamma+\beta.
	\end{equation}
	Since $f$ is positively expansive, by \eqref{eq-thm-intermediate-3-proof-1} we can fix $\varepsilon^*>0$ and $n^*$ sufficiently large such that $\varepsilon^*$ is less than the expansive constant and for any $n\geq n^*$, we have
	\begin{equation}\label{eq-thm-intermediate-3-proof-3}
		s_n(\varepsilon^*)\geq\exp(n(1+\tau+3\theta^*)\gamma)\quad\text{and}\quad\exp(n(1+\tau+4\theta^*)\gamma)>\exp(n(1+\tau+3\theta^*)\gamma)+1.
	\end{equation}
	Since $M(n,\varepsilon)$ is non-decreasing as $\varepsilon$ decreases, we have
	\begin{equation*}
		\tau=\sup_{\varepsilon>0}\liminf_{n\rightarrow\infty}\frac{M(n,\varepsilon)}{n}=\lim_{\varepsilon\rightarrow0}\liminf_{n\rightarrow\infty}\frac{M(n,\varepsilon)}{n}.
	\end{equation*}
	Take $\varepsilon_0>0$ sufficiently small such that
	\begin{itemize}
		\item $\varepsilon_0<\varepsilon^*/5$;
		\item for any $0<\varepsilon\leq\varepsilon_0$, one has
		\begin{equation*}
			\left|\tau-\liminf_{n\rightarrow\infty}\frac{M(n,\varepsilon)}{n}\right|<\theta^*.
		\end{equation*}
	\end{itemize}
	Take a positive integer $n_0\geq n^*$ sufficiently large such that
	\begin{equation*}
		\bigg|\tau-\frac{M(n_0,\varepsilon_0)}{n_0}\bigg|<2\theta^*.
	\end{equation*}
	Let $m_0=M(n_0,\varepsilon_0)$, then
	\begin{equation}\label{eq-thm-intermediate-3-proof-4}
		(\tau-2\theta^*)n_0\leq m_0\leq (\tau+2\theta^*)n_0.
	\end{equation}
	We retake larger $n_0$ such that
	\begin{equation}\label{eq-thm-intermediate-3-proof-5}
		s_{m_0}(\varepsilon_0)\leq\exp\big(m_0(h_{top}(f)+\theta^*)\big).
	\end{equation}
	For $j\in\mathbb{N}_0$, let $a(j)=j(n_0+m_0)$ and $b(j)=j(n_0+m_0)+n_0-1$.

	By \eqref{eq-thm-intermediate-3-proof-3}, we can choose an $(n_0,\varepsilon^*)$-separated set $\Omega$ such that
	\begin{equation}\label{eq-thm-intermediate-3-proof-6}
		\exp(n_0(1+\tau+3\theta^*)\gamma)\leq|\Omega|\leq\exp(n_0(1+\tau+4\theta^*)\gamma).
	\end{equation}
	For any $\xi=(x_0,x_1,\cdots)\in\Omega^{\mathbb{N}_0}$, define
	\begin{equation*}
		E^{\xi}_k=\bigcap_{j=0}^{k-1}\bigcap_{i=a(j)}^{b(j)}f^{-i}\Big(\overline{B(f^{i-a(j)}x_j,\varepsilon_0)}\Big).
	\end{equation*}
	From the non-uniform specification property, it follows that each $E^{\xi}_k$ is a non-empty closed set.
	Similarly, for any $l\in\mathbb{N}$, we can define $E^{\zeta}_k$ for $\zeta\in\Omega^l$ and $1\leq k\leq l$.

	For $l\in\mathbb{N}$, define
	\begin{equation*}
		F_l=\bigcup_{\zeta\in\Omega^l}E^{\zeta}_l
	\end{equation*}
	and $F=\bigcap_{l\geq1}F_l$.
	Then $F_l$ and $F$ are non-empty closed sets.
	Let $E^{\xi}=\bigcap_{k\geq1}E^{\xi}_k$, then $E^{\xi}$ is non-empty and closed.
	It is easy to see that $F=\bigcup_{\xi\in\Omega^{\mathbb{N}_0}}E^{\xi}$.
	Note that if $x\in E^{\xi}$, then for any $i\in\mathbb{N}_0$, we have $f^{i(n_0+m_0)}x\in E^{\sigma^i\xi}$, where $\sigma:\Omega^{\mathbb{N}_0}\rightarrow\Omega^{\mathbb{N}_0}$ is the shift map defined by $(\sigma\xi)_k=\xi_{k+1}$.
	As a consequence, $F$ is a non-empty closed $f^{(n_0+m_0)}$-invariant set.
	Define
	\begin{equation*}
		Y:=\bigcup_{k=0}^{n_0+m_0-1}f^k(F).
	\end{equation*}
	Then $Y$ is a closed $f$-invariant set.

	Now we show that $h_{top}(f|_Y)>\gamma$.
	For $l\in\mathbb{N}$, let $s_l(\kappa,Y)$ denote the maximal cardinality of $(l,\kappa)$-separated sets in $Y$.
	For any $l\in\mathbb{N}$, we will show that
	\begin{equation*}
		s_{l(n_0+m_0)}(\varepsilon_0,Y)\geq\exp(ln_0(1+\tau+3\theta^*)\gamma).
	\end{equation*}
	For any $\xi=(x_0,\cdots,x_{l-1})\in\Omega^l$, we fix an $x_{\xi}\in E^{\tilde{\xi}}\subset F\subset Y$, where $\tilde{\xi}=(\tilde{x}_0,\tilde{x}_1,\cdots)\in\Omega^{\mathbb{N}_0}$ satisfies $\xi=(\tilde{x}_0,\cdots,\tilde{x}_{l-1})$; then $x_{\xi}\in E^{\tilde{\xi}}_l=E^{\xi}_l$.
	For $\xi,\zeta\in\Omega^l$ written as $\xi=(x_0,\cdots,x_{l-1})$ and $\zeta=(z_0,\cdots,z_{l-1})$, if $\xi\neq\zeta$, then there exists $0\leq j\leq l-1$ such that $x_j\neq z_j$.
	From the definition of $E^{\xi}_l$, it follows that for any $a(j)\leq i\leq b(j)$, we have
	\begin{equation}\label{eq-thm-intermediate-3-proof-7}
		\begin{split}
			d(f^ix_{\xi},f^ix_{\zeta})&\geq d(f^{i-a(j)}x_j,f^{i-a(j)}z_j)-\Big(d(f^{i-a(j)}x_j,f^ix_{\xi})+d(f^ix_{\zeta},f^{i-a(j)}z_j)\Big)\\
			&\geq d(f^{i-a(j)}x_j,f^{i-a(j)}z_j)-2\varepsilon_0.
		\end{split}
	\end{equation}
	Since $\Omega$ is $(n_0,\varepsilon^*)$-separated, there is an $i\in\{a(j),\cdots,b(j)\}$ such that
	\begin{equation}\label{eq-thm-intermediate-3-proof-8}
		d(f^{i-a(j)}x_j,f^{i-a(j)}z_j)>\varepsilon^*>5\varepsilon_0.
	\end{equation}
	Combining \eqref{eq-thm-intermediate-3-proof-7} and \eqref{eq-thm-intermediate-3-proof-8}, we obtain that
	\begin{equation*}
		d(f^ix_{\xi},f^ix_{\zeta})>3\varepsilon_0.
	\end{equation*}
	Hence $\{x_{\xi}:\xi\in\Omega^l\}$ is an $(l(m_0+n_0),\varepsilon_0)$-separated set in $Y$.
	By \eqref{eq-thm-intermediate-3-proof-6}, we obtain that
	\begin{equation}\label{eq-thm-intermediate-3-proof-9}
		s_{l(n_0+m_0)}(\varepsilon_0,Y)\geq|\Omega^l|\geq\exp(ln_0(1+\tau+3\theta^*)\gamma).
	\end{equation}
	Combining \eqref{eq-thm-intermediate-3-proof-4} and \eqref{eq-thm-intermediate-3-proof-9}, we conclude that
	\begin{equation*}
		h_{top}(f|_Y)\geq\limsup_{l\rightarrow\infty}\frac{1}{l(n_0+m_0)}\log s_{l(n_0+m_0)}(\varepsilon_0,Y)\geq\frac{n_0(1+\tau+3\theta^*)}{n_0+m_0}\gamma\geq\frac{1+\tau+3\theta^*}{1+\tau+2\theta^*}\gamma>\gamma.
	\end{equation*}

	We proceed to show that $h_{top}(f|_Y)<\gamma+\beta$.
	For any $l\in\mathbb{N}$, we claim that
	\begin{equation}\label{eq-thm-intermediate-3-proof-10}
		s_{l(n_0+m_0)}(3\varepsilon_0,F)\leq\Big(\exp(n_0(1+\tau+4\theta^*)\gamma)\cdot s_{m_0}(\varepsilon_0)\Big)^l.
	\end{equation}
	Let $\Gamma\subset F$ be an $(l(n_0+m_0),3\varepsilon_0)$-separated set with $|\Gamma|=s_{l(n_0+m_0)}(3\varepsilon_0,F)$.
	For every $\xi\in\Omega^l$, let $\Gamma(\xi)=\Gamma\cap E_l^{\xi}$.
	For each $0\leq j\leq l-1$, let $S_j$ be a maximal $(m_0,\varepsilon_0)$-separated set in $X$.
	By the maximality of $S_j$, for every $z\in X$, fix $\pi_j(z)\in S_j$ such that $d_{m_0}(z,\pi_j(z))\leq\varepsilon_0$.
	Define
	\begin{equation*}
		\Phi_{\xi}:\Gamma(\xi)\longrightarrow\prod_{j=0}^{l-1}S_j
	\end{equation*}
	by
	\begin{equation*}
		\Phi_{\xi}(x)=\Big(\pi_j\big(f^{j(n_0+m_0)+n_0}x\big)\Big)_{j=0}^{l-1}.
	\end{equation*}
	We claim that $\Phi_{\xi}$ is injective.
	Suppose that $\Phi_{\xi}(x)=\Phi_{\xi}(y)$ for some $x,y\in\Gamma(\xi)$.
	Then for every $0\leq j\leq l-1$, we have
	\begin{equation*}
		\begin{split}
			d_{m_0}\big(f^{j(n_0+m_0)+n_0}x,f^{j(n_0+m_0)+n_0}y\big)&\leq d_{m_0}\Big(f^{j(n_0+m_0)+n_0}x,\pi_j\big(f^{j(n_0+m_0)+n_0}x\big)\Big)\\
			&\quad+d_{m_0}\Big(\pi_j\big(f^{j(n_0+m_0)+n_0}y\big),f^{j(n_0+m_0)+n_0}y\Big)\leq2\varepsilon_0.
		\end{split}
	\end{equation*}
	By the construction of $E_l^{\xi}$, we also have
	\begin{equation*}
		d(f^ix,f^iy)\leq2\varepsilon_0,\quad\forall\,0\leq j\leq l-1,\,a(j)\leq i\leq b(j).
	\end{equation*}
	Therefore,
	\begin{equation*}
		d_{l(n_0+m_0)}(x,y)\leq2\varepsilon_0,
	\end{equation*}
	which implies that $x=y$, since $\Gamma$ is $(l(n_0+m_0),3\varepsilon_0)$-separated.
	Hence $\Phi_{\xi}$ is injective, and thus
	\begin{equation*}
		|\Gamma(\xi)|\leq\prod_{j=0}^{l-1}|S_j|\leq\big(s_{m_0}(\varepsilon_0)\big)^l.
	\end{equation*}
	Since $\Gamma\subset\bigcup_{\xi\in\Omega^l}\Gamma(\xi)$, by \eqref{eq-thm-intermediate-3-proof-6}, we have
	\begin{equation*}
		s_{l(n_0+m_0)}(3\varepsilon_0,F)=|\Gamma|\leq|\Omega^l|\cdot\big(s_{m_0}(\varepsilon_0)\big)^l\leq\Big(\exp(n_0(1+\tau+4\theta^*)\gamma)\cdot s_{m_0}(\varepsilon_0)\Big)^l.
	\end{equation*}
	Hence \eqref{eq-thm-intermediate-3-proof-10} holds.
	By the monotonicity of $s_n(3\varepsilon_0,F)$ with respect to $n$, and combining \eqref{eq-thm-intermediate-3-proof-2}, \eqref{eq-thm-intermediate-3-proof-4}, \eqref{eq-thm-intermediate-3-proof-5} and \eqref{eq-thm-intermediate-3-proof-10}, we obtain that
	\begin{equation*}
		\begin{split}
			\limsup_{n\rightarrow\infty}\frac{1}{n}\log s_n(3\varepsilon_0,F)&=\limsup_{l\rightarrow\infty}\frac{1}{l(n_0+m_0)}\log s_{l(n_0+m_0)}(3\varepsilon_0,F)\\
			&\leq\frac{n_0}{n_0+m_0}(1+\tau+4\theta^*)\gamma+\frac{m_0}{n_0+m_0}\cdot\frac{1}{m_0}\log s_{m_0}(\varepsilon_0)\\
			&\leq\frac{1+\tau+4\theta^*}{1+\tau-2\theta^*}\gamma+\frac{(\tau+2\theta^*)(h_{top}(f)+\theta^*)}{1+\tau+2\theta^*}<\gamma+\beta.
		\end{split}
	\end{equation*}
	Since $(X,f)$ has expansive constant larger than $\varepsilon^*>5\varepsilon_0$, by \eqref{eq-pre-entropy-UC-1-expansive}, we obtain that
	\begin{equation*}
		h_{top}^{UC}(f,F)=\limsup_{n\rightarrow\infty}\frac{1}{n}\log s_n(3\varepsilon_0,F)<\gamma+\beta.
	\end{equation*}
	According to Proposition \ref{pre-entropy-prop1} and \cite[Proposition 2]{Bowen-1973}, we conclude that
	\begin{equation*}
		h_{top}(f|_Y)=h_{top}^B(f,F)\leq h_{top}^{UC}(f,F)<\gamma+\beta.
	\end{equation*}
	This completes the proof of Theorem \ref{thm-intermediate-3}.
\end{proof}

\section{Ergodic optimization}\label{Sect-optimization}
\begin{proof}[{\bf Proof of Theorem \ref{thm-ergodic-optimization}}]
	First assume that $0<\tau<\infty$.
	Let $(X,f)$ be the subshift constructed in the proof of Theorem \ref{thm-Bowen-upper-1} with $m=1$, and let $\nu_1$ be the unique invariant measure of $Z=Z_1$.
	Then $(X,f)$ satisfies the non-uniform periodic $\tau$-specification property, and $\nu_1$ is its unique ergodic measure of maximal entropy with $0<h_{\nu_1}(f)=h_{top}(f)<\infty$.
	Set
	\begin{equation*}
		K=\overline{\mathcal{M}_f^{erg}(X)\setminus\{\nu_1\}}.
	\end{equation*}
	In this construction, $\mathcal{M}_f^{erg}(Z)=\{\nu_1\}$ and $K=\overline{\mathcal{M}_0}$. Since the cylinder $[0]$ is clopen, the function $\varphi_0=-\chi_{[0]}$ belongs to $C(X)$. Equations \eqref{eq-thm-construction-proof-1} and \eqref{eq-thm-construction-proof-2} give
	\begin{equation}\label{eq-thm-erg-opt-finite-gap}
		\int_X\varphi_0\mathrm{d}\nu_1=0
		\quad\text{and}\quad
		\sup_{\nu\in K}\int_X\varphi_0\mathrm{d}\nu
		\leq-\frac{\tau}{1+\tau}=-1+\frac{1}{1+\tau}.
	\end{equation}
	Define
	\begin{equation*}
		U=\big\{\varphi\in C(X):\|\varphi-\varphi_0\|<\delta/2\big\}.
	\end{equation*}
	Then $U$ is a non-empty open neighborhood of $\varphi_0$. For every $\varphi\in U$, \eqref{eq-thm-erg-opt-finite-gap} yields
	\begin{equation*}
		\sup_{\nu\in K}\int_X\varphi\mathrm{d}\nu
		\leq-1+\frac{1}{1+\tau}+\|\varphi-\varphi_0\|
		<-1+\frac{1}{1+\tau}+\frac{\delta}{2},
	\end{equation*}
	whereas
	\begin{equation*}
		\int_X\varphi\mathrm{d}\nu_1
		\geq-\|\varphi-\varphi_0\|>-\frac{\delta}{2}.
	\end{equation*}
	Since $\delta<1-\frac{1}{1+\tau}$, we also have $-1+\frac{1}{1+\tau}+\delta/2<-\delta/2$.

	Now assume that $\tau=\infty$. Choose
	\begin{equation*}
		0<\eta<\min\{\delta/4,1/2\}.
	\end{equation*}
	For $\tau=\infty$, take $m=1$ in the construction used in the proof of Theorem \ref{thm-Bowen-upper-1}, and choose $\sigma>0$ such that
	\begin{equation*}
		\frac{1}{1+\sigma}<\eta.
	\end{equation*}
	As in the case $0<\tau<\infty$, the subshift $(X,f)$ satisfies the non-uniform periodic $\infty$-specification property and has a unique ergodic measure $\nu_1$ of maximal entropy with $0<h_{\nu_1}(f)=h_{top}(f)<\infty$.
	With
	\begin{equation*}
		K=\overline{\mathcal{M}_f^{erg}(X)\setminus\{\nu_1\}}=\overline{\mathcal{M}_0}
	\end{equation*}
	and $\varphi_0=-\chi_{[0]}$, equation \eqref{eq-thm-construction-proof-infinity-4} gives
	\begin{equation}\label{eq-thm-erg-opt-infinite-gap}
		\int_X\varphi_0\mathrm{d}\nu_1=0
		\quad\text{and}\quad
		\sup_{\nu\in K}\int_X\varphi_0\mathrm{d}\nu
		\leq-\frac{\sigma}{1+\sigma}
		=-1+\frac{1}{1+\sigma}<-1+\eta.
	\end{equation}
	Define
	\begin{equation*}
		U=\big\{\varphi\in C(X):\|\varphi-\varphi_0\|<\delta/4\big\}.
	\end{equation*}
	For every $\varphi\in U$, \eqref{eq-thm-erg-opt-infinite-gap} implies
	\begin{equation*}
		\sup_{\nu\in K}\int_X\varphi\mathrm{d}\nu
		<-1+\eta+\frac{\delta}{4}
		<-1+\frac{\delta}{2}=-1+\frac{1}{1+\tau}+\frac{\delta}{2},
	\end{equation*}
	and
	\begin{equation*}
		\int_X\varphi\mathrm{d}\nu_1>-\frac{\delta}{4}>-\frac{\delta}{2}.
	\end{equation*}
	Since $\delta<1-\frac{1}{1+\tau}=1$ in this case, we have $-1+\delta/2<-\delta/2$.

	In either case, fix $\varphi\in U$ and put
	\begin{equation*}
		c_\varphi=\int_X\varphi\mathrm{d}\nu_1,
		\qquad
		b_\varphi=\sup_{\nu\in K}\int_X\varphi\mathrm{d}\nu.
	\end{equation*}
	The preceding estimates give $b_\varphi<c_\varphi$. Let $\mu\in\mathcal{M}_f(X)$ and let $\kappa_\mu$ be its ergodic decomposition. Set $p=\kappa_\mu(\{\nu_1\})$. If $\mu\neq\nu_1$, then $p<1$, and every ergodic component other than $\nu_1$ belongs to $\mathcal{M}_f^{erg}(X)\setminus\{\nu_1\}\subset K$. Consequently,
	\begin{equation*}
	\begin{split}
		\int_X\varphi\mathrm{d}\mu
		&=p c_\varphi+
		\int_{\mathcal{M}_f^{erg}(X)\setminus\{\nu_1\}}
		\bigg(\int_X\varphi\mathrm{d}\omega\bigg)\mathrm{d}\kappa_\mu(\omega)\\
		&\leq p c_\varphi+(1-p)b_\varphi<c_\varphi.
	\end{split}
	\end{equation*}
	Thus $\mathcal{M}_{\max}(\varphi)=\{\nu_1\}$ for every $\varphi\in U$.

	Finally, let
	\begin{equation*}
		\mathscr{S}=\big\{\varphi\in C(X):h_\mu(f)<h_{top}(f)
		\text{ for some }\mu\in\mathcal{M}_{\max}(\varphi)\big\}.
	\end{equation*}
	For every $\varphi\in U$, the unique maximizing measure is $\nu_1$ and $h_{\nu_1}(f)=h_{top}(f)$. Hence $U\cap\mathscr{S}=\varnothing$. Since $U$ is non-empty and open, $\mathscr{S}$ is not dense in $C(X)$.
\end{proof}

\bigskip

\textbf{Acknowledgements.} W. Lin is supported by the National Natural Science Foundation of China (No. 124B2010). X. Tian is supported by the National Natural Science Foundation of China (No. 12471182) and  Natural Science Foundation of Shanghai (No. 23ZR1405800).

\end{document}